\documentclass{amsart}
\usepackage{amscd}
\usepackage{amssymb}
\usepackage{graphicx}
\usepackage{color}

\numberwithin{equation}{section}
\numberwithin{figure}{section}

\theoremstyle{plain}
\newtheorem{theorem}{Theorem}[section]
\newtheorem{lemma}[theorem]{Lemma}
\newtheorem{corollary}[theorem]{Corollary}
\newtheorem{proposition}[theorem]{Proposition}
\newtheorem{conjecture}[theorem]{Conjecture}
\theoremstyle{definition}
\theoremstyle{remark}
\newtheorem{remark}[theorem]{Remark}
\newtheorem{question}[theorem]{Question}
\newtheorem*{acknowledgments}{Acknowledgments}

\newcommand\nc{\newcommand}

\nc\ord{\operatorname{ord}}
\nc\id{\operatorname{id}}
\nc\End{\operatorname{End}}
\nc\Span{\operatorname{Span}}
\nc\ev{{\operatorname{ev}}}
\nc\coev{{\operatorname{coev}}}
\nc\Gal{\operatorname{Gal}}

\nc\tmodp{{\tau ^{\operatorname{mod} p}}}

\nc\cl{\operatorname{\sf cl}}
\nc\Inv{{\operatorname{Inv}}}

\nc\FI[2]{\begin{figure}
    \begin{center}\input{#1.pstex_t}\end{center}
    \caption{#2}
    \label{#1}
  \end{figure}}
\nc\FIG[3]{\begin{figure}
    \includegraphics[#3]{#1.eps}
    \caption{#2}
    \label{fig:#1}
    \end{figure}}
\nc\FF[3]{\begin{figure}
    \includegraphics[#3]{#1.eps}
    \caption{#2}
    \label{#1}
    \end{figure}}

\nc\ho{{\hat\otimes }}
\nc\zzzcolon {\colon\thinspace}

\nc\modB {{\mathcal B}}
\nc\C{{\mathbb C}}
\nc\modF {{\mathcal F}}
\nc\F{{\mathbb F}}
\nc\tF{\tilde F}
\nc\g{{\mathfrak g}}
\nc\modh {\hbar}
\nc\modI {{\mathcal I}}
\nc\modi {{\mathbf i}}
\nc\tJ{{\tilde J}}
\nc\modK {{\mathcal K}}
\nc\modM {{\mathcal M}}
\nc\modN {{\mathbb N}}
\nc\modO {{\mathcal O}}
\nc\modP {{\mathcal P}}
\nc\modQ {{\mathbb Q}}
\nc\modr {{\mathbf r}}
\nc\modR {{\mathcal R}}
\nc\R{{\mathbb R}}
\nc\modS {{\mathcal S}}
\nc\mods {\operatorname{ev}}
\nc\modU {{\mathcal U}}
\nc\Uq{\modU _q}
\nc\hUq{\hat\modU _q}
\nc\Uqv{{\modU _q^\ev}}
\nc\hUqv{{\hat\modU _q^\ev}}
\nc\V{{\mathsf V}}
\nc\modv {{\mathbf v}}
\nc\modZ {{\mathbb Z}}
\nc\calZ{{\mathcal Z}}

\nc\ModUh{{\operatorname{\sf Mod}}_{U_h}}
\nc\bUh{\ul{U_h}}
\nc\bU{\bar U}
\nc\bUqv{{\bar U_q^\ev}}
\nc\bUqvt{{(\bUqv)\tilde{}}}
\nc\bUqvtn{\bUqvt\;^{\tOn}}
\nc\bUqvtx[1]{\bUqvt\;^{\tO #1}}
\nc\bD{{\underline\Delta }}
\nc\bS{{\underline S}}
\nc\hf{{\frac12}}
\nc\tU{\tilde \modU }
\nc\tUq{\tU_q}
\nc\tUqv{\tUq^\ev}
\nc\tUqvn{(\tUqv)^{\tOn}}
\nc\tUqvx[1]{(\tUqv)^{\tO#1}}
\nc\Uqvx[1]{(\Uqv)^{\otimes #1}}
\nc\tUqx[1]{\tUq^{\tO#1}}
\nc\Uqvn{(\Uqv)^{\otimes n}}
\nc\hP{{\hat\modP }}
\nc\ul{\underline}
\nc\trr{\triangleright}
\nc\trrn{\trr_n}
\nc\BBq[2]{\{#1\}_{q,#2}}
\nc\bb[2]{\biggl[\begin{matrix}{#1}\\{#2}\end{matrix}\biggr]}
\nc\bbq[2]{{\bb{#1}{#2}}_q}
\nc\BB[2]{\{#1\}_{#2}}
\nc\zzzvert {\ |\ }
\nc\Uhho[1]{{U_h^{\ho #1}}}
\nc\Uhhon{{\Uhho n}}
\nc\Uhn{{\Uhho n}}
\nc\Uhx[1]{{\Uhho {#1}}}
\nc\tO{{\tilde\otimes }}
\nc\tOn{{\tilde\otimes n}}
\nc\LMO{{\operatorname{LMO}}}
\nc\BT{\operatorname{\sf BT}}
\nc\MMR{{\operatorname{MMR}}}
\nc\Zqq{\modZ [q,q^{-1}]}
\nc\Zvv{\modZ [v,v^{-1}]}
\nc\Zqh{\widehat{\modZ [q]}}
\nc\Qqh{\widehat{\modQ [q]}}
\nc\Za{\widehat{\modZ [\alpha ]}^\alpha }
\nc\Zx[1]{\widehat{\modZ [#1]}^{#1}}
\nc\vn[1]{\modv ^n_{#1}}
\nc\trqV{{\trq^V}}
\nc\trq{\tr_q}
\nc\lk{\operatorname{lk}}
\nc\tr{{\operatorname{tr}}} 
\nc\tP{\tilde P}
\nc\congto{\overset{\cong}{\rightarrow }}
\nc\floor[1]{\lfloor#1\rfloor}
\nc\zq{{\zeta ^{1/4}}}
\nc\Funz{\Fun^0(\modN ,\modZ _+)}
\nc\Fung{\Fun_g(\bFp^\times ,\bFp)}
\nc\Fun{\operatorname{Fun}}
\nc\Fp{\F_p}
\nc\Fpq{\Fp[q]}
\nc\Fpqh{\widehat{\F_p[q]}}
\nc\bFp{\bar\F_p}
\nc\half{{\frac12}}

\begin{document}

\title[A unified Witten-Reshetikhin-Turaev invariant]{A unified
 Witten-Reshetikhin-Turaev invariant\\for integral homology spheres}

\author{Kazuo Habiro}

\address{Research Institute for Mathematical Sciences\\ Kyoto
University\\ Kyoto\\ 606-8502\\ Japan}

\email{habiro@kurims.kyoto-u.ac.jp}

\thanks{This research was partially supported by the Japan Society for
  the Promotion of Science, Grant-in-Aid for Young Scientists (B),
  16740033.}

\date{May 12, 2006}

\begin{abstract}
  We construct an invariant $J_M$ of integral homology spheres $M$
  with values in a completion $\Zqh$ of the polynomial ring $\modZ [q]$
  such that the evaluation at each root of unity $\zeta $ gives the the
  $SU(2)$ Witten-Reshetikhin-Turaev invariant $\tau _\zeta (M)$ of $M$ at
  $\zeta $.  Thus $J_M$ unifies all the $SU(2)$ Witten-Reshetikhin-Turaev
  invariants of $M$.  As a consequence, $\tau _\zeta (M)$ is an algebraic
  integer.  Moreover, it follows that $\tau _\zeta (M)$ as a function on $\zeta $
  behaves like an ``analytic function'' defined on the set of roots of
  unity.  That is, the $\tau _\zeta (M)$ for all roots of unity are
  determined by a ``Taylor expansion'' at any root of unity, and also
  by the values at infinitely many roots of unity of prime power
  orders.  In particular, $\tau _\zeta (M)$ for all roots of unity are
  determined by the Ohtsuki series, which can be regarded as the
  Taylor expansion at $q=1$.
\end{abstract}

\maketitle

\section{Introduction}
\label{sec1}

In this paper we construct an invariant of integral homology spheres
which unifies the $SU(2)$ Witten-Reshetikhin-Turaev invariants at all
roots of unity, which we announced in a previous paper
\cite{H:rims2001}.

\subsection{The WRT invariant for integral homology spheres}
\label{sec:witt-resh-tura-1}

Witten \cite{Witten} introduced the notion of Chern-Simons path
integral which gives a quantum field theory interpretation of the
Jones polynomial \cite{Jones1,Jones2} and predicts the existence of
$3$-manifold invariants.  Using the quantum group $U_q(sl_2)$ at roots
of unity, Reshetikhin and Turaev \cite{Reshetikhin-Turaev:invariants}
gave a rigorous construction of $3$-manifold invariants, which are
believed to coincide with the Chern-Simons path integrals.  These
invariants are called the Witten-Reshetikhin-Turaev (WRT) invariants.

In the present paper, we focus on the WRT invariant for integral
homology spheres, i.e., closed $3$-manifolds $M$ with trivial first
homology groups.  If we fix such $M$, the WRT invariant $\tau _\zeta (M)\in \C$
is defined for each root of unity $\zeta $.  (Unlike the case of general
closed $3$-manifolds, one does not have to specify a fourth root of
$\zeta $.)  For an integer $r\ge 1$, $\tau _{\zeta _r}(M)$ for
$\zeta _r=\exp\frac{2\pi \sqrt{-1}}{r}$ is also denoted by $\tau _r(M)$.  In
the literature, usually $\tau _1(M)$ is not defined, but for our purpose
it is convenient to defined it as $1$.  For integral homology spheres,
the version $\tau '_r(M)$ introduced by Kirby and Melvin
\cite{Kirby-Melvin} defined for odd $r\ge 3$ is equal to $\tau _r(M)$.

Let $\calZ\subset \C$ denote the set
of all roots of unity.  Define the {\em WRT function} of $M$
\begin{equation*}
  \tau (M)\zzzcolon \calZ\rightarrow \C,
\end{equation*}
by $\tau (M)(\zeta )=\tau _\zeta (M)$.  The behavior of the function $\tau (M)$ is of
interest here.

For $\tau _\zeta (M)$ with $\zeta $ roots of unity of {\em odd prime (power)
orders}, there have been more extensive studies than for the other
cases.  Murakami \cite{Murakami1} proved that if $\zeta \in \calZ$ is of odd
prime order, then $\tau _\zeta (M)\in \modZ [\zeta ]$, hence it is an algebraic
integer.  Ohtsuki \cite{Ohtsuki1} extracted from the $\tau _\zeta (M)$ for
$\zeta $ of odd prime orders a power series invariant, known as the {\em
Ohtsuki series} of $M$
\begin{equation*}
  \tau ^O(M)=1+\sum_{n=1}^\infty \lambda _n(M)(q-1)^n\in \modQ [[q-1]].
\end{equation*}
Lawrence \cite{Lawrence:integrality} conjectured, and Rozansky
\cite{Rozansky} later proved, that $\tau ^O(M)\in \modZ [[q-1]]$ and that, for
each $\zeta \in \calZ$ of odd prime power order $p^e$, $\tau ^O(M)|_{q=\zeta }$
converges $p$-adically to $\tau _\zeta (M)$.  In this sense, the Ohtsuki
series {\em unifies} the WRT invariants at roots of unity of odd prime
power orders.  (For generalizations of the above-mentioned results to
rational homology spheres and to invariants associated to other Lie
groups, see
\cite{Le1,Le2,Masbaum-Roberts,Masbaum-Wenzl,Murakami2,Ohtsuki2,Takata-Yokota}.)

The proofs of the above-mentioned results depends heavily on the fact
that if $\zeta \in \calZ$ is of prime power order, then $\zeta -1$ is not a unit
in the ring $\modZ [\zeta ]$.  Otherwise, $\zeta -1$ is a unit in $\modZ [\zeta ]$, and
expansions in powers of $\zeta -1$ do not work.

\subsection{The ring $\Zqh$ of analytic functions on the
  set of roots of unity}
\label{sec:ring-zqh-cyclotomic}

The invariant $J_M$ of an integral homology sphere $M$ which we
construct in the paper takes values in a completion $\Zqh$ of the
polynomial ring $\modZ [q]$, which was introduced in \cite{H:rims2001} and
studied in \cite{H:cyclotomic}.  One of the simplest definitions of
$\Zqh$ is
\begin{equation*}
  \Zqh = \varprojlim_{n}\modZ [q]/((q)_n),
\end{equation*}
where we set
\begin{equation*}
  (q)_n = (1-q)(1-q^2)\cdots(1-q^n).
\end{equation*}

The ring $\Zqh$ may be regarded as the ring of ``analytic functions
defined on the set $\calZ$ of roots of unity''.  This statement is
justified by the following facts.  (The following overlaps those in
\cite{H:rims2001,H:cyclotomic}.)

First of all, {\em each element of $\Zqh$ can be evaluated at each
root of unity}.  That is, for each $\zeta \in \calZ$, the evaluation map
$\mods _\zeta \zzzcolon \modZ [q]\rightarrow \modZ [\zeta ]$, $f(q)\mapsto f(\zeta )$, induces a (surjective)
ring homomorphism
\begin{equation*}
  \mods _\zeta \zzzcolon \Zqh\rightarrow \modZ [\zeta ],
\end{equation*}
since $\mods _\zeta ((q)_n)=0$ if $n\ge \ord(\zeta )$.  It is often useful to write
$f(\zeta )=\mods _\zeta (f(q))$ for $f(q)\in \Zqh$.

Second, {\em each element of $\Zqh$ can be regarded as a
 (set-theoretic) function on the set of roots of unity.}  This means
 that each $f\in \Zqh$ is determined uniquely by the values
 $\mods _\zeta (f)\in \modZ [\zeta ]$ for all $\zeta \in \calZ$, or, equivalently, the
 function
\begin{equation}
  \label{e4}
  \mods _{\calZ}\zzzcolon \Zqh\rightarrow \prod_{\zeta \in \calZ}\modZ [\zeta ],\quad
  f(q)\mapsto (f(\zeta ))_{\zeta \in \calZ}.
\end{equation}
is injective \cite[Theorem 6.3]{H:cyclotomic}.  Here, in a natural
way, $\prod_{\zeta \in \calZ}\modZ [\zeta ]$ can be regarded as a subring of the
ring of $\C$-valued functions on $\calZ$.

Third, {\em each element of $\Zqh$ has a power series expansion in
  $q-\zeta $ for each root of unity, and any such power series determines
  $\Zqh$.}  In fact, for each $\zeta \in \calZ$, the inclusion
$\modZ [q]\subset \modZ [\zeta ][q]$ induces a ring homomorphism
\begin{equation}
  \label{e6}
  \iota _\zeta \zzzcolon \Zqh\rightarrow \modZ [\zeta ][[q-\zeta ]],
\end{equation}
since, for each $i\ge 0$, $(q)_n$ is divisible by $(q-\zeta )^i$ in
$\modZ [\zeta ][q]$ if $n\ge i\ord(\zeta )$.  Since $\iota _\zeta $ is injective
\cite[Theorem 5.2]{H:cyclotomic}, each $f\in \Zqh$ is determined by the
power series $\iota _\zeta (f)$.  $\iota _\zeta (f)$ may be regarded as the Taylor
expansion of $f$, see Section \ref{sec:power-seri-expans-1}.

Fourth, {\em each element of $\Zqh$ is completely determined by its
values on a subset $\calZ'\subset \calZ$ if $\calZ'$ has a limit point.
Otherwise, not completely.}  To explain what this means, we introduce
a topology on the set $\calZ$, which is different from the usual one
induced by the topology of $\C$.  Two elements $\zeta ,\xi \in \calZ$ are said
to be {\em adjacent} if $\zeta \xi ^{-1}$ is of prime power order, or,
equivalently, if $\zeta -\xi $ is not a unit in $\modZ [\zeta ,\xi ]$.  A subset
$\calZ'\subset \calZ$ is defined to be {\em open} if, for each $\zeta \in \calZ'$,
all but finitely many elements adjacent to $\zeta $ is contained in
$\calZ'$.  In this topology, an element
$\xi \in \calZ$ is a limit point of a subset $\calZ'\subset \calZ$ (i.e.,
$(U\setminus \{\xi \})\cap \calZ'\neq\emptyset$ for all neighborhood $U$ of $\xi $)
if and only if there are infinitely many $\zeta \in \calZ'$ adjacent to
$\xi $.  We have the following.

\begin{proposition}[{\cite[Theorem 6.3]{H:cyclotomic}}]
  \label{r85}
  If $\calZ'\subset \calZ$ has a limit point, then the ring homomorphism
  \begin{equation}
    \label{e61}
    \mods _{\calZ'}\zzzcolon \Zqh\rightarrow \prod_{\zeta \in \calZ'}\modZ [\zeta ],\quad
    f(q)\mapsto (f(\zeta ))_{\zeta \in \calZ},
  \end{equation}
  is injective.
\end{proposition}

If $\calZ'\subset \calZ$ has no limit point, then $\mods _{\calZ'}$ is not
injective, i.e., there is a non-zero ``analytic function'' $f\in \Zqh$
vanishing on $\calZ'$, see Proposition \ref{r56}.

The above-explained properties of $\Zqh$ are closely related to the
integrality of $\Zqh$.  In fact, the completion $\Qqh =
\varprojlim_{n}\modQ [q]/((q)_n)$ does not behave like $\Zqh$, see
\cite[Section 7.5]{H:cyclotomic}.

\subsection{A unified WRT invariant $J_M$ with values in $\Zqh$}
\label{sec:unif-wrt-invar-1}

The following is the main result of the present paper, and follows
from Theorems \ref{thm:41} and \ref{thm:65}.

\begin{theorem}
  \label{r40}
  There is an invariant $J_M\in \Zqh$ of an integral homology sphere $M$
  such that for any root of unity $\zeta $ we have
  \begin{equation}
    \label{e56}
    \mods _\zeta (J_M) = \tau _\zeta (M).
  \end{equation}
\end{theorem}

The properties of the ring $\Zqh$ explained in the last subsection
implies that the WRT function $\tau (M)$ may be regarded as an analytic
function defined on $\calZ$.  Let us describe some corollaries to
Theorem \ref{r40} and properties of the ring $\Zqh$.  (Part of the
discussion below overlaps those in \cite{H:rims2001,H:cyclotomic}.)

An immediate consequence of Theorem \ref{r40} is the following
generalization of Murakami's integrality result.

\begin{corollary}[Conjectured by Lawrence \cite{Lawrence:integrality}]
  \label{r14}
  For any integral homology sphere $M$ and for $\zeta \in \calZ$, we have
  $\tau _\zeta (M)\in \modZ [\zeta ]$.
\end{corollary}

Theorem \ref{r40} immediately implies {\em Galois equivariance} of
$\tau _\zeta (M)$.  Namely, we have
\begin{equation}
  \label{e67}
  \tau _{\alpha (\zeta )}(M)=\alpha (\tau _\zeta (M))
\end{equation}
for $\alpha \in \Gal(\modQ ^{ab}/\modQ )$ and $\zeta \in \calZ$.  Here $\modQ ^{ab}$ denotes
the maximal abelian extension of $\modQ $, which is generated over $\modQ $ by
all roots of unity.  ($\Gal(\modQ ^{ab}/\modQ )$ can be identified with the
automorphism group of the group $\calZ\cong\modQ /\modZ $.)  It is well known
that the Galois equivariance of $\tau _\zeta (M)$ is implied by Reshetikhin
and Turaev's definition.  Theorem \ref{r40} reexplains it in an
apparent way.

Proposition \ref{r85} implies the following.

\begin{theorem}
  \label{r48}
  The invariant $J_M$ is determined by the WRT function $\tau (M)$.
  (Thus $J_M$ and $\tau (M)$ have the same strength in distinguishing two
  integral homology spheres.)  Moreover, both $J_M$ and $\tau (M)$ are
  determined by the values of $\tau _\zeta (M)$ for $\zeta \in \calZ'$, where
  $\calZ'\subset \calZ$ is any infinite subset with a limit point in the
  sense explained in Section \ref{sec:ring-zqh-cyclotomic}.
\end{theorem}

The invariant $J_M$ unifies not only the WRT invariants but also the
Ohtsuki series $\tau ^O(M)$.  Namely, we have (Theorem \ref{r52})
\begin{equation}
  \label{e12}
  \iota _1(J_M)=\tau ^O(M).
\end{equation}

Injectivity of $\iota _\zeta $ in \eqref{e6} for each $\zeta \in \calZ$ implies that
the $J_M$ and hence $\tau (M)$ are determined by the power series
expansion $\iota _\zeta (J_M)\in \modZ [\zeta ][[q-\zeta ]]$.  In particular, $J_M$ is
determined by $\tau ^O(M)$, in view of \ref{e12}.  Thus $J_M$ and
$\tau ^O(M)$ have the same strength in distinguishing integral homology
spheres.  As a consequence, the Le-Murakami-Ohtsuki invariant
\cite{LMO} determines $J_M$ and $\tau (M)$, since it determines $\tau ^O(M)$
(see \cite{Ohtsuki4}).

For further properties of $J_M$, see Sections \ref{sec12} and
\ref{sec13}.

\subsection{Organization of the paper}
\label{sec:organization-paper}

In Section \ref{sec2}, we first recall the definition of the quantized
enveloping algebra $U_h=U_h(sl_2)$ of the Lie algebra $sl_2$, which is
a $h$-adically complete Hopf $\modQ [[h]]$-algebra, and then introduce
$\modZ [q,q^{-1}]$-subalgebras $\Uq$, where $\modZ [q,q^{-1}]$ is regarded as
a subring of $\modQ [[h]]$ by setting $q=\exp h$.  The Hopf algebra
structure of $U_h$ induces a Hopf algebra structure on $\Uq$.  We also
introduce a $\modZ [q,q^{-1}]$-subalgebra $\Uqv$, which is the even part
of $\Uq$ with respect to a natural $(\modZ /2\modZ )$-grading of $\Uq$.
We define completions $\tUq$ and $\tUqv$ of the algebras $\Uq$ and
$\Uqv$, and also completed tensor products of copies of $\tUq$ and
$\tUqv$.  $\tUq$ is equipped with a complete Hopf algebra structure
induced by the Hopf algebra structure of $\Uq$.

In Section \ref{sec3}, we first recall the ribbon Hopf algebra
structure for $U_h$, and then the braided Hopf algebra structure
$\bUh$, canonically defined for $U_h$.  The main observation in
this section (Theorem \ref{thm:26}) is that $\tUqv$ is equipped with a
braided Hopf algebra structure inherited from that for
$\bUh$.

In Section \ref{sec4}, we first recall from \cite{H:universal} the
notion of bottom tangles.  An $n$-component bottom tangle $T$ is a
tangle in a cube consisting of $n$ arc components whose endpoints lie
in a line in the bottom square of the cube in such a way that between
the two endpoints of each arc there are no endpoints of other arcs.
Then we adapt the universal invariants for bottom tangles associated
to the ribbon Hopf algebra to the case of $U_h$.  The universal
invariant $J_T$ of $T$ takes values in the invariant part
$\Inv(\Uhhon)$ of the $n$-fold completed tensor product $\Uhhon$ of
$U_h$, where the invariant part is considered with respect to the
standard tensor product left $U_h$-module structure defined using the
left adjoint action of $U_h$.  The main observation in this section
(Theorem \ref{thm:8}) is that if $T$ is an $n$-component,
algebraically-split, $0$-framed bottom tangle, then $J_T$ is contained
in the invariant part $\modK _n$ of the $n$-fold completed tensor product
$\tUqvn$ of $\tUqv$.  (Recall that a link $L$ is {\em algebraically
split} if the linking number of any two distinct components are zero.)
The proof of Theorem \ref{thm:8} uses the braided Hopf algebra
structure of $\tUqv$.

For a $0$-framed bottom knot (i.e., $1$-component bottom tangle) $T$,
the universal invariant $J_T$ of $T$ takes values in the center
$Z(\tUqv)$ of $\tUqv$.  In view of a result proved in \cite{H:uqsl2},
this implies that we can express $J_T$ as an infinite sum
$\sum_{p\ge 0}a_p(T)\sigma _p$, where $a_p(T)\in \Zqq$, $p\ge 0$, and where
$\sigma _p\in Z(\Uqv)$, $p\ge 0$, are a certain basis of $Z(\Uqv)\subset Z(\tUqv)$.
Each $a_p(T)$ gives a $\Zqq$-valued invariant of $T$, which may be
regarded also as an invariant of the closure $\cl(T)$ of $T$.

In Section \ref{sec5}, we recall the definition of the colored Jones
polynomial $J_L(W_1,\ldots,W_n)$ of a framed link $L=L_1\cup \cdots\cup L_n$, where
each component $L_i$ of $L$ is colored by a finite-dimensional,
irreducible representation $W_i$ of $U_h(sl_2)$.  Then we recall from
\cite{H:universal} a formulation of the colored Jones polynomial using
universal invariant of bottom tangle:
\begin{equation*}
  J_L(W_1,\ldots,W_n)=(\trq^{W_1}\otimes \cdots\otimes \trq^{W_n})(J_T),
\end{equation*}
where $T$ is a bottom tangle whose closure is $L$, and where
$\trq^{W_i}\zzzcolon U_h\rightarrow \modQ [[h]]$ is the quantum trace in $W_i$.

Recall that for each $d\ge 0$ there is exactly one $(n+1)$-dimensional,
irreducible representation of $U_h$ up to isomorphism, denoted by
$\V_d$.  We need extensions of the colored Jones polynomials for
framed links whose components are colored by linear combinations (over
$\modQ (v)$, $\modZ [v,v^{-1}]$, etc., where $v=q^{1/2}=\exp\frac h2$) of the
$\V_d$, which are defined naturally by multilinearity.

In Section \ref{sec6}, we relate the universal invariant
$J_T\in Z(\tUqv)$ and the $\Zqq$-valued invariants $a_p(T)$, $p\ge 0$, of
a bottom knot $T$ defined in Section \ref{sec4} to the colored Jones
polynomials of the closure $\cl(T)$ of $T$.  Theorem \ref{r18}
identifies $a_p(T)$ with $J_{\cl(T)}(P''_p)$, where $P''_p$ is a
$\modQ (v)$-linear combination of $\V_0,\V_1,\ldots,\V_p$.

In Section \ref{sec7}, we give some remarks on the universal invariant
$J_T\in Z(\tUqv)$ for a bottom knot.  First of all, $Z(\tUqv)$ is
identified with a completion $\Lambda $ of $\modZ [q,q^{-1},t+t^{-1}]$.  For a
knot $K=\cl(T)$ with $T$ a bottom knot, we set $J_K(t,q)=J_T\in \Lambda $ by
abuse of notation, which we call the {\em two-variable colored Jones
invariant} of $K$.  The normalized colored Jones polynomial
$J_K(\V_n)/J_{\text{unknot}}(\V_n)\in \modZ [q,q^{-1}]$ is equal to the
specialization $J_K(q^{n+1},q)$ of $J_K(t,q)$.  The specialization
$J_K(1,q)\in \Zqh$ can be regarded as a universal form for the Kashaev
invariants of $K$.  We relate the invariant $J_K(t,q)$ to Rozansky's
integral version of the Melvin-Morton expansion of the colored Jones
polynomials of $K$, and give several conjectures which generalizes
Rozansky's rationality theorem.

In Section \ref{sec8}, we consider invariants of algebraically-split
links.  We define a $\modZ [q,q^{-1}]$-algebra $\modP $ which is spanned by
certain normalizations $\tilde P'_n$ of $P''_n$, and define a
completion $\hP$ of $\modP $.  We show that for any $n$-component,
algebraically-split, $0$-framed link $L$ and for any elements
$x_1,\ldots,x_n\in \hP$, there is a well-defined element
$J_L(x_1,\ldots,x_m)\in \Zqh$, see Corollary \ref{r87}.  In the proof,
results proved in the previous sections, such as Theorem \ref{thm:8},
are used.

In Section \ref{sec9}, we define an element $\omega $ in the ring $\hP$.
Theorem \ref{thm:99} states that for any $(m+1)$-component,
algebraically-split, $0$-framed link $L_1\cup \cdots\cup L_m\cup K$ such that $K$
is an unknot, and for $x_1,\ldots,x_m\in \hP$, we have
\begin{equation*}
  J_{L\cup K}(x_1,\ldots,x_m,\omega ^{\mp1})=J_{L_{(K,\pm 1)}}(x_1,\ldots,x_m).
\end{equation*}
Here $L_{(K,\pm 1)}$ is the framed link in $S^3$ obtained from
$L_1\cup \cdots\cup L_m$ by $\pm 1$-framed surgery along $K$.

In Section \ref{sec10}, we prove the existence of an invariant
$J_M\in \Zqh$ of integral homology sphere $M$ (Theorem \ref{thm:41}).
We outline the proof below.  Recall that $M$ can be expressed as the
result of surgery along an algebraically-split framed link
$L=L_1\cup \cdots\cup L_m$ in $S^3$ with framings $f_1,\ldots,f_m\in \{\pm 1\}$.  Then
$J_M$ is defined by
\begin{equation*}
  J_M:= J_{L^0}(\omega ^{-f_1},\omega ^{-f_2},\ldots,\omega ^{-f_m}),
\end{equation*}
where $L^0$ is the framed link obtained from $L$ by changing all the
framings to $0$.  By Corollary \ref{r87}, we have $J_M\in \Zqh$.  To
prove that $J_M$ does not depend on the choice of $L$, we use the
twisting property of $\omega $ (Theorem \ref{thm:99}) and a refined version
of Kirby's calculus for algebraically-split, $\pm 1$-framed links (see
Theorem \ref{r32}), which was conjectured by Hoste \cite{Hoste} and
proved in \cite{H:kirby1}.  The proof of Theorem \ref{thm:41} does not
involve any existence proof of the WRT
invariants $\tau _\zeta (M)$ at roots of unity $\zeta $, hence can be regarded as
a new, unified proof for the existence for $\tau _\zeta (M)$, after
establishing the specialization property \eqref{e56}.

In Section \ref{sec11}, we prove this specialization property (Theorem
 \ref{thm:65}).  We also give an alternative proof of the existence of
 $J_M$ which uses the existence of $\tau _\zeta (M)$ but does not use Theorem
 \ref{r32}.

In Section \ref{sec12}, we make several observations and give some
applications.  In Section~\ref{sec:conn-sum-orient}, we observe the
behavior of $J_M$ under taking connected sums and
orientation-reversal.  In Section \ref{sec:determination-jm-wrt}, we
observe the failure of an approach to the conjecture that the
WRT invariants $\tau _\zeta (M)$ at any infinitely many
roots of unity determine $J_M$.  In Section
\ref{sec:power-seri-expans}, we study the power series invariants
$\iota _\zeta (J_M)$, including the Ohtsuki series $\tau ^O(M)=\iota _1(J_M)$.  In
Section \ref{sec:divisibility-jm-1}, we give some divisibility results
for $J_M-1$, etc., implied by well-known results for $\tau _\zeta (M)$ for
$\zeta $ of small orders $1,2,3,4,6$, and give some applications of these
results to the coefficients of the Ohtsuki series and the power series
$\iota _{-1}(J_M)$.  We also state a conjecture about the values of the
eighth WRT invariant $\tau _8(M)$.

In Section \ref{sec13}, we first observe that for any complex number
$\alpha $, there is a formal specialization of $J_M$ at $q=\alpha $.  Motivated
by this observation, for each prime $p$, we define a $p$-adic analytic
version $\tau ^p(M)$ of the WRT function, which is
a $p$-adic analytic function from the unit circle in the field $\C_p$
of complex $p$-adic numbers into the valuation ring of $\C_p$.  The
mod $p$ reduction of $\tau ^p(M)$, denoted by $\tmodp(M)$, is defined
on the group of units, $\bar\F_p^\times $, in the algebraic closure
$\bar\F_p$ of the field $\F_p$ of $p$ elements, and takes values in
$\bar\F_p$.

In Section \ref{sec14}, we compute some examples of $J_M$ for integral
homology spheres obtained as the result of surgery along the Borromean
rings $A$ with framings $1/a,1/b,1/c$ with $a,b,c\in \modZ $, and some
related knot and link invariants.  First we compute the colored Jones
polynomials of the Borromean rings.  Then we compute the powers of the
ribbon element in $U_h$ and the powers of the twist element~$\omega $.
This enables us to compute the invariants of the result of surgery
from the Borromean rings by performing surgery along some (possibly
all) of the three components by framings in $\{1/m\zzzvert m\in \modZ \}$.

In Section \ref{sec15}, we first generalize the universal invariant
$J_K$ of a knot in $S^3$ to knots in integral homology spheres
(Theorem \ref{r79}).  Using this invariant, we prove that if two
integral homology spheres $M$ and $M'$ are related by surgery along a
knot with framing $1/m$ with $m\in \modZ $, then $J_M$ and $J_{M'}$ are
congruent modulo $q^{2m}-1$ (Theorem \ref{r9}).  This result suggests
that it would be natural to consider a generalization of Ohtsuki's
theory of finite type invariants of integral homology spheres,
involving $(1/m)$-surgeries ($m\in \modZ $) along knots.

In Section \ref{sec16}, we give some remarks.  In Section
\ref{sec:wrt-invariants-as}, we discuss the relationships between the
unified WRT invariant $J_M$ and another approach to unify the WRT
invariants by realizing them as limiting values of holomorphic
functions on the disk $|q|<1$.  In Sections \ref{sec:general} and
\ref{sec:rati-homol-spher}, we mention generalizations of $J_M$ to
simple Lie algebras and rational homology spheres.  In Section
\ref{sec:homology-cylinders}, we announce a generalization of $J_M$ to
certain cobordisms of surfaces, which includes homology cylinders.

\section{The algebra $U_h(sl_2)$ and its subalgebras}
\label{sec2}

In this section, we recall the definition and some properties of the
quantized enveloping algebra $U_h(sl_2)$.  Then we define subalgebras
$\Uq$ and $\Uqv$ of $U_h$, as well as their completions $\tUq$ and
$\tUqv$, which we studied in \cite{H:uqsl2}.

\subsection{$q$-integers}
\label{sec:q-integers}
Let $h$ be an indeterminate, and set
\begin{gather*}
  v=\exp\frac h2\in \modQ [[h]],\quad
  q=v^2=\exp h\in \modQ [[h]].
\end{gather*}
We have $\Zqq\subset \Zvv\subset \modQ [[h]]$.

We use two systems of $q$-integer notations.  One is the
``$q$-version'':
\begin{gather*}
  \{i\}_q = q^i-1,\quad
  \BBq{i}{n} = \{i\}_q\{i-1\}_q\cdots\{i-n+1\}_q,\quad
  \{n\}_q! =\BBq{n}{n},\\
  [i]_q=\{i\}_q/\{1\}_q,\quad
  [n]_q! = [n]_q[n-1]_q\cdots[1]_q,\quad
  \bbq{i}{n} = \BBq{i}{n}/\{n\}_q!,
\end{gather*}
for $i\in \modZ $, $n\ge 0$.  These are elements in $\Zqq$.  (In later
sections, we also use $(q)_n=(-1)^n\{n\}_q!$.)  The other is the
``balanced $v$-version'':
\begin{gather*}
  \{i\} = v^i-v^{-i},\quad
  \BB{i}{n} = \{i\}\{i-1\}\cdots\{i-n+1\},\quad
  \{n\}!=\BB{n}{n},\\
  [i]=\{i\}/\{1\},\quad
  [n]! = [n][n-1]\cdots[1],\quad
  \bb{i}{n} = \BB{i}{n}/\{n\}!,
\end{gather*}
for $i\in \modZ $, $n\ge 0$.  These are elements of $\Zqq\sqcup v\Zqq\subset \Zvv$.
These two families of notations are the same up to multiplication by
powers of $v$.  The former system is useful in clarifying that
formulas are defined over $\modZ [q,q^{-1}]$.  The latter is sometimes
useful in clarifying that formulas have symmetry under conjugation
$v\leftrightarrow v^{-1}$.

\subsection{The quantized enveloping algebra $U_h$}
\label{sec:quantized-enveloping-algebra}
We define $U_h=U_h(sl_2)$ as the $h$-adically complete
$\modQ [[h]]$-algebra, topologically generated by the elements $H$, $E$,
and $F$, satisfying the relations
\begin{equation*}
  HE-EH=2E,\quad HF-FH=-2F,\quad
  EF-FE=
  \frac{K-K^{-1}}{v-v^{-1}},
\end{equation*}
where we set
\begin{gather*}
  K= v^H = \exp\frac{hH}2.
\end{gather*}

The algebra $U_h$ has a complete Hopf algebra structure with the
comultiplication $\Delta \zzzcolon U_h\rightarrow U_h\ho U_h$, the counit $\epsilon \zzzcolon U_h\rightarrow \modQ [[h]]$
and the antipode $S\zzzcolon U_h\rightarrow U_h$ defined by
\begin{gather*}
    \Delta (H)=H\otimes 1+1\otimes H,\quad\epsilon (H)=0,\quad S(H)=-H,\\
    \Delta (E)=E\otimes 1+K\otimes E,\quad\epsilon (E)=0,\quad S(E)=-K^{-1}E,\\
    \Delta (F)=F\otimes K^{-1}+1\otimes F,\quad\epsilon (F)=0,\quad S(F)=-FK.
\end{gather*}
(Here $\ho$ denotes the $h$-adically completed tensor product.)

For $p\in \modZ $, let $\Gamma _p(U_h)$ denote the complete $\modQ [[h]]$-submodule
of $U_h$ topologically spanned by the elements $F^iH^jE^k$ with
$i,j,k\ge 0$, $k-i=p$.  This gives a topological $\modZ $-graded algebra
structure for $U_h$
\begin{equation*}
  U_h = \hat\oplus_{p\in \modZ }\Gamma _p(U_h).
\end{equation*}
(Here $\hat\oplus$ denotes $h$-adically completed direct sum.)  The
elements of $\Gamma _p(U_h)$ are said to be {\em homogeneous} of {\em
degree} $p$.  For a homogeneous element $x$ of $U_h$, the degree of
$x$ is denoted by $|x|$.

\subsection{The subalgebras $\Uq$ and $\Uqv$ of $U_h$}
\label{sec:subalgebras-uq-uqv}
Set
\begin{gather*}
  e= (v-v^{-1})E,\\
  F^{(n)} = F^n/[n]!,\\
  \tF^{(n)} = F^nK^n/[n]_q! = v^{-\hf n(n-1)}F^{(n)}K^n
\end{gather*}
for $n\ge 0$.  Let $\Uq$ denote the $\Zqq$-subalgebra of $U_h$ generated
by $K$, $K^{-1}$, $e$, and $\tF^{(n)}$ for $n\ge 1$.
The definition of $\Uq$ here is equivalent to that in \cite[Section
11]{H:uqsl2}.

Let $\Uqv$ denote the $\Zqq$-subalgebra of $\Uq$ generated by $K^2$,
$K^{-2}$, $e$, and $\tF^{(n)}$ for $n\ge 1$.  ($\Uqv$ is the same as
${\mathcal G}_0\Uq$ in \cite{H:uqsl2}.)  $\Uq$ is equipped with a
$(\modZ /2\modZ )$-graded $\Zqq$-algebra structure
\begin{equation}
  \label{eq:56}
  \Uq= \Uqv \oplus K\Uqv.
\end{equation}

Later we need the following formulas in $\Uq$.
\begin{gather}
  \label{eq:58}
  Ke=qeK,\quad
  K\tF^{(n)}=q^{-n}\tF^{(n)}K,\\
  \label{eq:59}
  \tF^{(m)}\tF^{(n)}=q^{-mn}\bbq{m+n}{m} \tF^{(m+n)},\\
  \label{eq:60}
  e^m\tF^{(n)}= \sum_{p=0}^{\min(m,n)}
    q^{-n(m-p)}{\bbq mp} \tF^{(n-p)}\{H-m-n+2p\}_{q,p} {e^{m-p}}.
\end{gather}
Here, for $i\in \modZ $ and $p\ge 0$, we set
\begin{equation*}
  \{H+i\}_{q,p} = \{H+i\}_q\{H+i-1\}_q\cdots\{H+i-p+1\}_q,
\end{equation*}
where
\begin{equation*}
  \{H+j\}_q = q^{H+j}-1 = q^jK^2-1
\end{equation*}
for $j\in \modZ $.  Note that $\{H+j\}_q,\{H+i\}_{q,p}\in \Zqq[K^2,K^{-2}]$.

The following is a ``$q$-version'' of \cite[Proposition 3.1]{H:uqsl2},
which can be easily proved using \eqref{eq:58}--\eqref{eq:60}.

\begin{lemma}
  \label{r39}
  $\Uq$ (resp. $\Uqv$) is freely spanned over $\Zqq$ by the elements
  $\tF^{(i)}K^je^k$ (resp. $\tF^{(i)}K^{2j}e^k$) with $i,k\ge 0$ and
  $j\in \modZ $.
\end{lemma}

Recall from \cite{H:uqsl2} that $\Uq$ inherits from $U_h$ a Hopf
 $\Zqq$-algebra structure, which can be easily verified using the
 following formulas, which we also need later.
\begin{gather}
  \label{eq:53}
  \Delta (K^i) = K^i\otimes K^i,\quad   S^{\pm 1}(K^i) = K^{-i},\\
  \label{eq:54}
  \Delta (e^n) = \sum_{j=0}^n {\bb nj}_q e^{n-j}K^j\otimes e^j,\quad
  \Delta (\tF^{(n)})=\sum_{j=0}^n \tF^{(n-j)}K^j\otimes \tF^{(j)},\\
  \label{eq:55}
  S^{\pm 1}(e^n)=(-1)^n q^{\hf n(n\mp1)} K^{-n} e^n,\quad
  S^{\pm 1}(\tF^{(n)}) =(-1)^n q^{-\hf n(n\mp1)} K^{-n}\tF^{(n)},\\
  \label{eq:4}
  \epsilon (K^i)=1,\quad \epsilon (e^n)=\epsilon (\tF^{(n)}) = \delta _{n,0}.
\end{gather}

For $n\ge 0$, the $n$-output comultiplication $\Delta ^{[n]}\zzzcolon U_h\rightarrow \Uhhon$ is
defined inductively by $\Delta ^{[0]}=\epsilon $, and $\Delta ^{[n+1]}=(\Delta ^{[n]}\otimes \id)\Delta $
for $n\ge 0$.  For $x\in U_h$ and $n\ge 1$, we write
\begin{equation*}
  \Delta ^{[n]}(x)=\sum x_{(1)}\otimes \cdots\otimes x_{(n)}.
\end{equation*}

\subsection{Adjoint action}
\label{sec:adjoint-action}
Let $\trr\zzzcolon U_h\ho U_h \rightarrow U_h$ denote the (left) adjoint action defined
by
\begin{equation*}
  \trr(x\otimes y)=x \trr y = \sum x_{(1)} y S(x_{(2)})
\end{equation*}
for $x,y\in U_h$.  We regard $U_h$ as a left $U_h$-module via the
adjoint action.  Since we have $\Uq\trr\Uq\subset \Uq$, we may regard $\Uq$
as a left $\Uq$-module.

For each homogeneous element $x\in U_h$, we have
\begin{gather}
  \label{e17}
  K^i\trr x = q^{i|x|} x\quad \text{for $i\in \modZ $},\\
  \label{e24}
  e^n\trr x =
  \sum_{j=0}^n (-1)^j q^{\hf j(j-1)+j|x|} {\bb nj}_q e^{n-j}x e^j
  \quad \text{for $n\ge 0$},\\
  \label{e64}
  \tF^{(n)}\trr x =
  \sum_{j=0}^n (-1)^j q^{-\hf j(j-1)+j|x|}\tF^{(n-j)}x\tF^{(j)}
  \quad \text{for $n\ge 0$}.
\end{gather}

\begin{proposition}
  \label{thm:19}
  $\Uqv$ is a left $\Uq$-submodule of $\Uq$.
\end{proposition}

\begin{proof}
  It follows from \eqref{e17}--\eqref{e64} that if $x$
  is a homogeneous element of $\Uqv$, then we have $y\trr x\in \Uqv$ for
  $y=K^i,e^n,\tF^{(n)}$ with $i\in \modZ $, $n\ge 0$.  Since these elements
  generate $\Uq$, we have the assertion.
\end{proof}

\subsection{Filtrations and completions}
\label{sec:filtration-completion}

Here we introduce filtrations for the algebras $\Uq$ and $\Uqv$, and
also for the tensor powers $\Uqvn$ of $\Uqv$.  These filtrations
produce the associated completions.  The definitions below are
equivalent to those in \cite{H:uqsl2}.

First, we consider the filtration and the completion for $\Uq$.  For
$p\ge 0$, set
\begin{equation*}
  \modF _p(\Uq)=\Uq e^p\Uq,
\end{equation*}
the two-sided ideal in $\Uq$ generated by $e^p$.  Let $\tUq$ denote
the ``completion in $U_h$'' of $\Uq$ with respect to the decreasing
filtration $\{\modF _p(\Uq)\}_{p\ge 0}$, i.e., $\tUq$ is the image of the
homomorphism
\begin{equation}
  \label{eq:6}
  \varprojlim_{p\ge 0} \Uq/\modF _p(\Uq) \rightarrow  U_h
\end{equation}
induced by $\Uq\subset U_h$.  (Conjecturally, \eqref{eq:6} is injective, see
\cite[Conjecture 7.2]{H:uqsl2}.) Clearly, $\tUq$ is a
$\Zqq$-subalgebra of $U_h$.

Second, we consider $\Uqv$.  For $p\ge 0$, set
\begin{equation*}
  \modF _p(\Uqv)=\modF _p(\Uq)\cap \Uqv=\Uqv e^p\Uqv.
\end{equation*}
We have a $(\modZ /2\modZ )$-grading $\modF _p(\Uq)=\modF _p(\Uqv)\oplus K\modF _p(\Uqv)$.
Let $\tUqv$ denote the ``completion in $U_h$'' of $\Uqv$ with respect
to $\{\modF _p(\Uqv)\}_{p\ge 0}$, i.e., the image of
\begin{equation*}
  \varprojlim_{p\ge 0}\Uqv/\modF _p(\Uqv)\rightarrow U_h.
\end{equation*}
Note that $\tUqv$ is a $\Zqq$-subalgebra of $\tUq$, and $\tUq$ has a
$(\modZ /2\modZ )$-graded $\Zqq$-algebra structure
\begin{equation*}
  \tUq=\tUqv\oplus K\tUqv
\end{equation*}
induced by \eqref{eq:56}.

Note that the elements of $\tUq$ (resp. $\tUqv$) are the elements of
$U_h$ that can be expressed as infinite sums
\begin{math}
  \sum_{i=0}^\infty \sum_{j=1}^{N_i} x_{i,j}e^iy_{i,j},
\end{math}
where $N_0,N_1,\ldots\ge 0$, and $x_{i,j},y_{i,j}\in \Uq$
(resp. $x_{i,j},y_{i,j}\in \Uqv$) for $i\ge 0$, $1\le j\le N_i$.

Now, we define the filtration for $\Uqvn$, $n\ge 1$, by
\begin{equation*}
  \modF _p(\Uqvn)
  =\sum_{i=1}^n\Uqvx{(i-1)}\otimes \modF _p(\Uqv)\otimes \Uqvx{(n-i)}\subset \Uqvn.
\end{equation*}
Thus, an element of $\Uqvn$ is in the $p$th filtration if and only if
it is expressed as a sum of terms each having at least one tensor
factor in the $p$th filtration.  Define the ``completed tensor
product'' $\tUqvn= \tUqv\tO\cdots\tO\tUqvn$ to be the ``completion in
$\Uhn$'' of $\Uqvn$ with respect to this filtration, i.e., the image
of the homomorphism
\begin{equation*}
  \varprojlim_{p\ge 0}\Uqvn/\modF _p(\Uqvn)\rightarrow \Uhn.
\end{equation*}
For $n=0$, it is natural to set
\begin{equation*}
  \modF _p(\Uqvx{0})=\modF _p(\Zqq)=
  \begin{cases}
    \Zqq&\text{if $p=0$,}\\
    0&\text{otherwise.}
  \end{cases}
\end{equation*}
Thus we have
\begin{equation*}
  (\tUqv)^{\tO 0} = \Zqq.
\end{equation*}

In what follows, we will also need the filtrations and completions of
other iterated tensor products of $\Uq$ and $\Uqv$, whose definitions
should be obvious from the above definitions.  For example,
$\Uq\otimes \Uqv$ has a filtration defined by
\begin{equation*}
  \modF _p(\Uq\otimes \Uqv)=\modF _p(\Uq)\otimes \Uqv+\Uq\otimes \modF _p(\Uqv).
\end{equation*}
Moreover, $\tUq\tO\tUqv$ is defined to be the image of
$\varprojlim_{p\ge 0}\modF _p(\Uq\otimes \Uqv)\rightarrow U_h^{\ho2}$.

\subsection{The Hopf algebra structure for $\tUq$}
\label{sec:Hopf-algebra-structure}

In this subsection, we show that the $\Zqq$-algebra $\tUq$ inherits
from $U_h$ a complete Hopf algebra structure over $\Zqq$.  (This fact
is observed in \cite{H:uqsl2} as a corollary to the case of a
$\Zvv$-form $\tU\cong\tUq\otimes _{\Zqq}\Zvv$.  However, it is convenient to
provide a direct proof here.)

The Hopf $\Zqq$-algebra structure of $\Uq$ induces a complete Hopf
algebra structure (with invertible antipode) over $\Zqq$ of $\tUq$,
since \eqref{eq:53}--\eqref{eq:4} imply
\begin{gather*}
  \Delta (\modF _p(\Uq))\subset \modF _{\floor{\frac{p+1}2}}(\Uq^{\otimes 2}),\\
  \epsilon (\Uq)\subset \Zqq,\quad \epsilon (\modF _1(\Uq))=0, \\
  S^{\pm 1}(\modF _p(\Uq))\subset \modF _p(\Uq),
\end{gather*}
for $p\ge 0$.  Here $\floor{\frac{p+1}2}$ denotes the largest integer smaller
than or equal to $\frac{p+1}2$.  The structure morphisms
\begin{equation*}
  \Delta \zzzcolon \tUq\rightarrow \tUq^{\tO2},\quad \epsilon \zzzcolon \tUq\rightarrow \Zqq,\quad S\zzzcolon \tUq\rightarrow \tUq
\end{equation*}
of $\tUq$ are equal to the restrictions of those of $U_h$ to $\tUq$.

A consequence of the above fact is that if $f\zzzcolon \Uhx i\rightarrow \Uhx j$,
$i,j\ge 0$, is a $\modQ [[h]]$-module homomorphism obtained from finitely
many copies of $\id_{U_h}$, $P_{U_h,U_h}$, $\mu $, $\eta $, $\Delta $, $\epsilon $ and
$S^{\pm 1}$ by taking completed tensor products and compositions, then we
have $f(\tUqx i)\subset \tUqx j$.  Here $P_{U_h,U_h}\zzzcolon \Uhx2\rightarrow \Uhx2$ is
defined by $P_{U_h,U_h}(\sum x\otimes y)=\sum y\otimes x$.

It follows that the adjoint action $\trr\zzzcolon \Uq\otimes \Uq\rightarrow \Uq$ induces a
left action
\begin{equation*}
  \trr \zzzcolon  \tUq\tO\tUq \rightarrow \tUq,
\end{equation*}
which is equal to the restriction of $\trr\zzzcolon U_h\ho U_h\rightarrow U_h$.  By
Proposition \ref{thm:19}, $\trr$ restricts to
\begin{gather*}
  \trr \zzzcolon  \tUq\tO\tUqv \rightarrow \tUqv.
\end{gather*}
Thus, $\tUqv$ is a left $\tUq$-submodule of $\tUq$.

\section{Braided Hopf algebra structure for $\tUqv$}
\label{sec3}

In this section, we recall a ribbon Hopf algebra structure for $U_h$
and show that the associated braided Hopf algebra structure for $U_h$
induces that for $\tUqv$.

\subsection{Ribbon structure for $U_h$}
\label{sec:ribbon-structure-uh}
The Hopf algebra $U_h$ has a ribbon Hopf algebra structure as follows.
 The universal $R$-matrix and its inverse are given by
\begin{gather}
  R = D\Bigl(\sum_{n\ge 0} v^{n(n-1)/2}
  \frac{(v-v^{-1})^n}{[n]!} F^n\otimes E^n\Bigr),\\
  \begin{split}
  R^{-1}
  &= \Bigl(\sum_{n\ge 0} (-1)^n v^{-n(n-1)/2}
  \frac{(v-v^{-1})^n}{[n]!}  F^n\otimes E^n\Bigr)D^{-1},
  \end{split}
\end{gather}
where
\begin{equation*}
  D=v^{\hf H\otimes H}=\exp (\frac h4 H\otimes H)\in U_h^{\ho2}.
\end{equation*}

In what follows, we use the following notations.
\begin{equation*}
  R=\sum\alpha \otimes \beta ,\quad R^{-1}=\sum\bar\alpha \otimes \bar\beta 
  \bigl(=\sum S(\alpha )\otimes \beta \bigr).
\end{equation*}

The ribbon element and its inverse are given by
\begin{gather*}
  \modr  = \sum S(\alpha )K^{-1}\beta =\sum\beta KS(\alpha ),\quad
  \modr ^{-1} = \sum \alpha K\beta =\sum\beta K^{-1}\alpha .
\end{gather*}
The associated grouplike element $\kappa \in U_h$ defined by
\begin{equation*}
  \kappa =\bigl(\sum S(\beta )\alpha \bigr)\modr ^{-1}
\end{equation*}
satisfies $\kappa =K^{-1}$.

We also use the following notations.
\begin{gather*}
  D = \sum D_{[1]}\otimes D_{[2]},\quad
  D^{-1} = \sum \bar D_{[1]}\otimes \bar D_{[2]}.
\end{gather*}
The following properties of $D$ are freely used in what follows.
\begin{gather}
  \sum D_{[2]}\otimes D_{[1]}=D,\quad
  (\Delta \otimes 1)(D)=D_{13}D_{23},\\
  (\epsilon \otimes 1)(D)=1,\quad
  (S\otimes 1)(D)=D^{-1},\\
  \label{e31}
  D^{\pm 1}(1\otimes x)=(K^{\pm |x|}\otimes x)D\quad \text{for homogeneous $x\in U_h$},
\end{gather}
where $D_{13}=\sum D_{[1]}\otimes 1\otimes D_{[2]}$ and
$D_{23}=\sum1\otimes D_{[1]}\otimes D_{[2]}=1\otimes D$.

We can easily obtain the following formulas.
\begin{gather}
  \label{eq:11}
  \begin{split}
  R
  &=D\Bigl(\sum_{n\ge 0} q^{\hf n(n-1)}\tF^{(n)}K^{-n}\otimes e^n\Bigr),
  \end{split}\\
  \label{e30}
  \begin{split}
    R^{-1}
    &=D^{-1}\Bigl(\sum_{n\ge 0} (-1)^n\tF^{(n)}\otimes K^{-n}e^n\Bigr),
  \end{split}\\
  \label{eq:65}
  \modr =\sum_{n\ge 0}(-1)^n\tF^{(n)}v^{-\hf H(H+2)}e^n,\\
  \label{eq:66}
  \modr ^{-1}=\sum_{n\ge 0}q^{\hf n(n-1)}
  \tF^{(n)}K^{-2n}v^{\hf H(H+2)}e^n.
\end{gather}
We have
$R^{\pm 1}\in D^{\pm 1}(\tUq\tO\tUq)$, and $\modr ^{\pm 1}\in v^{\mp\hf H(H+2)}\tUqv$.

\subsection{Braided Hopf algebra structure for $U_h$}
\label{sec:transmutation}

Let $\ModUh$ denote the category of $h$-adically complete left
$U_h$-modules and continuous left $U_h$-module homomorphisms.  The
category $\ModUh$ is equipped with a standard braided category
structure, where the braiding $\psi _{V,W}\zzzcolon V\otimes W\rightarrow W\otimes V$ of two objects
$V$ and $W$ in $\ModUh$ is defined by
\begin{equation*}
  \psi _{V,W}(v\otimes w) = \sum \beta w\otimes \alpha v\quad \text{for $v\in V$, $w\in W$}.
\end{equation*}
The inverse $\psi _{V,W}^{-1}\zzzcolon W\otimes V\rightarrow V\otimes W$ of $\psi _{V,W}$ is given by
\begin{equation*}
  \psi _{V,W}^{-1}(w\otimes v) = \sum S(\alpha ) v\otimes \beta w\quad \text{for
  $v\in V$, $w\in W$}.
\end{equation*}
We regard $U_h$ as an object of $\ModUh$, equipped
with the adjoint action.  For simplicity, we write
$\psi =\psi _{U_h,U_h}$.  Thus the braiding and its inverse for $U_h$ satisfy
\begin{gather*}
  \psi (x\otimes y) = \sum( \beta \trr y)\otimes (\alpha \trr x),\\
  \psi ^{-1}(x\otimes y) = \sum (S(\alpha )\trr y)\otimes (\beta \trr x)
\end{gather*}
for $x,y\in U_h$.

Let $\bUh=(U_h,\mu ,\eta ,\bD,\epsilon ,\bS)$ denote the transmutation
\cite{Majid:algebras,Majid:foundations} of $U_h$, which is a standard
braided Hopf algebra structure in $\ModUh$ associated to $U_h$.  Here
\begin{equation*}
  \mu \zzzcolon U_h^{\ho2}\rightarrow U_h,\quad
  \eta \zzzcolon \modQ [[h]]\rightarrow U_h,\quad
  \epsilon \zzzcolon U_h\rightarrow \modQ [[h]],
\end{equation*}
are the structure morphisms of $U_h$, and
\begin{equation*}
  \bD\zzzcolon U_h\rightarrow U_h^{\ho2}, \quad \bS\zzzcolon U_h\rightarrow U_h
\end{equation*}
are the twisted versions of comultiplication and antipode defined by
\begin{gather}
  \bD(x) = \sum x_{(1)}S(\beta )\otimes (\alpha \trr x_{(2)}),\\
  \label{eq:27}
  \bS(x) = \sum \beta \; S(\alpha \trr x)
\end{gather}
for $x\in U_h$.  The inverse $\bS^{-1}\zzzcolon U_h\rightarrow U_h$ of $\bS$ is given by
\begin{equation}
  \label{eq:28}
  \bS^{-1}(x) = \sum S^{-1}(\alpha \trr x) \beta .
\end{equation}

\subsection{Braided Hopf algebra structure for $\tUqv$}
\label{sec:braid-hopf-algebra}

Unlike $\tUq$, the even part $\tUqv$ of $\tUq$ does not inherit a Hopf
algebra structure from $U_h$.  However, we have the following.

\begin{theorem}
  \label{thm:26}
  The braided Hopf algebra structure of $\bUh$ induces a braided
  Hopf algebra structure with invertible antipode for $\tUqv$.  In
  other words, for
  \begin{equation*}
    f\in \{\psi ,\psi ^{-1},\mu ,\eta ,\bD,\epsilon ,\bS,\bS^{-1}\}
  \end{equation*}
  with $f\zzzcolon U_h^{\ho i}\rightarrow U_h^{\ho j}$ ($i,j\ge 0$) we have
  \begin{equation*}
    f((\tUqv)^{\tO i})\subset (\tUqv)^{\tO j},
  \end{equation*}
  and the induced map $f\zzzcolon (\tUqv)^{\tO i}\rightarrow (\tUqv)^{\tO j}$ is
  continuous.
\end{theorem}

\begin{corollary}
  \label{thm:20}
  Suppose that $f\zzzcolon U_h^{\ho i}\rightarrow U_h^{\ho j}$, $i,j\ge 0$, is a
  $\modQ [[h]]$-module homomorphism obtained from finitely many copies of
  \begin{gather*}
    1_{U_h}\zzzcolon U_h\rightarrow U_h,\quad \psi ^{\pm 1}\zzzcolon U_h^{\ho2}\rightarrow U_h^{\ho2},\quad
    \mu \zzzcolon U_h^{\ho2}\rightarrow U_h,\quad \eta \zzzcolon \modQ [[h]]\rightarrow U_h,\\
    \bD\zzzcolon U_h\rightarrow U_h^{\ho2},\quad \epsilon \zzzcolon U_h\rightarrow \modQ [[h]],\quad
    \bS^{\pm 1}\zzzcolon U_h\rightarrow U_h
  \end{gather*}
  by taking iterated tensor products and compositions.  Then we have
  \begin{equation*}
      f(\tUqvx i)\subset \tUqvx j.
  \end{equation*}
\end{corollary}

Theorem \ref{thm:26} follows immediately from the following.

\begin{proposition}
  \label{r51}
  If $f\in \{\psi ,\psi ^{-1},\mu ,\eta ,\epsilon ,\bS,\bS^{-1}\}$ with $f\zzzcolon U_h^{\ho
  i}\rightarrow U_h^{\ho j}$ ($i,j\ge 0$), then we have
  \begin{equation}
    f(\modF _p(\tUqvx i))\subset \modF _p(\tUqvx j)
  \end{equation}
  for $p\ge 0$.  Moreover, we have
  \begin{equation}
    \label{eq:39}
    \bD(\modF _p(\tUqv))\subset \modF _{\floor{\frac{p+1}2}}(\tUqvx2)
  \end{equation}
  for $p\ge 0$.
\end{proposition}

\begin{proof}

 The assertion is obvious for $f=\mu ,\eta ,\epsilon $.

The case $f=\psi ^{\pm 1}$ follows from the following formulas
\begin{gather}
  \psi (x\otimes y)
  =\sum_{n=0}^\infty  q^{-\hf n(n+1)+(|x|-n)|y|}
  (e^n\trr y)\otimes (\tF^{(n)}\trr x),\\
  \psi ^{-1}(x\otimes y)
  =\sum_{n=0}^\infty (-1)^n q^{-(|x|+n)|y|}
  (\tF^{(n)}\trr y)\otimes (e^n\trr x)
\end{gather}
for homogeneous elements $x,y\in U_h$.

Consider the case $f=\bS^{\pm 1}$.
Using \eqref{eq:27} and \eqref{eq:28}, we obtain
\begin{gather}
  \label{eq:16}
  \bS(x)=\sum_{n=0}^\infty q^{-\hf n(n+1)}
  e^nK^{|x|-n}S(\tF^{(n)}\trr x),\\
  \label{eq:17}
  \bS^{-1}(x)=\sum_{n=0}^\infty q^{\hf n(n-1)-n|x|}
  S^{-1}(\tF^{(n)}\trr x)K^{|x|-n}e^n.
\end{gather}
Using these formulas, we see easily that
$\bS(\modF _p(\tUqv))\subset \modF _p(\tUq)$.  Hence, it remains to show that if
$x,y\in \Uqv$ are homogeneous, then each term in \eqref{eq:16} and
\eqref{eq:17} is contained in $\Uqv$.  This follows, since for any
homogeneous $x\in \Uqv$ we have
\begin{math}
  K^{|x|}S^{\pm 1}(x)\in \Uqv
\end{math}
by \eqref{eq:53} and \eqref{eq:55}.

Finally, we prove \eqref{eq:39}.  By computation,
\begin{equation}
  \label{eq:18}
  \bD(x)
  =\sum_{n=0}^\infty (-1)^n q^{-n}
  x_{(1)}K^{-|x_{(2)}|}e^n\otimes (\tF^{(n)}\trr x_{(2)}),
\end{equation}
where $\Delta (x)=\sum x_{(1)}\otimes x_{(2)}$.  Using this formula, we easily
see that $\bD(\modF _p(\tUqv))\subset \modF _{\floor{\frac{p+1}2}}(\tUq\otimes \tUqv)$.
It suffices to show that if $x\in \Uqv$ is homogeneous, then each term
in \eqref{eq:18} is contained in $\Uqvx2$.  This follows, since we
have
\begin{gather}
  \label{eq:22}
  \Delta (\Uqv)\subset \Uqv\otimes (\Uqv)_0+K\Uqv\otimes (\Uqv)_1,
\end{gather}
where $(\Uqv)_0$ (resp. $(\Uqv)_1$) denotes the $\Zqq$-submodule of
$\Uqv$ spanned by the homogeneous elements of even (resp. odd)
degrees.  The inclusion \eqref{eq:22} can be easily verified using
\eqref{eq:53} and \eqref{eq:54}.
\end{proof}

\begin{remark}
  \label{r15}
  $\tUq$ also has a braided Hopf algebra structure inherited from that
  of $\bUh$, i.e., Theorem \ref{thm:26} holds if we replace
  $\tUqv$ with $\tUq$.  Moreover, the $(\modZ /2\modZ )$-grading for $\tUq$ is
  compatible with the braided Hopf algebra structure:
  \begin{gather*}
    \psi ^{\pm 1}(K^i\tUqv\tO K^j\tUqv)\subset K^j\tUqv\tO K^i\tUqv,\\
    \bD(K^i\tUqv)\subset K^i\tUqv \tO K^i\tUqv,\quad
    \bS^{\pm 1}(K^i\tUqv)\subset K^i\tUqv,
  \end{gather*}
  for $i,j\in \{0,1\}$.
\end{remark}

\section{Universal $sl_2$ invariant of bottom tangles}
\label{sec4}

In this section, we recall the definition of the universal invariant
of bottom tangles, and prove necessary results.

\subsection{Bottom tangles}
\label{sec:bottom-tangles}

Here we recall from \cite{H:universal} the notion of bottom tangles.

An $n$-component {\it bottom tangle} $T=T_1\cup \cdots\cup T_n$ is a framed
tangle in a cube, which is drawn as a diagram in a rectangle as usual,
consisting of $n$ arcs $T_1,\ldots,T_n$ such that for each $i=1,\ldots,n$ the
component $T_i$ runs from the $2i$th endpoint on the bottom to the
$(2i-1)$st endpoint on the bottom, where the endpoints are counted
from the left.  For example, see Figure \ref{fig:bt} (a).
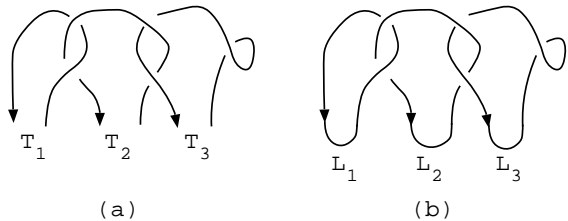
\FIG{bt}{(a) A $3$-component bottom tangle $T=T_1\cup T_2\cup T_3$.  (b) Its
closure $\cl(L)=L_1\cup L_2\cup L_3$.}{height=30mm} (Here and in what
follows, we use the blackboard framing convention.)

For each $n\ge 0$, let $\BT_n$ denote the set of the isotopy classes of
$n$-component bottom tangles.  Set $\BT=\bigcup_{n\ge 0}\BT_n$.

The {\em closure} $\cl(T)$ of $T$ is the $n$-component, oriented,
ordered framed link in $S^3$, obtained from $T$ by pasting a
``$\cup $-shaped tangle'' to each component of $L$, as depicted in Figure
\ref{fig:bt} (b).  For any oriented, ordered framed link $L$, there is
a bottom tangle whose closure is isotopic to $L$.

The {\em linking matrix} of a bottom tangle $T=T_1\cup \cdots\cup T_n$ is
defined as that of the closure $T$.  Thus the linking number of $T_i$
and $T_j$, $i\neq j$, is defined as the linking number of the
corresponding components in $\cl(T)$, and the framing of $T_i$ is
defined as the framing of the closure of $T_i$.

A link or a bottom tangle is called {\em algebraically-split} if the
linking matrix is diagonal.

For $n\ge 0$, let $\BT^0_n$ denote the subset of $\BT_n$ consisting of
algebraically-split, $0$-framed bottom tangles.  Set
$\BT^0=\bigcup_{n\ge 0}\BT^0_n\subset \BT$.

\subsection{Universal $sl_2$ invariant of bottom tangles}
\label{sec:univ-sl2-invar}

For each ribbon Hopf algebra $H$, there is a ``universal invariant''
of links and tangles from which one can recover the operator
invariants, such as the colored Jones polynomials.  Such universal
invariants has been studied in
\cite{Lawrence:89,Lawrence:90,Reshetikhin:89,Lee:90,Ohtsuki:93,Kauffman:93,Kerler:97,Kauffman-Radford:01}.
Here we need only the case of bottom tangles, which is described in
\cite{H:universal}.

For $T=T_1\cup \cdots\cup T_n\in \BT_n$, we define the {\em universal $sl_2$
  invariant} $J_T\in U_h^{\ho n}$ of $T$ as follows.  We choose a
  diagram for $T$, which is obtained from copies of {\em fundamental
  tangle}, see Figure \ref{fig:fundamental}, by pasting horizontally
  and vertically.
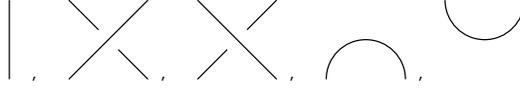
\FIG{fundamental}{Fundamental tangles: vertical line, positive and
  negative crossings, local minimum and local maximum.  Here the
  orientations are arbitrary.}{height=12mm}
For each copy of fundamental tangle in the diagram of $T$, we put
elements of $U_h$ with the rule described in Figure
\ref{fig:fundamental2}.
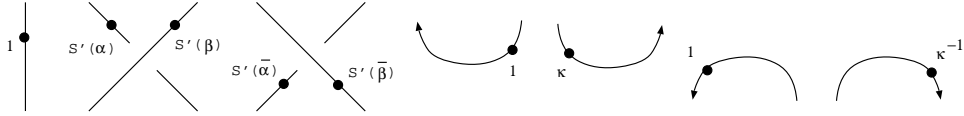
\FIG{fundamental2}{How to put elements of $U_h$
on the strings.  For each string in the positive and the negative
crossings, ``$S'$'' should be replaced with $\id$ if the string is
oriented downward, and by $S$ otherwise.}{height=15mm}
We set
\begin{equation*}
  J_T =\sum J_{(T_1)}\otimes \cdots\otimes J_{(T_n)}\in U_h^{\ho n},
\end{equation*}
where for each $i=1,\ldots,n$, the $i$th tensorand $J_{(T_i)}$ is defined
to be the product of the elements put on the component $T_i$.  Here
the elements read off along each component are written from right to
left.  Then $J_T$ does not depend on the choice of diagram, and
defines an isotopy invariant of bottom tangles.

\subsection{Universal $sl_2$-invariant of algebraically-split,
  $0$-framed bottom tangles}
\label{sec:invariant}

For any left $U_h$-module $W$, let $\Inv(W)$ denote the {\em invariant
  part} of $W$, defined by
\begin{equation*}
  \Inv(W) =\{w\in W\zzzvert  x\cdot w=\epsilon (x)w\quad \forall x\in U_h\}.
\end{equation*}

Recall that we regard $U_h$ as a left $U_h$-module via the adjoint
action.  For $n\ge 0$, the completed tensor product $\Uhn$ is equipped
with a left $U_h$-module structure $\trrn$ in the standard way: For
$x=\sum x_1\otimes \cdots\otimes x_n\in \Uhn$ and $y\in U_h$ we have
\begin{equation*}
  y\trrn x = \sum (y_{(1)}\trr x_1)\otimes \cdots\otimes (y_{(n)}\trr x_n).
\end{equation*}
In particular, $\Uhx0=\modQ [[h]]$ is given the trivial left $U_h$-module
structure.  For any subset $X\subset \Inv(U_h^{\ho n})$, we set
\begin{equation*}
  \Inv(X) = \Inv(U_h^{\ho n})\cap X.
\end{equation*}

For $n\ge 0$, we set
\begin{equation*}
  \modK _n= \Inv(\tUqvn)\subset \tUqvn.
\end{equation*}
One can easily see that $\modK _n$ is the $\tUq$-invariant part of the
$\tUq$-module $\tUqvn$, i.e., we have
\begin{equation*}
  \modK _n =\{x\in \tUqvn\zzzvert y\trrn x=\epsilon (y)x\quad \text{for all $y\in \tUq$}\}.
\end{equation*}

The main result of this subsection is the following.

\begin{theorem}
  \label{thm:8}
  If $T\in \BT^0_n$, $n\ge 0$, then we have $J_T\in \modK _n$.
\end{theorem}

To prove Theorem \ref{thm:8}, we use the following two results from
\cite{H:universal}.

\begin{proposition}[$U_h$ case of {\cite[Proposition 8.2]{H:universal}}]
  \label{r19}
  For any $n$-component bottom tangle $T\in \BT_n$, we have
  $J_T\in \Inv(U_h^{\ho n})$.
\end{proposition}

\begin{proposition}[$U_h$ case of {\cite[Corollary 9.15]{H:universal}}]
  \label{lem:1}
  Let $X_n\subset U_h^{\ho n}$, $n\ge 0$, be subsets satisfying the following
  conditions.
  \begin{enumerate}
  \item $1\in X_0$, $1\in X_1$, and $J_B\in X_3$.  Here $B\in \BT^0_3$ is the
    {\em Borromean tangle} depicted in Figure
    \ref{fig:borromean}.\FIG{borromean}{The Borromean tangle
    $B\in \BT^0_3$.}{height=15mm}
  \item If $x\in X_l$ and $y\in X_m$ with $l,m\ge 0$, then $x\otimes y\in X_{l+m}$.
  \item For $p,q\ge 0$ and $f\in \{\psi ^{\pm 1},\mu ,\bD,\bS\}$ with
    $f\zzzcolon U_h^{\ho i}\rightarrow U_h^{\ho j}$, we have
    \begin{equation*}
      (1^{\otimes p}\otimes f\otimes 1^{\otimes q})(X_{p+i+q})\subset X_{p+j+q}.
    \end{equation*}
  \end{enumerate}
  Then, for any $T\in \BT^0_n$, we have $J_T\in X_n$.
\end{proposition}

\begin{proof}[Proof of Theorem \ref{thm:8}]
  By Proposition \ref{r19}, it suffices to show that $J_T\in \tUqvx n$.
  We have only to verify the conditions in Proposition \ref{lem:1},
  where we set $X_n=\tUqvx n$.

  The condition (2) is obvious.

  The condition (3) follows from Corollary \ref{thm:20}.

  To prove (1), it suffices to prove $J_B\in \tUqvx3$.  Using Figure
  \ref{fig:borromean}, we obtain
  \begin{equation*}
    J_B
    = \sum
    \bar\alpha _3 \beta _1         \alpha _3 S^2(\bar\beta _1)\otimes 
    \bar\alpha _1 \beta _2         \alpha _1 S^2(\bar\beta _2)\otimes 
    \bar\alpha _2 S^{-2}(\beta _3) \alpha _2 \bar\beta _3,
  \end{equation*}
  where $R=\sum \alpha _i\otimes \beta _i$ and $R^{-1}=\sum \bar\alpha _i\otimes \bar\beta _i$ for
  $i=1,2,3$.  Using \eqref{eq:11} and \eqref{e30}, we obtain
    \begin{equation*}
    \begin{split}
      J_B
      &=
      \sum_{m_1,m_2,m_3,n_1,n_2,n_3\ge 0}
      (-1)^{n_1+n_2+n_3}
      q^{-\hf m_1(m_1+1)-\hf m_2(m_2+1)-\hf m_3(m_3+1)}\\&\quad
      \tF^{(n_3)}\bar D'_3 e^{m_1}D''_1        D'_3\tF^{(m_3)} S^2(\bar D''_1e^{n_1})\otimes 
      \tF^{(n_1)}\bar D'_1 e^{m_2}D''_2         D'_1\tF^{(m_1)} S^2(\bar
      D''_2e^{n_2})\\
      &\quad \otimes \tF^{(n_2)}\bar D'_2 S^{-2}(e^{m_3}D''_3) D'_2\tF^{(m_2)} \bar D''_3e^{n_3},
    \end{split}
  \end{equation*}
    where $D=\sum D'_i\otimes D''_i$ and $D^{-1}=\sum \bar D'_i\otimes \bar D''_i$
    for $i=1,2,3$.  We slide the tensor factors of the copies of
    $D^{\pm 1}$ using \eqref{e31} so that these copies cancel at the cost
    of inserting powers of $K$.  Thus we obtain
    \begin{equation}
      \label{eq:7}
      \begin{split}
	J_B
	&=
	\sum_{m_1,m_2,m_3,n_1,n_2,n_3\ge 0}q^{m_3+n_3}
	(-1)^{n_1+n_2+n_3}
	q^{\sum_{i=1}^3( -\hf m_i(m_i+1)-n_i+m_im_{i+1}-2m_in_{i-1})}\\
	&\quad \tF^{(n_3)} e^{m_1}   \tF^{(m_3)} e^{n_1}K^{-2m_2}\otimes 
	\tF^{(n_1)} e^{m_2}    \tF^{(m_1)} e^{n_2}K^{-2m_3}\otimes 
	\tF^{(n_2)} e^{m_3}    \tF^{(m_2)} e^{n_3}K^{-2m_1},
      \end{split}
    \end{equation}
    where the index $i$ should be considered modulo $3$.  Each term in
    \eqref{eq:7} is in $\Uqvx3$.  For any $p\ge 0$, all but finitely
    many terms in \eqref{eq:7} involve $e^r$ with $r\ge p$, and therefore
    are contained in $\modF _p(\Uqvx3)$.  Therefore, we have
    $J_B\in \tUqvx3$.
\end{proof}

\subsection{Universal $sl_2$ invariant of bottom knots}
\label{sec:case-n=1}

By a {\em bottom knot}, we mean a $1$-component bottom tangle.  In
what follows, we assume that bottom knots are given $0$-framing.  Thus
the set of bottom knots (up to isotopy) is $\BT^0_1$.  By Theorem
\ref{thm:8}, for any bottom knot $T\in \BT^0_1$, we have
\begin{equation}
  \label{e34}
  J_T\in \Inv(\tUqv)=Z(\tUqv).
\end{equation}

Let us recall from \cite{H:uqsl2} the structure of $Z(\tUqv)$.
(In \cite{H:uqsl2}, $\tUqv$ is denoted by $\mathcal{G}_0(\tUq)$.)
Set
\begin{equation*}
  C = (v-v^{-1})^2 FE+vK+v^{-1}K^{-1}\in Z(U_h),
\end{equation*}
which is a well-known central element.  We have
\begin{equation*}
  C= (v-v^{-1}) \tF^{(1)}K^{-1}e+vK+v^{-1}K^{-1}
  \in vK\Uqv\cap Z(U_h).
\end{equation*}
Hence $C^2\in \Uqv\cap Z(U_h)=Z(\Uqv)$.  As a $\Zqq$-algebra, $Z(\Uqv)$ is
freely generated by $C^2$, i.e., we have $Z(\Uqv)\cong\Zqq[C^2]$.  For
$p\ge 0$, set
\begin{equation*}
  \sigma _p = \prod_{i=1}^p(C^2-(q^i+2+q^{-i}))\in Z(\Uqv),
\end{equation*}
which is a monic polynomial of degree $p$ in $C^2$.  Therefore,
\begin{equation*}
  Z(\Uqv) = \Span_{\Zqq}\{\sigma _p\zzzvert p\ge 0\}.
\end{equation*}

As for the center of the completion $\tUqv$, we have the following.

\begin{theorem}[{\cite[Theorem 11.2]{H:uqsl2}}]
  \label{thm:10}
  The isomorphism $Z(\Uqv)\cong \Zqq[C^2]$ induces an isomorphism
  \begin{equation*}
    Z(\tUqv) \cong \varprojlim_{p\ge 0} \Zqq[C^2]/(\sigma _p).
  \end{equation*}
  Thus, each element in $Z(\tUqv)$ is uniquely expressed as an infinite sum
  $\sum_{p\ge 0} a_p \sigma _p$, where $a_p\in \Zqq$ for $p\ge 0$.
\end{theorem}

This implies the following.

\begin{theorem}
  \label{thm:23}
  If $T$ is a bottom knot, then $J_T$ is uniquely expressed as
  \begin{equation}
    \label{eq:74}
    J_T= \sum_{p\ge 0} a_p(T)\sigma _p,
  \end{equation}
  where $a_p(T)\in \Zqq$ for $p\ge 0$.
\end{theorem}

Note that the $a_p(T)$ are invariants of a bottom knot $T$.  In
Section \ref{sec6}, we give a formula which express $a_p(T)$ using the
colored Jones polynomials of the closure of~$T$.

\begin{remark}
  \label{r6}
  There is an obvious one-to-one correspondence between bottom knots
  and {\em string knots} (i.e., a string link consisting of just one
  arc component running from the above to the bottom).  A bottom knot
  and the corresponding string knot have the same value of the
  universal invariant.  Therefore, we have the result announced in
  \cite[Theorem 2.1]{H:rims2001}, \cite[Theorem 1.2]{H:uqsl2}, which
  is the string knot version of Theorem \ref{thm:23}.
\end{remark}

\section{Colored Jones polynomials}
\label{sec5}

In this section, we recall the definition of the colored Jones
polynomials of framed links.

\subsection{Finite-dimensional representations of $U_h$}
\label{sec:representations}
By a {\em finite-dimensional representation} of $U_h$, we mean a left
$U_h$-module which is free of finite rank as a $\modQ [[h]]$-module.

It is well known that, for each $n\ge 0$, there is exactly one
irreducible finite-dimensional representation $\V_n$ of rank $n+1$ up
to isomorphism, which corresponds to the $(n+1)$-dimensional
irreducible representation of $sl_2$.

The structure of $\V_n$ is as follows.  Let $\vn0\in \V_n$ denote a
highest weight vector of $\V_n$, which is characterized by $E\vn0=0$,
$H\vn0=n\vn0$ and $U_h\vn0=\V_n$.  It is useful to define the other
basis elements by
\begin{equation*}
  \vn i= \tF^{(i)}\vn0\quad \text{for $i=1,\ldots,n$}.
\end{equation*}
Then the action of $\Uq$ on $\V_n$ is given by
\begin{gather}
  K^{\pm 1} \vn i= v^{n-2i}\vn i,\\
  e^m \vn i= \{n-i+m\}_{q,m} \vn{i-m},\\
  \tF^{(m)} \vn i= q^{-mi}\bbq{i+m}{m}\vn{i+m},
\end{gather}
for $i=0,\ldots,n$ and $m\ge 0$, where we understand $\vn i=0$ unless
$0\le i\le n$.

\subsection{Colored Jones polynomials}
\label{sec:color-jones-polyn}

Let $L=L_1\cup \cdots\cup L_m$ be an $m$-component, framed, oriented, ordered
link, and let $W_1,\ldots,W_m\ge 0$ be finite-dimensional representations of
$U_h$.  We consider each $W_i$ as a {\em color} attached to the
component $L_i$.  Then the colored Jones polynomial $J_L(W_1,\ldots,W_m)$
of the colored link $(L;W_1,\ldots,W_m)$ is defined as follows.

We choose a diagram of $L$ which is obtained by pasting copies of the
fundamental tangles, as in the definition of the universal $sl_2$
invariant in Section \ref{sec:univ-sl2-invar}.  To each copy of a
fundamental tangle in $L$, we associate a left $U_h$-module
homomorphism.  To a vertical line contained in $L_i$, we associate to
the left $U_h$-module $W_i$ if the line is oriented downward, and the
dual $W_i^*$ of $W_i$ if the line is oriented upward.  To a positive
crossing, we associate a braiding operator
\begin{equation*}
  \psi _{W,W'}\zzzcolon W\otimes W'\rightarrow W'\otimes W, \quad x\otimes y\mapsto \sum\beta y\otimes \alpha x,
\end{equation*}
where $W$ (resp. $W'$) are the left $U_h$-modules associated to the
string which connects the upper left corner and the lower right corner
(resp. the upper right corner and the lower left corner).  To a
negative crossing, we associate
\begin{equation*}
  \psi _{W',W}^{-1}\zzzcolon W\otimes W'\rightarrow W'\otimes W, \quad x\otimes y\mapsto \sum\bar\alpha y\otimes \bar\beta x,
\end{equation*}
where $W$ and $W'$ are as above.  To the last four tangles in Figure
\ref{fig:fundamental2}, we associate the following operators
\begin{gather*}
  \ev_{W}\zzzcolon W^*\otimes W\rightarrow \modQ [[h]], \quad f\otimes x\mapsto f(x),\\
  \ev'_{W}\zzzcolon W\otimes W^*\rightarrow \modQ [[h]], \quad x\otimes f\mapsto f(\kappa x),\\
  \coev_{W}\zzzcolon \modQ [[h]]\rightarrow W\otimes W^*, \quad 1\mapsto \sum_ix_i\otimes x^i,\\
  \coev'_{W}\zzzcolon \modQ [[h]]\rightarrow W^*\otimes W, \quad 1\mapsto \sum_ix^i\otimes \kappa ^{-1}x_i,
\end{gather*}
respectively, where $W=W_i$ if $L_i$ is the component of $L$
containing the fundamental tangle, and where the $x_i$ are a basis of
$W$ and the $x^i\in W^*$ are the dual basis.

By tensoring and composing the operators, we obtain a left
$U_h$-module homomorphism from $\modQ [[h]]$ to $\modQ [[h]]$, and we define
$J_L(W_1,\ldots,W_m)$ to be the trace of this homomorphism.

Note that ``$J_L$'' alone denotes the universal $sl_2$ invariant of
$L$.  The above definition of $J_L(W_1,\ldots,W_n)$ is an abuse of
notation, but it should not cause confusion.

It is well known that $J_L\in q^{p/4}\Zqq\subset \modZ [q^{\pm 1/4}]$, where
$p\in \modZ /4\modZ $ depends only on the linking matrix of $L$ and on
$n_1,\ldots,n_m$.  The invariant is normalized so that
$J_{\emptyset}()=1$, and $J_U(\V_n)=[n+1]$ for $n\ge 0$, where $U$
denotes the $0$-framed unknot.

\subsection{Quantum trace and the colored Jones polynomial}
\label{sec:quantum-trace}

If $V$ is a finite-dimensional representation of $U_h$, then the
{\em quantum trace} $\trqV(x)$ in $V$ of an element $x\in U_h$ is
defined by
\begin{equation*}
  \trqV(x) = \tr^V(\rho _V(\kappa x)) \Bigl(= \tr^V(\rho _V(K^{-1}x))\Bigr)\in \modQ [[h]],
\end{equation*}
where $\rho _V\zzzcolon U_h\rightarrow \End(V)$ denotes the left action of $U_h$ on $V$,
and $\tr^V\zzzcolon \End(V)\rightarrow \modQ [[h]]$ denotes the trace in $V$.  We have a
continuous left $H$-module homomorphism
\begin{equation*}
  \trqV\zzzcolon U_h\rightarrow \modQ [[h]].
\end{equation*}

Let $L=L_1\cup \cdots\cup L_m$ be an $m$-component, ordered, framed oriented
link in $S^3$.  Choose $T\in \BT_m$ such that $\cl(T)$ is isotopic to
$L$.  As explained in \cite[Section 1.2]{H:universal} we have
\begin{equation}
  \label{e1}
  J_L(W_1,\ldots,W_m)
  =(\tr_q^{W_1}\otimes \tr_q^{W_2}\otimes \cdots\otimes \tr_q^{W_m})(J_T)
\end{equation}
for finite dimensional representations $W_1,\ldots,W_m$ of $U_h$.

\subsection{Representation rings}
\label{sec:representation-ring-}

For a commutative ring $A$ with unit, let $\modR _A$ denote the $A$-algebra
\begin{equation*}
  \modR _A = \Span_A\{\V_n\zzzvert n\ge 0\}.
\end{equation*}
with the multiplication induced by tensor product.  Since each
finite-dimensional representation of $U_h$ is a direct sum of copies
of $\V_n$, $n\ge 0$, we may regard $\modR _A$ as the {\em representation
ring} of $U_h$ over $A$.

By the well-known isomorphism of left $U_h$-modules
\begin{equation*}
  \V_m\otimes \V_n \cong
  \V_{|m-n|}\oplus\V_{|m-n|+2}\oplus\cdots\oplus\V_{m+n},
\end{equation*}
we have the identity in $\modR _A$
\begin{equation*}
  \V_m\V_n
  = \V_{|m-n|}+\V_{|m-n|+2}+\cdots+\V_{m+n}.
\end{equation*}
As an $A$-algebra, $\modR _A$ is freely generated by $\V_1$, i.e.,
$\modR _A\cong A[\V_1]$.
Thus we identify $\modR _A$ with $A[\V_1]$.

For $y=\sum_n a_n \V_n$ ($a_n\in \modQ [[h]]$) and $x\in U_h$, we set
\begin{equation*}
  \trq^{y}(x) = \sum_{n}a_n\trq^{\V_n}(x).
\end{equation*}
Thus we have bilinear maps
\begin{gather*}
  \trq^{-}(-)\zzzcolon \modR _{\Zvv}\times U_h\rightarrow \modQ [[h]],\\
  \trq^{-}(-)\zzzcolon \modR _{\modQ (v)}\times U_h\rightarrow \modQ ((h)),
\end{gather*}
where $\modQ ((h))$ denote the quotient field of $\modQ [[h]]$.

Similarly, for a link $L$ we have multilinear maps
\begin{gather*}
  J_L\zzzcolon \modR _{\Zvv}\times \cdots\times \modR _{\Zvv}\rightarrow  \modZ [v^{1/2},v^{-1/2}],\\
  J_L\zzzcolon \modR _{\modQ (v)}\times \cdots\times \modR _{\modQ (v)}\rightarrow  \modQ (v^{1/2}).
\end{gather*}
If the framing of every component of $L$ is even, then $J_L$ above
take values in $\Zvv$, $\modQ (v)$, respectively.

\section{Knots}
\label{sec6}

In this section, we study the relationships between the universal
invariant for a bottom knot $T$ and the colored Jones polynomials for
the closure of $T$.

\subsection{The elements $P_n$ and $P''_n$}
\label{sec:elements-pn-pn}

For each $n\ge 0$ set
\begin{gather*}
  P_n = \prod_{i=0}^{n-1}(\V_1 - v^{2i+1}-v^{-2i-1}) \in \modR _{\Zvv},\\
  P''_n = P_n/\BB{2n+1}{2n}\in \modR _{\modQ (v)}.
\end{gather*}

Since $P_n$ is a monic polynomial in $\V_1$ of degree $n$ for each
$n\ge 0$, it follows that the $P_n$ form a basis of $\modR _{\Zvv}$.  We have
the following base change formulas.

\begin{lemma}
  \label{r16}
  For $n\ge 0$, we have
  \begin{gather}
    \label{eq:97}
    P_n=\sum_{i=0}^n(-1)^{n-i}\frac{[2i+2]}{[n+i+2]}
    \bb{2n+1}{n+1+i}\V_i,\\
    \label{eq:100}
    \V_n =\sum_{i=0}^n \bb{n+i+1}{2i+1} P_i.
  \end{gather}
\end{lemma}

\begin{remark}
  \label{r17}
  A formula essentially the same as \eqref{eq:100} formulated using
  the Kauffman bracket skein module of a solid torus has appeared in
  \cite[Proposition 5.1.]{H:iias2000}, \cite[Equation 47]{Masbaum}.
  (The former contained an incorrect sign.)
\end{remark}

\begin{proof}
  The proofs of \eqref{eq:97} and \eqref{eq:100} are by inductions on $n$
  using
  \begin{gather*}
    P_n=P_{n-1}(\V_1-(v^{2n-1}+v^{-2n+1}))\quad (n\ge 1),\\
    \V_n=\V_1\V_{n-1}-\V_{n-2}\quad (n\ge 2),
  \end{gather*}
  respectively.  For an alternative proof of \eqref{eq:100}, see
  Remark \ref{thm:91}.
\end{proof}

\subsection{The reduced Jones polynomials of knots}
\label{sec:basis-pn-}

We prove the following result in Section \ref{sec:proof-prop-refr22}.

\begin{proposition}
  \label{r22}
  For $m,n\ge 0$ we have
  \begin{gather}
    \trq^{P''_m}(\sigma _n) = \delta _{m,n}.
  \end{gather}
\end{proposition}

Using Proposition \ref{r22} we obtain the following.

\begin{theorem}
  \label{r18}
  For a bottom knot $T\in \BT_1$ with closure $K=\cl(T)$, we have
  \begin{equation}
    \label{e9}
    J_T = \sum_{n\ge 0}J_K(P''_n)\sigma _n.
  \end{equation}
\end{theorem}

\begin{proof}
  We express $J_T$ as in Theorem \ref{thm:23}.  Then, applying
  $\trq^{P''_m}$ to \eqref{eq:74}, we have by Proposition \ref{r22}
  \begin{equation*}
    \trq^{P''_m}(J_T)= \sum_{p\ge 0} a_p(T)\trq^{P''_m}(\sigma _p)
    =a_p(T).
  \end{equation*}
  Therefore, we have \eqref{e9}.
\end{proof}

Theorems \ref{thm:23} and \ref{r18} implies that for a knot $K$, we
have $J_K(P''_n)\in \Zqq$ for $n\ge 0$.  We call $J_K(P''_n)$ the {\em
$n$th reduced Jones polynomial} of $K$.  It is clear that
$J_K(P''_0)=1$.  For $n=1$, we have
\begin{equation*}
  J_K(P''_1) = (J_K(\V_1)-[2])/\{2\}\{3\}.
\end{equation*}
The quantity $-q^{-2}J_K(P''_1)$ is known as the {\em reduced Jones
polynomial} of $K$, and appears in the original paper of Jones
\cite[Proposition 12.5]{Jones2}.  Some examples of $J_K(P''_n)$ are
given in Section \ref{sec:surg-borr-rings}.

It is clear that the universal invariant $J_T$ determines the colored
Jones polynomials $J_K(\V_n)$, hence the reduced Jones polynomials
$J_K(P''_n)$.  Conversely, Theorem \ref{r18} implies that the
universal invariant $J_T$ is determined by the $J_K(P''_n)$, hence by
the $J_K(\V_n)$ as is presumably well known.

By \eqref{eq:100}, we have
\begin{equation}
  \label{e89}
  J_K(\V_n) =
  \sum_{i=0}^n\frac{\BB{n+1+i}{2i+1}}{\{1\}}J_K(P''_i),
\end{equation}
where the sum may be replaced with $\sum_{i=0}^\infty $ since
$\BB{n+1+i}{2i+1}=0$ for $i>n$.  (Conversely, one can use
\eqref{eq:97} to obtain a formula for the reduced Jones polynomials in
terms of the colored Jones polynomials.)  \eqref{e89} implies the
following.

\begin{proposition}
  \label{r36}
  If $K$ is a knot, then for each $n\ge 0$ the $n$th colored Jones
  polynomial $J_K(\V_n)$ is determined modulo $(\{2n+1\}_{2n})$ by
  $J_K(\V_0),J_K(\V_1),\ldots,J_K(\V_{n-1})$.
\end{proposition}

\subsection{Proof of Proposition \ref{r22}}
\label{sec:proof-prop-refr22}

\subsubsection{The homomorphism $\xi \zzzcolon \modR _{\modQ [[h]]}\rightarrow Z(U_h)$}
\label{sec:isomorphism-ztu}

Let $c_+$ denote the $2$-component bottom tangle depicted in Figure
\ref{fig:c+}.  \FIG{c+}{The bottom tangle $c_+\in \BT_2$.}{height=16mm}
We have
\begin{equation*}
  J_{c_+}= (S\otimes 1)(R_{21}R)=\sum S(\alpha )S(\beta ')\otimes \alpha '\beta \in \Inv(\Uhx2),
\end{equation*}
where $R=\sum\alpha \otimes \beta =\sum\alpha '\otimes \beta '$.

Define a continuous $\modQ [[h]]$-algebra homomorphism
\begin{math}
  \xi \zzzcolon \modR _{\modQ [[h]]} \rightarrow  Z(U_h)
\end{math}
by
\begin{equation*}
  \xi (y) :=  (1\otimes \trq^y)(J_{c_+})
  =\sum S(\alpha )S(\beta ')\trq^y(\alpha '\beta )
\end{equation*}
for $y\in \modR _{\modQ [[h]]}$.  Indeed, we can verify
\begin{equation}
  \label{e35}
  \xi (VV')=\xi (V)\xi (V')
\end{equation}
graphically as depicted in Figure \ref{fig:T01}.  \FIG{T01}{A graphical
proof of \eqref{e35}.}{height=22mm}  We can also verify $\xi (\V_1)=C$ by
computation.  Since the $C^i$, $i\ge 0$, are linearly independent in
$Z(U_h)$, it follows that $\xi $ is injective.

\subsubsection{The Hopf link pairing $\langle ,\rangle \zzzcolon \modR _{\Zvv}\times \modR _{\Zvv}\rightarrow \Zvv$}
\label{sec:pairing-hhh}

Define a symmetric bilinear form
\begin{equation*}
  \langle ,\rangle \zzzcolon \modR _{\modQ [[h]]}\times \modR _{\modQ [[h]]}\rightarrow \modQ [[h]]
\end{equation*}
by
\begin{equation}
  \label{eq:85}
  \langle x,y\rangle :=(\trq^x\otimes \trq^y)(J_{c_+})= \trq^x(\xi (y)) = J_H(x,y),
\end{equation}
for $x,y\in U_h$, where $H=H_1\cup H_2=\cl(c_+)$ denotes the $0$-framed
Hopf link with linking number $-1$.

It is well known (see
\cite{Reshetikhin-Turaev:invariants,Kirby-Melvin}) that for $m,n\ge 0$,
we have
\begin{equation*}
  \langle \V_m, \V_n\rangle  = [(m+1)(n+1)].
\end{equation*}
Therefore, $\langle ,\rangle $ restricts to a bilinear form
\begin{equation*}
  \langle ,\rangle \zzzcolon \modR _{\Zvv}\times \modR _{\Zvv}\rightarrow \Zvv.
\end{equation*}
We also need the induced bilinear form
\begin{equation*}
  \langle ,\rangle \zzzcolon \modR _{\modQ (v)}\times \modR _{\modQ (v)}\rightarrow \modQ (v).
\end{equation*}

\subsubsection{The elements $S_n$}
\label{sec:elements-checkn}

Let $\modS $ denote the $\Zqq$-subalgebra of $\modR _{\Zqq}$ generated by
$\V_1^2$.  We have $\modS =\Span_{\Zqq}\{\V_{2i}\zzzvert i\ge 0\}$.

For $n\ge 0$, set
\begin{equation*}
  S_n
  =\prod_{i=1}^{n}(\V_1^2-(v^i+v^{-i})^2)
  =\prod_{i=1}^{n}(\V_2-(q^i+1+q^{-i}))
  \in \modS .
\end{equation*}
(The elements corresponding to $S_n$ in the Kauffman bracket skein
module of solid torus have been introduced in \cite{H:iias2000}.)
Clearly, we have
\begin{equation}
  \label{e90}
  \xi (S_n)=\sigma _n
\end{equation}
for $n\ge 0$.  Observe that $S_n$ is a monic polynomial in $\V_1^2$ of
degree $n$.  Therefore, $\modS $ is spanned over $\Zqq$ by the elements
$S_n$, $n\ge 0$.

\begin{proposition}
  \label{r23}
  If $m,n\ge 0$, then we have
  \begin{gather}
    \label{e21} \langle P_m, S_n\rangle  = \delta _{m,n}\BB{2m+1}{2m}.
  \end{gather}
  (As a consequence, the map $\langle ,\rangle \zzzcolon \modR _{\Zvv}\times \modS \rightarrow \Zvv$ is nondegenerate.)
\end{proposition}

\begin{proof}
  For $k\ge 0$, set $\V'_k=\V_k/[k+1]$.  The map
  \begin{equation*}
    \langle -,\V'_k\rangle \zzzcolon \modR _{\Zvv}\rightarrow \Zvv,\quad  x\mapsto \langle x,\V'_k\rangle 
  \end{equation*}
  is a $\Zvv$-algebra homomorphism.

  We will prove
  \begin{equation}
    \label{eq:37} \langle P_m, \V_{2n}\rangle  =[2n+1]\{n+m\}_{2m}.
  \end{equation}
  For $p\ge 0$, we have
  \begin{equation*}
    \langle \V_1-v^{2p+1}-v^{-2p-1},\V'_{2n}\rangle 
    =v^{2n+1}+v^{-2n-1}-v^{2p+1}-v^{-2p-1}=\{n-p\}\{n+p+1\}.
  \end{equation*}
  Since $\langle -,\V'_{2n}\rangle $ is a $\Zvv$-algebra homomorphism, we have
  \begin{equation*}
    \langle P_m,\V'_{2n}\rangle 
    =\prod_{p=0}^{m-1}\langle \V_1 - v^{2p+1}-v^{-2p-1},\V'_{2n}\rangle 
    =\prod_{p=0}^{m-1}\{n-p\}\{n+p+1\}=\{n+m\}_{2m}.
  \end{equation*}
  Therefore, we have \eqref{eq:37}.  Similarly, we can prove
  \begin{equation}
    \label{eq:40} \langle \V_m,S_n\rangle  = \{m+n+1\}_{2n+1}/\{1\}.
  \end{equation}

  By \eqref{eq:37}, $\langle P_m,\V_{2n}\rangle =0$ when $0\le n<m$.  Therefore, we
  have $\langle P_m,S_n\rangle =0$ if $0\le n<m$.  By \eqref{eq:40},
  $\langle \V_m,S_n\rangle =0$ when $0\le m<n$.  Hence, we have $\langle P_m,S_n\rangle =0$ when
  $0\le m<n$.  By \eqref{eq:37},
  \begin{equation*}
    \langle P_m,S_m\rangle =\langle P_m,\V_{2m}\rangle =[2m+1]\BB{2m}{2m}=\BB{2m+1}{2m},
  \end{equation*}
  since $S_m-\V_{2m}$ is a linear combination of $\V_{2i}$ for
  $i=0,1,\ldots,m-1$.  This completes the proof.
\end{proof}

\begin{remark}
  \label{thm:91}
  One can prove \eqref{eq:100} as follows.  Assume
  $\V_n=\sum_{j=0}^na_{n,j}P_j$, where $a_{n,j}\in \Zvv$ for $j=0,\ldots,n$.
  Applying $\langle -,S_i\rangle $ to the both sides, we obtain
  $\langle \V_n,S_i\rangle =\sum_{j=0}^n a_{n,j}\langle P_j,S_i\rangle $.  By \eqref{eq:40}
  and \eqref{e21}, we have
  $\{n+i+1\}_{2i+1}/\{1\}=a_{n,i}\BB{2i+1}{2i}$, hence
  $a_{n,i}=\bb{n+i+1}{2i+1}$.
\end{remark}

\subsubsection{Proof of Proposition \ref{r22}}
\label{sec:proof-prop-refr22-1}

For each $n\ge 0$ we have $\xi (S_n) = \sigma _n$.  Therefore, by \eqref{e90}, \eqref{eq:85}
and \eqref{e21}, we have
\begin{equation*}
  \trq^{P_m}(\sigma _n)
  =\trq^{P_m}(\xi (S_n))
  =\langle P_m,S_n\rangle 
  =\delta _{m,n}\BB{2m+1}{2m},
  \end{equation*}
hence the assertion.

\section{Remarks on knot invariants}
\label{sec7}

In this section we discuss some consequences of the results in Section
\ref{sec6}.  The reader may skip this section in the first
reading, since it is not necessary for the proof of the existence of
the invariants $J_M$ for integral homology spheres.

Throughout this section, let $T$ be a bottom knot and let $K=\cl(T)$
be the closure of $T$.  By abuse of notation, we set
$J_K=J_T\in Z(\tUqv)$.

\subsection{The two-variable colored Jones invariant}
\label{sec:color-jones-polyn-2}

Part of the following overlaps \cite[Section 3.3]{H:rims2001}, where
some results are stated without proofs.

It follows from the results in the previous sections that we have
\begin{equation*}
  J_K=J_T=\sum_{n\ge 0}J_K(P''_n)\sigma _n\in Z(\tUqv)
  \cong\varprojlim_{k}\modZ [q,q^{-1}][C^2]/(\sigma _k).
\end{equation*}
Thus $J_K$ may be regarded as an invariant with two variables $q$ and
$C$.  It is useful to introduce variables $t$ and $\alpha $ satisfying
\begin{equation*}
  \alpha ^2=C^2-4=t+t^{-1}-2.
\end{equation*}
Then $Z(\tUqv)$ may be identified with
\begin{equation*}
  \Lambda := \varprojlim_{k}\modZ [q,q^{-1}][t+t^{-1}]/(\sigma _k) \cong
  \varprojlim_{k}\modZ [q,q^{-1}][\alpha ^2]/(\sigma _k),
\end{equation*}
where
\begin{equation*}
  \sigma _k
  =\prod_{i=1}^k(t+t^{-1}-q^i-q^{-i})
  =\prod_{i=1}^k(\alpha ^2-q^i-q^{-i}+2)
\end{equation*}

\begin{remark}
  \label{r27}
  $\Lambda $ can be naturally regarded as a subring of
  \begin{equation*}
    \tilde\Lambda :=\varprojlim_{k}\modZ [q,q^{-1},t,t^{-1}]/(\sigma _k)
    \cong\varprojlim_k\modZ [q,q^{-1},t,t^{-1}]/
      \Bigl(\prod_{-k\le i\le k,\ i\neq0}(t-q^i)\Bigr).
  \end{equation*}
  In fact, $\Lambda $ consists of the elements of $\tilde\Lambda $ which are
  invariant under the involutive continuous ring automorphism of
  $\tilde\Lambda $ defined by $f(q,t)\mapsto f(q,t^{-1})$.
\end{remark}

By abuse of notation, we write
\begin{equation*}
  J_K(t,q)=J_K\in \Lambda .
\end{equation*}

For $i\in \modZ $, let
\begin{equation*}
  \theta _i\zzzcolon \modZ [q,q^{-1}][t+t^{-1}] \rightarrow  \modZ [q,q^{-1}]
\end{equation*}
be the unique $\modZ [q,q^{-1}]$-algebra homomorphism such that
$\theta _i(t+t^{-1})=q^i+q^{-i}$.  We have
\begin{equation}
  \label{e80}
  \theta _i(\sigma _k)=\prod_{j=1}^k(q^i+q^{-i}-q^j-q^{-j}).
\end{equation}

If $i\neq0$, then $\theta _i$ induces a homomorphism
\begin{equation*}
  \hat\theta _i\zzzcolon \Lambda  \rightarrow  \modZ [q,q^{-1}],
\end{equation*}
since it follows from \eqref{e80} that $\theta _i(\sigma _j)=0$ for all
$j\ge |i|$.  Since we have $\theta _i=\theta _{-i}$, it follows that
$\hat\theta _i=\hat\theta _{-i}$.  By \eqref{e89}, it follows that for $i>0$
\begin{equation}
  \label{e51}
  \theta _i(J_K)=J_K(q^i,q) = J_K(\V'_{i-1}),
\end{equation}
where $\V'_{i-1}=\V_{i-1}/[i]$.  Thus $J_K(q^i,q)\in \modZ [q,q^{-1}]$
is the $(i-1)$st normalized colored Jones polynomial.  Hence, it may
be natural to call $J_K(t,q)\in \Lambda $ the {\em two-variable colored Jones
invariant} of $K$.  For the comparison with the well-known,
two-variable power series expansion of the colored Jones polynomials,
known as the Melvin-Morton expansion, see Section
\ref{sec:melv-mort-expans} below.

Now, consider the case $i=0$.  $\theta _0$ induces a continuous ring
homomorphism
\begin{equation*}
  \hat\theta _0\zzzcolon \Lambda  \rightarrow  \Zqh, \quad f(t,q)\mapsto f(1,q),
\end{equation*}
since it follows from \eqref{e80} that
$\theta _0(\sigma _k)=(-1)^k(\{k\}!^2)\in \modZ [q,q^{-1}]$ is divisible by
${(q)_k}^2$.  For $\zeta \in \calZ$ with $r=\ord(\zeta )$, by \eqref{e51}, we
have
\begin{equation*}
  \mods _\zeta (J_K(1,q))=J_K(1,\zeta )=J_K(\zeta ^r,\zeta )=\mods _\zeta (J_K(q^r,q))
  =\mods _\zeta (J_K(\V'_{r-1})).
\end{equation*}
The right-hand side is known as the {\em Kashaev invariant} of $K$ at
$\zeta $ \cite{Kashaev,Murakami-Murakami}.  Hence, $J_K(1,q)\in \Zqh$ may be
regarded as a {\em unified Kashaev invariant}.  Essentially the same
observation has been made in \cite[Section 3.3]{H:rims2001} (without
proof).  Using a different method, Huynh and Le \cite{Huynh-Le} have
given a proof of the existence of the unified Kashaev invariant in the
above sense.

\subsection{Melvin-Morton expansions and Rozansky's rationality
  theorem}
\label{sec:melv-mort-expans}

Rozansky \cite{Rozansky3} proved that the colored Jones polynomials of
a knot $K$ can be repackaged into one invariant
$J^\MMR_K\in \modZ [[q-1,\alpha ^2]]$, which is an integer-coefficient version of
the Melvin-Morton expansion \cite{Melvin-Morton}.  (In fact, Rozansky
proved much stronger results there, see below.)

Let
\begin{equation*}
  \gamma \zzzcolon \Lambda \rightarrow \modZ [[q-1,\alpha ^2]]
\end{equation*}
be the ring homomorphism induced by $\id_{\modZ [q,q^{-1},\alpha ^2]}$.  This
is well-defined since $\sigma _{2k}\in (\alpha ^{2k},(q-1)^{2k})$ in
$\modZ [q,q^{-1},\alpha ^2]$ for $k\ge 0$.  As mentioned in \cite[Section
3.3]{H:rims2001}, $\gamma (J_K)\in \modZ [[q-1,\alpha ^2]]$ is equal to the Rozansky's
integral version of the Melvin-Morton expansion.

The first half of the following proposition was mentioned in
\cite{H:rims2001} without proof, and implies that $J_K$ and $J^\MMR_K$
have the same strength.

\begin{proposition}
  \label{r72}
  The map $\gamma $ is injective and non-surjective.
\end{proposition}

\begin{proof}
  Let us prove injectivity.
  Let $f=\sum_{n=0}^\infty f_n(q)\sigma _n\in \Lambda $, where $f_n(q)\in \modZ [q,q^{-1}]$ for
  $n\ge 0$, and suppose that $\gamma (f)=0$.  We can express each $\sigma _n$ as
  \begin{equation*}
    \sigma _n= \sum_{k=0}^n g_{n,k}(q)(2-q-q^{-1})^{2n-2k}\alpha ^{2k},
  \end{equation*}
  where $g_{n,k}(q)\in \modZ [q,q^{-1}]$.  Hence, we obtain
  \begin{equation*}
    \gamma (f)=\sum_{k=0}^\infty 
    \Bigl(\sum_{n=k}^\infty f_n(q)g_{n,k}(q)(2-q-q^{-1})^{2n-2k}\Bigr)\alpha ^{2k}.
  \end{equation*}
  Since $\gamma (f)=0$, we have
  \begin{equation}
    \label{e82}
    \sum_{n=k}^\infty f_n(q)g_{n,k}(q)(2-q-q^{-1})^{2n-2k}=0
  \end{equation}
  for all $k\ge 0$.

  One can show that $g_{n,k}(q)$ is not divisible by $q-1$.  (Indeed,
  one can use induction on $n$ to show that $g_{n,k}(1)\in \modZ $ is
  positive.)  Using this fact and \eqref{e82}, one can prove using
  induction that for any $n,d\ge 0$, $f_n(q)$ is divisible by $(q-1)^d$.
  Hence $f_n(q)=0$ for all $n$.

  Non-surjectivity of $\gamma $ can be verified in several ways.  For
  example, let $u\zzzcolon \modZ [[q-1,\alpha ^2]]\rightarrow \modZ [[q-1]]$ be the continuous
  $\modZ [[q-1]]$-algebra homomorphism such that $u(\alpha ^2)=0$.  Then we
  have
  \begin{equation*}
    u\gamma (\Lambda )=\modZ [q,q^{-1}]\subsetneq\modZ [[q-1]]=u(\modZ [[q-1,\alpha ^2]]).
  \end{equation*}
  It follows that $\gamma $ is not surjective.
\end{proof}

Rozansky \cite{Rozansky3} also proved that the coefficient of
$(q-1)^n$ of the integral Melvin-Morton expansion of a knot $K$ is a
rational function in $\alpha ^2$ with denominator a power of the Alexander
polynomial $\Delta _K(t)\in \modZ [t+t^{-1}]$ of $K$.  (This result generalizes a
conjecture of Melvin and Morton \cite{Melvin-Morton}, proved in
\cite{Bar-Natan-Garoufalidis}.)  More precisely, Rozansky proved that
$J^\MMR_K$ has the following power series expansion
\begin{equation*}
  J_K^\MMR = \sum_{n\ge 0}\frac{P_{K,n}(t)}{\Delta _K(t)^{2n+1}}(q-1)^n
  \in \modZ [t+t^{-1},\frac{1}{\Delta _K(t)}][[q-1]],
\end{equation*}
where $P_{K,n}(t)\in \modZ [t+t^{-1}]$ for $n\ge 0$.  The $P_{K,n}(t)$ are
uniquely determined by $K$ and $n$.  Thus $\gamma (J_K)$ is contained in a
much smaller ring than $\modZ [[q-1,\alpha ^2]]$.  Melvin and Morton's
conjecture (essentially) states that
\begin{equation*}
  J_K^\MMR\bigl|_{q=1}=\frac{1}{\Delta _K(t)}.
\end{equation*}
Thus we have $P_{K,0}(t)=1$.

The following conjecture (supported by some computer calculations)
generalizes Rozansky's rationality theorem.

\begin{conjecture}
  \label{r44}
  Let $K$ be a knot and let $d\zzzcolon \modN \rightarrow \{0,1,2,\ldots \}$ be a function which
  vanishes for all but finitely many elements of $\modN $.  Set
  \begin{equation*}
    \Delta ^d_K(t) = \prod_{r\in \modN ,d(r)>0}\Delta _K(t^r)^{2d(r)-1}\in \modZ [t+t^{-1}].
  \end{equation*}
  Then we have
  \begin{equation*}
    \Delta ^d_K(t)J_K \in \modZ [t+t^{-1},q,q^{-1}] + \Phi _d(q)\Lambda .
  \end{equation*}
\end{conjecture}

Conjecture \ref{r44} implies that $J_K^\MMR$ is contained in the image
of the natural (injective) homomorphism
  \begin{equation*}
    \varprojlim_k \modZ [q,q^{-1},t+t^{-1},\frac{1}{\Delta _K(t^i)}\quad  (i\ge 1)]/((q)_k)
    \rightarrow \modZ [[q-1,\alpha ^2]].
  \end{equation*}

Conjecture \ref{r44} also implies the following conjecture, in which
the case $\zeta =1$ is Rozansky's rationality theorem.

\begin{conjecture}
  \label{r58}
  Let $\zeta \in \calZ$ with $r=\ord(\zeta )$.  Let
  \begin{equation*}
    \gamma _\zeta \zzzcolon \Lambda \rightarrow \modZ [\zeta ][[\alpha ^2,q-\zeta ]]
  \end{equation*}
  be the homomorphism induced by
  $\modZ [q,q^{-1},\alpha ^2]\subset \modZ [\zeta ][q,q^{-1},\alpha ^2]$.    Define
  $P_{K,\zeta ,k}(t)\in \modZ [\zeta ][[\alpha ^2]]$ for $k\ge 0$ by
  \begin{equation*}
    \gamma _\zeta (J_K)=\sum_{k\ge 0}\frac{P_{K,\zeta ,k}(t)}{\Delta _K(t^r)^{2k+1}}(q-\zeta )^k.
  \end{equation*}
  Then we have $P_{K,\zeta ,k}(t)\in \modZ [\zeta ][t+t^{-1}]$.
\end{conjecture}

Using the injectivity of $\gamma =\gamma _1$, one can prove that the
homomorphism $\gamma _\zeta $ is injective.

Conjecture \ref{r58} implies that
\begin{equation*}
  J_K(t,\zeta )=\frac{P_{K,\zeta ,0}(t)}{\Delta _K(t^r)}
  \in \frac{1}{\Delta _K(t^r)}\modZ [\zeta ][t+t^{-1}].
\end{equation*}

Some computer calculations support the following conjecture.

\begin{conjecture}
  \label{r64}
  We have
  \begin{equation*}
    J_K(t,-1)=\frac{\Delta _K(-t)}{\Delta _K(t^2)}.
  \end{equation*}
  Thus we have $P_{K,-1,0}(t)=\Delta _K(-t)$ in the notation of Conjecture
  \ref{r58}.
\end{conjecture}

Presumably, $P_{K,\zeta _r,0}(t)$ with $\zeta _r=\exp\frac{2\pi \sqrt{-1}}{r}$
is the $r$th Akutsu-Deguchi-Ohtsuki invariant
\cite{Akutsu-Deguchi-Ohtsuki}, after suitable change of variables.
(Recall that the second Akutsu-Deguchi-Ohtsuki invariant of a knot is
the Alexander polynomial \cite{MurakamiJ}.)  If this is true, then it
follows that the Akutsu-Deguchi-Ohtsuki invariants of a knot $K$
determine $J_K$, hence also the colored Jones polynomials $J_K(\V_m)$,
$m\ge 0$.

\section{Algebraically-split links}
\label{sec8}

In this section, we prove an integrality result for the colored Jones
polynomials of algebraically-split, $0$-framed links.  We also give
some conjectures for boundary links in the last subsection.

\subsection{The algebra $\modP $ and its completion $\hat\modP $}
\label{sec:algebra--its-1}

For $n\ge 0$, set
\begin{gather*}
  P'_n = \frac{P_n}{\{n\}!}\in \modR _{\modQ (v)},\\
  \tP'_n = v^{-\hf n(n-1)}P'_n\in \modR _{\modQ (v)}.
\end{gather*}
We have $P'_0=1$.  Define a $\Zqq$-submodule $\modP $ of $\modR _{\modQ (v)}$ by
\begin{equation*}
  \modP  = \Span_{\Zqq}\{\tP'_n\zzzvert n\ge 0\}.
\end{equation*}

\begin{lemma}
  \label{r26}
  $\modP $ is a $\Zqq$-subalgebra of $\modR _{\modQ (v)}$.
\end{lemma}

\begin{proof}
  By induction, we obtain
  \begin{equation*}
    P_mP_n = \sum_{i=0}^{\min(m,n)}
    \frac{\BB mi\BB ni\BB{m+n}i}{\{i\}!}P_{m+n-i}.
  \end{equation*}
  Hence,
  \begin{equation}
    \label{eq:69}
    P'_mP'_n = \sum_{i=0}^{\min(m,n)}
    \frac{\{m+n\}!}{\{i\}!\{m-i\}!\{n-i\}!}P'_{m+n-i}.
  \end{equation}
  The coefficient in the right-hand side is contained in $v^d\Zqq$,
  where
  \begin{math}
    d\equiv\binom{m+n+1}2-\binom{i+1}2-\binom{m-i+1}2-\binom{n-i+1}2
    \equiv \binom{i+1}2\pmod2.
  \end{math}
  Since
  \begin{math}
    \binom m2+\binom n2 -\binom{m+n-i}2
      \equiv \binom{i+1}2\pmod2,
  \end{math}
  we have $v^{\binom m2}P'_mv^{\binom n2}P'_n\in \modP $.  Therefore, we have
  $\modP \modP \subset \modP $.  Since we have $1=P'_0\in \modP $, we have the assertion.
\end{proof}

For $k\ge 0$, set
\begin{equation*}
  \modP _k = \Span_{\Zqq}\{\tP'_n\zzzvert n\ge k\}\subset \modP .
\end{equation*}
We have
\begin{equation*}
  \modP =\modP _0\supset \modP _1\supset \modP _2\supset \cdots.
\end{equation*}
It follows from \eqref{eq:69} that each $\modP _k$ is an ideal in $\modP $.
Set
\begin{equation*}
  \hP = \varprojlim_{k\ge 0} \modP /\modP _k,
\end{equation*}
which inherits a $\Zqq$-algebra structure from $\modP $.  Each element
$x\in \hP$ can be uniquely expressed as an infinite sum $x=\sum_{k=0}^\infty 
x_k\tP'_k$, where $x_k\in \Zqq$ for $k\ge 0$.  Since
$\bigcap_{k\ge 0}\modP _k=\{0\}$, we may regard $\modP $ as a subalgebra of
$\hP$.

\subsection{Integrality for algebraically-split, $0$-framed links}
\label{sec:maps-jlhphphaq}
The following theorem is proved in Section
\ref{sec:proof-theorem-refr29}.

\begin{theorem}
  \label{r29}
  Let $L$ be an $m$-component, algebraically-split, $0$-framed link.
  For $i=1,\ldots,m$, let $x_i\in \modP _{k_i}$, $k_i\ge 0$.  Then we have
  \begin{equation*}
    J_L(x_1,\ldots,x_m) \in \frac{\BBq{2k+1}{k+1}}{\{1\}_q}\Zqq,
  \end{equation*}
  where $k=\max(k_1,\ldots,k_m)$.
\end{theorem}

Theorem \ref{r29} immediately implies the following.

\begin{corollary}
  \label{r87}
  Let $L$ be an $m$-component, algebraically-split, $0$-framed link.
  Then the $\modQ (v)$-multilinear map
  $J_L\zzzcolon \modR _{\modQ (v)}\times \cdots\times \modR _{\modQ (v)}\rightarrow \modQ (v)$ restricts to a
  $\Zqq$-multilinear map
  \begin{equation*}
    J_L\zzzcolon \modP \times \cdots\times \modP \rightarrow \Zqq,
  \end{equation*}
  which induces a $\Zqq$-multilinear map
  \begin{equation*}
    J_L\zzzcolon \hP\times \cdots\times \hP\rightarrow \Zqh.
  \end{equation*}
\end{corollary}

\begin{remark}
  \label{r30}
  It follows from \eqref{eq:100} that for an $m$-component,
  algebraically-split, $0$-framed link $L$, the colored Jones
  polynomial $J_L(\V_{n_1},\ldots,\V_{n_m})$, $n_1,\ldots,n_m\ge 0$, can be
  expressed as a linear combination of the
  $J_L(P'_{k_1},\ldots,P'_{k_m})$, $k_1,\ldots,k_m\ge 0$ as follows:
  \begin{equation}
    \label{e37}
    \begin{split}
      &J_L(\V_{n_1},\ldots,\V_{n_m})
      =\sum_{k_1=0}^{n_1}\cdots\sum_{k_m=0}^{n_m}
      \prod_{i=1}^m\Bigl(\bb{n_i+k_i+1}{2k_i+1}\{k_i\}!\Bigr)
      J_L(P'_{k_1},\ldots,P'_{k_m}).
    \end{split}
  \end{equation}
  Here the sums may be replaced with infinite sums $\sum_{k_i=0}^\infty $.
\end{remark}

\subsection{Proof of Theorem \ref{r29}}
\label{sec:proof-theorem-refr29}

The following lemma is proved in the next subsection.

\begin{lemma}
  \label{thm:63}
  If $x\in \Uqv$ and $y\in \modP $, then we have $\trq^y(x)\in \Zqq$.
\end{lemma}

Using Lemma \ref{thm:63}, we obtain the following generalization of
\eqref{e34}.

\begin{proposition}
  \label{r25}
  If $T$ is an $m$-component, algebraically-split, $0$-framed bottom
  tangle, and if $x_2,\ldots,x_m\in \modP $, then we have
  \begin{equation*}
    (1\otimes \trq^{x_2}\otimes \cdots\otimes \trq^{x_m})(J_T)\in Z(\tUqv).
  \end{equation*}
\end{proposition}

\begin{proof}
  By Lemma \ref{thm:63},
  \begin{equation}
    \label{e26}
    (1\otimes \trq^{x_2}\otimes \cdots\otimes \trq^{x_m})((\Uqv)^{\otimes m})\subset \Uqv.
  \end{equation}
  Let $N\ge 0$ be such that
  \begin{equation}
    x_2,\ldots,x_m
    \in \Span_{\modQ (v)}\{P'_0,P'_1,\ldots,P'_N\}
    =\Span_{\modQ (v)}\{\V_0,\V_1,\ldots,\V_N\}.
  \end{equation}
  If $p>N$ then we have $\trq^{x_i}(\modF _p(\Uqv))=0$ for $i=2,\ldots,m$,
  since $e^p$ acts as zero on $\V_0,\V_1,\ldots,\V_N$.  Therefore, we have
  \begin{equation}
    \label{e27}
    (1\otimes \trq^{x_2}\otimes \cdots\otimes \trq^{x_m})(\modF _p((\Uqv)^{\otimes m}))
    \subset \modF _p(\Uqv).
  \end{equation}
  Now, \eqref{e26}, \eqref{e27} and Theorem \ref{thm:8} imply the
  assertion.
\end{proof}

\begin{proof}[Proof of Theorem \ref{r29}]
  We may assume $k_1=k=\max(k_1,\ldots,k_m)$ without loss of generality.
  (Otherwise, change the order of components of $L$.)  Choose a bottom
  tangle $T\in \BT_m$ such that $\cl(T)=L$.  By Proposition \ref{r25},
  \begin{equation*}
    J_L(x_1,\ldots,x_m)
    =(\trq^{x_1}\otimes \cdots\otimes \trq^{x_m})(J_T)
    =\trq^{x_1}(y),
  \end{equation*}
  where we set $y=(1\otimes \trq^{x_2}\otimes \cdots\otimes \trq^{x_m})(J_T)$, which is
  contained $Z(\tUqv)$ by Proposition \ref{r25}.  By the arguments in
  Section \ref{sec:case-n=1},
  \begin{equation*}
    y=\sum_{p\ge 0} a_p\sigma _p,
  \end{equation*}
  where $a_p\in \Zqq$ for $p\ge 0$.  Moreover, we have
  \begin{equation*}
    x_1=\sum_{p=k}^N b_p\tP'_p,
  \end{equation*}
  where $b_p\in \Zqq$ and $N\ge k$.  Therefore, we have
  \begin{equation*}
    \trq^{x_1}(y)
    =\sum_{p=k}^N a_pb_p\trq^{\tP'_p}(\sigma _p).
  \end{equation*}
  By Proposition \ref{r22}, we have
  $\trq^{\tP'_p}(\sigma _p)\in \frac{\BBq{2p+1}{p+1}}{\{1\}_q}\Zqq$, hence
  the assertion.
\end{proof}

\subsection{Proof of Lemma \ref{thm:63}}
\label{sec:integr-quant-trac}

To prove Proposition \ref{thm:63}, we need some lemmas.

\begin{lemma}
  \label{lem:2}
  If $y,y'\in \modR _{\modQ (v)}$ and $x\in U_h$, then we have
  \begin{equation*}
    \trq^{yy'}(x) = \sum \trq^{y}(x_{(1)})\trq^{y'}(x_{(2)}),
  \end{equation*}
  where $\Delta (x)=\sum x_{(1)}\otimes x_{(2)}$.
\end{lemma}

\begin{proof}
  The assertion follows from the well-known identity for any two
  finite-dimensional representations $W,W'$
  \begin{equation*}
    \trq^{W\otimes W'}(x)=(\trq^W\otimes \trq^{W'})\Delta (x).
  \end{equation*}
\end{proof}

\begin{lemma}
  \label{lem:6}
  (1) If $n,l,l'\ge 0$, $l\neq l'$ and $j\in \modZ $, then we have
  $\trq^{P_n}(\tF^{(l)}K^{2j}e^{l'})=0$.

  (2) If $0\le n<l$ and $j\in \modZ $, then we have
      $\trq^{P_n}(\tF^{(l)}K^{2j}e^l)=0$.

  (3) If $0\le l\le n$ and $j\in \modZ $, then we have
  \begin{equation}
    \label{eq:46}
    \trq^{P_n}(\tF^{(l)}K^{2j}e^l)
    =v^n q^{-nj+2lj+l^2-ln} \{n\}_q!\{n-l\}_q!
    \bbq{j+n-1}{n-l} \bbq{j-1}{n-l}.
  \end{equation}
\end{lemma}

\begin{proof}
  The assertion (1) follows from the fact that $l\neq l'$ implies that
  $\tF^{(l)}K^{2j}e^{l'}$ acts nilpotently on each $\V_n$, $n\ge 0$.

  The assertion (2) follows from the fact that $P_n$ is a linear
  combination of $\V_0,\ldots,\V_n$, and that if $0\le n<l$, then $e^l$
  acts as $0$ on $\V_0,\ldots,\V_n$.

  We prove (3) by induction on $n$.  The case $n=l=0$ follows from
  \begin{equation*}
    \trq^{P_0}(K^{2j})=\epsilon (K^{2j})=1.
  \end{equation*}
  Let $n\ge 1$.  Set $p_n=\V_1-v^{2n-1}-v^{-2n+1}$, so that
  $P_n=P_{n-1}p_n$.  Using Lemma \ref{lem:2}, we have
  \begin{equation}
    \label{eq:63}
    \trq^{P_n}(\tF^{(l)}K^{2j}e^l)
    =\trq^{P_{n-1}p_n}(\tF^{(l)}K^{2j}e^l)
    =\sum\trq^{P_{n-1}}(x_{(1)})\trq^{p_n}(x_{(2)}),
  \end{equation}
  where
  \begin{equation*}
    \begin{split}
      \sum x_{(1)}\otimes x_{(2)}
      &=\Delta (\tF^{(l)}K^{2j}e^l)\\
      &=\sum_{r=0}^l\sum_{s=0}^l
      \bbq{l}{s}
      (\tF^{(l-r)}K^rK^{2j}e^{l-s}K^s\otimes \tF^{(r)}K^{2j}e^s).
    \end{split}
  \end{equation*}
  Since $p_n$ is a linear combination of $\V_0$ and $\V_1$, we have
  $\trq^{p_n}(F^{(r)}K^{2j}e^s)=0$ unless $(r,s)=(0,0)$ or $(1,1)$.
  Therefore, the right-hand side of \eqref{eq:63} equals
  \begin{equation}
    \label{eq:36}
      \trq^{P_{n-1}}(\tF^{(l)}K^{2j}e^{l})\trq^{p_n}(K^{2j})
       + q^{-l+1}[l]_q
      \trq^{P_{n-1}}(\tF^{(l-1)}K^{2j+2}e^{l-1})\trq^{p_n}(\tF^{(1)}K^{2j}e).
  \end{equation}
  By straightforward computation,
  \begin{gather*}
    \trq^{p_n}(K^{2j}) = v^{-2j+1}\{j-n\}_q\{j+n-1\}_q,\\
    \trq^{p_n}(\tF^{(1)}K^{2j}e) = v^{2j+1}\{1\}_q.
  \end{gather*}
  Using these identities and the inductive hypothesis, we can verify
  that \eqref{eq:36} equals the right-hand side of \eqref{eq:46}.
  This completes the proof.
\end{proof}

\begin{proof}[Proof of Lemma \ref{thm:63}]
  It suffices to prove that if $n\ge 0$ and $x\in \Uqv$, then we have
  $\trq^{P'_n}(x)\in v^{\hf n(n-1)}\Zqq$.  We may assume without loss of
  generality that $x=\tF^{(l)}K^{2j}e^k$ with $j\in \modZ $, $k,l\ge 0$.  By
  Lemma \ref{lem:6}, $\trq^{P_n}(x)\in v^n\{n\}_q!\Zqq=v^{\hf
  n(n-1)}\{n\}!\Zqq$.  Therefore, we have the assertion.
\end{proof}

\subsection{Remarks on boundary links and boundary bottom tangles}
\label{sec:remarks-bound-links}

In this subsection, we give some remarks on the case of boundary
links.  The reader may skip this subsection since it is not necessary
in the rest of the paper.

Let $\bUqv$ denote the $\modZ [q,q^{-1}]$-subalgebra of $\Uqv$ generated
by the elements $K^2$, $K^{-2}$, $e$ and $f:=(q-1)FK$.  Define a
decreasing filtration $\modF _k(\bUqv)$, $k\ge 0$, by
\begin{equation*}
  \modF _k(\bUqv)=\bUqv\cap \modF _p(\Uqv)=\bUqv\cap \modF _p(\Uq).
\end{equation*}
Let $\bUqvt$ denote the completion in $U_h$ of $\bUqv$ with respect to
the filtration $\modF _k(\bUqv)$, i.e., we set
\begin{equation*}
  \bUqvt=\operatorname{Image}\Bigl(\varprojlim_k\bUqv/\modF _k(\bUqv)\rightarrow U_h\Bigr).
\end{equation*}
Clearly, $\bUqvt$ is a $\modZ [q,q^{-1}]$-subalgebra of $\tUqv$.  In the
similar way to the case of $\tUqvn$, we can define the completed
tensor products $\bUqvt\;^{\tOn}=\bUqvt\;\tO\cdots\tO\bUqvt\subset \Uhn$.
Let $\Inv(\bUqvtn)$ denote the invariant part of $\bUqvtn$.  We have
\begin{equation*}
  \Inv(\bUqvtx1)=\Inv(\bUqvt)=Z(\bUqvt)=Z(\tUqv).
\end{equation*}
Here one can easily verify the last identity using \cite[Section
9]{H:uqsl2}.

A bottom tangle $T=T_1\cup \cdots\cup T_n$ with $0$-framings in a cube $[0,1]^3$
is called a {\em boundary bottom tangle} if there are $n$ disjoint,
connected, oriented surface $F_1,\ldots,F_n$ in $[0,1]^3$ such that for
each $i=1,\ldots,n$ the boundary $\partial F_i$ of $F_i$ is the union of $T_i$
and the line segment in $\partial [0,1]^3$ bounded by the two endpoints of
$T_i$.  (See also \cite[Section 9.3]{H:universal}.)

The following conjecture generalizes the case $n=1$ of Theorem
\ref{thm:8}.  (Note that a bottom knot is a $1$-component, boundary
bottom tangle.)

\begin{conjecture}
  \label{r81}
  Let $T\in \BT_n$ be an $n$-component, boundary bottom tangle.
  Then we have $J_T\in \Inv(\bUqvtn)$.
\end{conjecture}

Here we give an outline of a possible proof of Conjecture \ref{r81}.
We have not worked out all the details yet, but the following idea
seems useful.  The main tool is \cite[Theorem 9.9]{H:universal}, which
reduces the problem to showing that $Y^{\otimes g}(J_{T'})\in \bUqvtx g$ for
any bottom tangle $T'\in \BT_{2g}$ (with any linking matrix), $g\ge 0$.
Here $Y\zzzcolon U_h\ho U_h\rightarrow U_h$ is the commutator morphism for the
transmutation of $U_h$ (see \cite[Section 9.3]{H:universal} for the
definition).

Conjecture \ref{r81} implies the following variant of Theorem
\ref{r29}.

\begin{conjecture}
  \label{r20}
  Let $L$ be an $m$-component, $0$-framed, boundary link.  For
  $i=1,\ldots,m$, let $x_i\in \modP _{k_i}$, $k_i\ge 0$.  Then we have
  \begin{equation}
    \label{e88}
    J_L(x_1,\ldots,x_m)
    \in \frac{\BBq{2k_1+1}{k_1+1}}{\{1\}_q}I_{k_2}I_{k_3}\cdots I_{k_m},
  \end{equation}
  Here, for $k\ge 0$, $I_k$ is the ideal in $\Zqq$ generated by the
  elements $\{k-l\}_q!\{l\}_q!$, $l=0,\ldots,k$.  (One can permute
  $k_1,\ldots,k_m$ in \eqref{e88}.)
\end{conjecture}

\begin{proof}[Proof that Conjecture \ref{r81} implies Conjecture \ref{r20}]
  We can express $L$ as the closure of an $m$-component, boundary
  bottom tangle $T$.  We have
  \begin{equation*}
    y:=(\id\otimes \trq^{x_2}\otimes \cdots\otimes \trq^{x_m})(J_T)
    \in I_{k_2}I_{k_3}\cdots I_{k_m}(\bUqv)\tilde{},
  \end{equation*}
  since $J_T$ is an infinite sum of elements of $(\bUqv)^{\otimes m}$ by
  Conjecture \ref{r81} and we have $\trq^{x_i}(\bUqv)\subset I_{k_i}$ by
  Lemma \ref{lem:6}.  By \eqref{e1}, we have
  \begin{equation*}
    J_L(x_1,\ldots,x_m)=(\trq^{x_1}\otimes \cdots\otimes \trq^{x_m})(J_T)=\trq^{x_1}(y).
  \end{equation*}
  The rest of the proof is similar to that for Theorem \ref{r29}.
\end{proof}

\begin{remark}
  \label{r21}
  One can prove partial results of Conjecture \ref{r20} by focusing on
  the divisibility of link invariants by powers of $q-1$ or $q+1$.  It
  is not difficult to prove using the theory of Goussarov-Vassiliev
  finite type link invariants that if $L$ is an $m$-component,
  boundary link, then $J_L(P_{k_1},\ldots,P_{k_m})$ is divisible by
  $(q-1)^{2(k_1+\cdots+k_m)}$ (hence $J_L(P'_{k_1},\ldots,P'_{k_m})$ is
  divisible by $(q-1)^{k_1+\cdots+k_m}$).  One can also prove that
  $J_L(P_{k_1},\ldots,P_{k_m})$ is divisible by $(q+1)^{k_1+\cdots+k_m}$.  The
  latter assertion follows from the fact that the Jones polynomial
  (with colors $\V_1$) of $n$-component, boundary link is divisible by
  $(q+1)^n$, which we plan to prove in a paper in preparation
  \cite{H:spanning} using skein theory.
\end{remark}

\section{Twists}
\label{sec9}

In this section, we introduce an element $\omega \in \hP$ satisfying a
``twisting property''.

\subsection{The twist element $\omega \in \hP$}
\label{sec:twist-element-hr-1}

Define two elements $\omega _+,\omega _-\in \hP$ by
\begin{equation*}
  \omega _\pm =\sum_{n=0}^\infty (\pm 1)^n v^{\pm \hf n(n+3)}P'_n.
\end{equation*}

Set $\modS _{\modQ (v)}=\modS \otimes _{\Zqq}\modQ (v)\subset \modR _{\modQ (v)}$.  By Proposition \ref{r23},
the bilinear map $\langle ,\rangle \zzzcolon \modP \times \modS _{\modQ (v)}\rightarrow \modQ (v)$ induces a bilinear map
\begin{equation*}
  \langle ,\rangle \zzzcolon \hP\times \modS _{\modQ (v)}\rightarrow \modQ (v).
\end{equation*}

\begin{proposition}
  \label{thm:38}
  For each $x\in \modS _{\modQ (v)}$, we have
  \begin{gather}
    \label{eq:67}
    \langle \omega _\pm ,x\rangle =J_{U_\pm }(x),
  \end{gather}
  where $U_\pm $ is a $\pm 1$-framed unknot.
\end{proposition}

\begin{proof}
  We may assume $x=\V'_{2k}=\V_{2k}/[2k+1]$, $k\ge 0$, without loss of
  generality.

  By \eqref{eq:37}, we have
  \begin{equation}
    \label{e40}
      \langle \omega _\pm ,\V'_{2k}\rangle 
      =\sum_{n\ge 0}(\pm 1)^nv^{\pm \hf n(n+3)}\langle P'_n,\V'_{2k}\rangle 
      =\sum_{n=0}^k (\pm 1)^nv^{\pm \hf n(n+3)}\frac{\{k+n\}_{2n}}{\{n\}!}.
  \end{equation}

  Since $\V_{2k}$ is an irreducible $U_h$-module, we have
  $\modr ^{\mp1}\modv _i^{2k}=J_{U_\pm }(\V'_{2k})\modv _i^{2k}$ for all
  $i=0,\ldots,2k$.
  (Recall that $\modr \in U_h$ denotes the ribbon element, and $\modv ^{2k}_i$
  for $i=0,\ldots,2k$ denote the basis elements of the irreducible left
  $U_h$-module $\V_{2k}$.)  Consider the case $i=k$.
  By straightforward calculations using \eqref{eq:65} and
  \eqref{eq:66}, we obtain
  \begin{equation*}
    \modr ^{\mp1}\modv ^{2k}_k
    =\sum_{n=0}^k(\pm 1)^nv^{\pm \hf n(n+3)}\frac{\BB{k+n}{2n}}{\{n\}!}\modv ^{2k}_k.
  \end{equation*}
  Hence, $J_{U_\pm }(\V'_{2k})$ is equal to the right-hand side of \eqref{e40}.
  Hence, we have \eqref{eq:67} for $x=\V'_{2k}$.
\end{proof}

\begin{proposition}
  \label{thm:42}
  $\omega _+$ and $\omega _-$ are inverse to each other in the algebra $\hP$.
\end{proposition}

\begin{proof}
  A direct proof using \eqref{eq:69} is possible.  Here we give
  another proof.

  By Proposition \ref{thm:38} and the fact that
  $\langle -,\V'_{2p}\rangle \zzzcolon \modS _{\modQ (v)}\rightarrow \modQ (v)$ is a $\modQ (v)$-algebra
  homomorphism, we have
  \begin{equation*}
    \langle \omega _+\omega _-,\V'_{2p}\rangle 
    =\langle \omega _+,\V'_{2p}\rangle \langle \omega _-,\V'_{2p}\rangle 
    =q^{p(p+1)}q^{-p(p+1)}
    =1=\langle 1,\V'_{2p}\rangle 
  \end{equation*}
  for each $p\ge 0$.  Hence, we have $\langle \omega _+\omega _-,x\rangle =\langle 1,x\rangle $ for all
  $x\in \modS _{\modQ (v)}$.  Since the map
  \begin{equation*}
    \hP \rightarrow \Fun(\modS _{\modQ (v)},\modQ (v)),\quad x\mapsto (y\mapsto \langle x,y\rangle ),
  \end{equation*}
  is injective (by Proposition \ref{r23}), it follows that
  $\omega _+\omega _-=1$.  (Here $\Fun(\modS _{\modQ (v)},\modQ (v))$ denotes the set of
  functions on $\modS _{\modQ (v)}$ with values in $\modQ (v)$.)
\end{proof}

Set $\omega =\omega _+$.  We have $\omega _-=\omega ^{-1}$.

\begin{remark}
  \label{r28}
  The twist element $\omega $ was announced first in \cite{H:iias2000} in a
  formulation using skein theory, and then in \cite{H:rims2001} in the
  present form (without proofs).  Later, Masbaum \cite{Masbaum} gave a
  proof using skein theory.
\end{remark}

\subsection{Twisting theorem}
\label{sec:twisting-theorems}

\begin{theorem}
  \label{thm:99}
  Let $L_1\cup \cdots\cup L_m\cup K$ be an $(m+1)$-component, algebraically-split,
  $0$-framed link such that $K$ is an unknot.  Set $L=L_1\cup \cdots\cup L_m$,
  and let $L_{(K,\pm 1)}$ denote the framed link in $S^3$ obtained from
  $L$ by $\pm 1$-framed surgery along $K$.  Then, for any
  $x_1,\ldots,x_m\in \hP$, we have
  \begin{equation}
    \label{eq:81}
    J_{L\cup K}(x_1,\ldots,x_m,\omega ^{\mp1})=J_{L_{(K,\pm 1)}}(x_1,\ldots,x_m).
  \end{equation}
\end{theorem}

\begin{proof}
  By Theorem \ref{r29}, it suffices to prove \eqref{eq:81} for
  $x_1,\ldots,x_m\in \modP $.  In fact, we will show \eqref{eq:81} for
  $x_1,\ldots,x_m\in \modR _{\modQ (v)}$.  (Well-definedness of the left-hand side
  of \eqref{eq:81} in this case follows from the argument below.)
  Apply the standard ``fusing'' argument to the strands running
  through the unknot $K$, see Figure
  \ref{X01}.{\def\tI{=\sum}\FI{X01}{A fusing identity.}}  In the
  left-hand side, the strands running through $K$ are paired in such a
  way that each pair contains two segments from the same component of
  $L$ in the opposite orientations.  By fusing these strands, we
  obtain a linear combination of trivalent colored graphs, where each
  term has exactly one edge $e$ running through $K$, colored by an
  even number.  Hence, there is $y\in \modS _{\modQ (v)}$ such that
  \begin{equation*}
    J_{L\cup K}(x_1,\ldots,x_m,\omega ^{\mp1}) =J_H(y,\omega ^{\mp1}) =\langle \omega ^{\mp1},y\rangle ,
  \end{equation*}
  where $H$ denotes the Hopf link with $0$-framings.  By Proposition
  \ref{thm:38}, the right-hand side is equal to $J_{U_\mp}(y)$, where
  $U_\mp$ denotes the $\mp1$-framed unknot.  Since we have
  \begin{equation*}
    J_{U_\mp}(y)=J_{L_{(K,\pm 1)}}(x_1,\ldots,x_m),
  \end{equation*}
  we have \eqref{eq:81}.
\end{proof}

\section{The invariant $J_M$}
\label{sec10}

In this section, we construct an invariant $J_M\in \Zqh$ of an integral
homology sphere $M$.

\subsection{Admissible framed links and the Hoste moves}
\label{sec:hoste-moves}

A framed link $L$ in $S^3$ is said to be {\em admissible} if $L$ is
algebraically split and $\pm 1$-framed.  Surgery on $S^3$ along an
admissible framed link yields an integral homology sphere.
Conversely, each integral homology sphere $M$ can be obtained as the
result of surgery on $S^3$ along an admissible framed link.

In the definition of the Reshetikhin-Turaev invariant
\cite{Reshetikhin-Turaev:invariants} and its variants, one uses
Kirby's theorem \cite{Kirby} or its variant by Fenn and Rourke
\cite{Fenn-Rourke}.  Fenn and Rourke's theorem states that two framed
links in $S^3$ yields orientation-preserving homeomorphic results of
surgery if and only if they are related by a sequence of isotopies and
{\em Fenn-Rourke moves}, where a Fenn-Rourke move is either surgery on
a unknotted, $\pm 1$-framed component, or the inverse operation.  See
Figure \ref{fig:FR}.  \FIG{FR}{A Fenn-Rourke move.}{height=30mm}

We use the following version of Fenn-Rourke's theorem, which was
conjectured by Hoste \cite{Hoste}.

\begin{theorem}[\cite{H:kirby1}]
  \label{r32}
  Two admissible framed links $L$ and $L'$ in $S^3$ yield
  orientation-preserving homeomorphic results of surgery if and only
  if $L$ and $L'$ are related by a sequence of isotopies and {\em
  Hoste moves}.  Here a Hoste move is defined to be a Fenn-Rourke move
  between two admissible framed links.
\end{theorem}

\subsection{Definition of $J_M$}
\label{sec:defin-unif-sl2}

Let $M$ be an integral homology sphere.  Choose an admissible framed
link $L=L_1\cup \cdots\cup L_m$ in $S^3$ such that $S^3_L\cong M$.  For
$i=1,\ldots,m$, let $f_i=\pm 1$ denote the framing of $L_i$.  Let
$L^0=L^0_1\cup \cdots\cup L^0_m$ denote the link $L$ with $0$ framings.  Set
\begin{equation}
  \label{eq:42}
  J_M= J_{L^0}(\omega ^{-f_1},\omega ^{-f_2},\ldots,\omega ^{-f_m})\in \Zqh.
\end{equation}

\begin{theorem}
  \label{thm:41}
  For an integral homology sphere $M$, the element $J_M\in \Zqh$ defined
  in \eqref{eq:42} does not depend on the choice of $L$.  Hence, the
  correspondence $M\mapsto J_M$ defines an invariant of integral
  homology spheres with values in $\Zqh$.
\end{theorem}

\begin{proof}
  Let $I(L)$ denote the right-hand side of \eqref{eq:42}.  Clearly,
  $I(L)$ is an isotopy invariant of admissible framed links.  By
  Theorem \ref{r32}, it suffices to prove that $I(L)$ is invariant
  under Hoste moves.  Observe that $I(L)$ is invariant under
  permutation of components.  Therefore, it suffices to show that if the
  last component $L_m$ of $L$ is an unknot, then
  \begin{equation*}
    J_{L^0}(\omega ^{-f_1},\ldots,\omega ^{-f_{m-1}},\omega ^{-f_m})
    =J_{(L^0\setminus L^0_m)_{L_m}}(\omega ^{-f_1},\ldots,\omega ^{-f_{m-1}}),
  \end{equation*}
  where $(L^0\setminus L^0_m)_{L_m}$ is obtained from $L^0\setminus L^0_m$ by
  performing surgery along $L_m$ with framing $f_m$.  This identity
  follows from Theorem \ref{thm:99}.  Hence the assertion.
\end{proof}

In Section \ref{sec:an-alternative-proof}, we give an alternative
proof of Theorem \ref{thm:41}.

\subsection{Some remarks}
\label{sec:some-remarks}

If $M$ and $L$ are as above, then we have the following formula for $J_M$:
\begin{equation}
  \label{e2}
  J_M
  = \sum_{k_1,\ldots,k_m=0}^\infty 
  \Bigl(\prod_{i=1}^m(-f_i)^{k_i}v^{\frac{-f_i}{2} k_i(k_i+3)}\Bigr)
  J_{L^0}(P'_{k_1},\ldots,P'_{k_m}).
\end{equation}
In particular, if $L$ is a knot (i.e., $m=1$), then we have
\begin{equation}
  \label{eq:132}
    J_M
    = J_{L^0}(\omega ^{\mp1})
    = \sum_{k=0}^\infty (\mp1)^kv^{\mp\hf k(k+3)}
    \frac{\BB{2k+1}{k+1} }{\{1\}} J_{L^0}(P''_k).
\end{equation}

In \eqref{e2}, the term for $k_1=\cdots=k_m=0$ is
\begin{equation*}
  J_{L^0}(P'_0,\ldots,P'_0)=J_{L^0}(1,\ldots,1)=1.
\end{equation*}
It follows from Theorem \ref{r29} that the other terms are divisible
by $(q^2-1)(q^3-1)/(q-1)$ in $\Zqq$.  Hence, we have the following
result, which we improve in Section \ref{sec:divisibility-results}.

\begin{lemma}
  \label{r47}
  Let $M$ be an integral homology sphere.  Then $J_M-1$ is divisible
  by $(q^2-1)(q^3-1)/(q-1)$ in $\Zqh$.
\end{lemma}

\begin{remark}
  \label{r24}It follows from Theorem \ref{r29} that
  \begin{equation*}
    J_M\mod I_k\in \Zqh/I_k\cong \modZ [q]/I_k,
  \end{equation*}
  where $I_k=\frac{\BBq{2k+1}{k+1}}{\{1\}_q}\Zqq$, is a finite linear
  combination of $J_{L^0}(P'_{k_1},\ldots,P'_{k_m})$ for
  $k_1,\ldots,k_m\in \{0,1,\ldots,k-1\}$, hence is a finite linear combination
  of the colored Jones polynomial $J_{L^0}(\V_{k_1},\ldots,\V_{k_m})$ for
  $k_1,\ldots,k_m\in \{0,1,\ldots,k-1\}$.  In particular, for $k=2$, $J_M\mod
  I_2$ is a linear combination of the Jones polynomials (with colors
  $\V_1$) of all the sublinks of $L^0$.
\end{remark}

\section{Specializations at roots of unity}
\label{sec11}

In this section, we prove the following.

\begin{theorem}
  \label{thm:65}
  Let $M$ be an integral homology sphere and let $\zeta $ be a primitive
  $r$th root of unity, $r\in \modN $.  Then we have
  \begin{equation}
    \label{eq:130}
    \mods _\zeta (J_M) = \tau _\zeta (M),
  \end{equation}
  where $\tau _\zeta (M)$ denotes the $sl_2$ WRT
  invariant of $M$ at $\zeta $.
\end{theorem}

\subsection{Definition of $\tau _\zeta (M)$}
\label{sec:definition-m}

We briefly recall the definition of the $sl_2$
WRT invariant $\tau _\zeta (M)$ of a closed
$3$-manifold $M$ at  $\zeta \in \calZ$.  For the details, see
\cite{Reshetikhin-Turaev:invariants,Kirby-Melvin}.

If $\zeta =1$, then we set $\tau _\zeta (M)=1$ for any closed $3$-manifold $M$.

Let $r\ge 2$ and let $\zeta ^{1/4}$ be a primitive $4r$th root of unity.
Thus $\zeta $ is a primitive $r$th root of unity.  (The following
construction depends not only on $\zeta $ but also $\zeta ^{1/2}$ or
$\zeta ^{1/4}$, but we suppress them from the notation.  For integral
homology spheres, $\tau _\zeta (M)$ depends only on $\zeta $ and $M$, not on
$\zeta ^{1/4}$.)

Let
\begin{equation*}
  A_r=\modZ [q^{1/4},q^{-1/4}]_{(\Phi _{4r}(q^{1/4}))},
\end{equation*}
denote the localization of the ring $\modZ [q^{1/4},q^{-1/4}]$ by the
principal ideal generated by the $4r$th cyclotomic polynomial
$\Phi _{4r}(q^{1/4})$ in $q^{1/4}$ (which is equal to
$\Phi _{2r}(q^{1/2})$).  We have
\begin{equation*}
  A_r=\Bigl\{\frac{f(q^{1/4})}{g(q^{1/4})}\in \modQ (q^{1/4})
    \;\Bigl|\;f(q^{1/4}),g(q^{1/4})\in \modZ [q^{1/4}],\; g(\zeta ^{1/4})\neq0 \Bigr\}.
\end{equation*}
We extend $\mods _\zeta \zzzcolon \Zqq\rightarrow \modZ [\zeta ]$ to the $\modQ $-algebra homomorphism
\begin{equation*}
  \mods _\zeta \zzzcolon A_r\rightarrow \modQ (\zeta ^{1/4})
\end{equation*}
defined by $\mods _\zeta (f(q^{1/4}))=f(\zeta ^{1/4})$.

Set
\begin{equation*}
  \Omega _r = \sum_{i=0}^{r-2}[i+1]\V_i\in \modR _{\Zvv}.
\end{equation*}
Let $L$ be an $m$-component link in $S^3$.  Set
\begin{gather*}
  I_r(L)=J_L(\Omega _r,\ldots,\Omega _r) \in \modZ [q^{1/4},q^{-1/4}],\\
  I_\zeta (L)=\mods _\zeta (I_r(L)) \in \modZ [\zeta ^{1/4},\zeta ^{-1/4}].
\end{gather*}
For an unknot $U_\pm $ of framing $\pm 1$, we have $I_\zeta (U_\pm )\neq0$.  Set
\begin{equation}
  \label{e43}
  \tau _\zeta (L)=
  \frac{I_\zeta (L)}{I_\zeta (U_+)^{\sigma _+(L)}I_\zeta (U_-)^{\sigma _-(L)}}
  \in \modQ (\zeta ^{1/4}),
\end{equation}
where $U_\pm $ denotes an unknot with framing $\pm 1$, and $\sigma _+(L)$
(resp. $\sigma _-(L)$) denotes the number of positive (resp. negative)
eigenvalues of the linking matrix of $L$.  It is known that $\tau _\zeta (L)$
is invariant under Kirby moves, i.e., stabilization and handle slides.
Hence, by setting $\tau _\zeta (M)=\tau _\zeta (L)$ for a $3$-manifold $M$ with
$S^3_L\cong M$, we obtain a $3$-manifold invariant $\tau _\zeta (M)$.

\subsection{WRT invariants of integral homology spheres}
\label{sec:wrt-invar-integr}

Define $t\zzzcolon \modR _{\modQ (q^{1/4})}\rightarrow \modR _{\modQ (q^{1/4})}$ to be the unique
$\modQ (q^{1/4})$-module isomorphism, called the {\em twist operator},
characterized by
\begin{equation*}
  J_U(t^{\pm 1}(x))= J_{U_\pm }(x)\quad \text{for $x\in \modR _{\modQ (q^{1/4})}$},
\end{equation*}
where $U$ is the $0$-framed unknot.  We have
$t(\V_n)=q^{\frac{n(n+2)}{4}}\V_n$.

Set
\begin{equation*}
  \Omega _{r,\pm 1} = t^{\mp1}(\Omega _r)/I_r(U_\mp).
\end{equation*}
An important property of $\Omega _{r,\pm 1}$ is the following:
\begin{equation}
  \label{e46}
  \mods _\zeta (\langle \Omega _{r,\pm 1},x\rangle ) = \mods _\zeta (J_{U_\pm }(x))
\end{equation}
for $x\in \modR _{\modQ (q^{1/4})}$.

For admissible framed links, we can simplify \eqref{e43} as follows.
Let $L=L_1\cup \cdots\cup L_m$ be an admissible framed link with framings
$f_1,\ldots,f_m\in \{\pm 1\}$.  Let $L^0=L^0_1\cup \cdots\cup L^0_m$ denote $L$ with $0$
framings.  Then we have
\begin{equation*}
  \tau _\zeta (M)=\tau _\zeta (L)=\mods _\zeta (J_{L^0}(\Omega _{r,-f_1},\ldots,\Omega _{r,-f_m})).
\end{equation*}

\subsection{Proof of Theorem \ref{thm:65}}
\label{sec:proof-theor-refth}

First consider the case $r=1$, $\zeta =1$.  We have $\tau _1(M)=1$ by the
definition.  We also have $\mods _1(J_M)=1$ by Lemma \ref{r47}.
Hence, we have \eqref{eq:130} in this case.

Let $r\ge 2$.  Let $L$ be as in the last subsection.  By the definition
of $J_M$, we have
\begin{equation}
  \label{e49}
  \mods _\zeta (J_M)= \mods _\zeta (J_{L^0}(\omega ^{-f_1},\omega ^{-f_2},\ldots,\omega ^{-f_m}))
\end{equation}

Set $d=\floor{\frac{r-2}{2}}$.  It follows from Theorem~\ref{r29} that
if $x_1,\ldots,x_m\in \modP $ and we have $x_i\in \modP _{d+1}$ for some
$i\in \{1,\ldots,m\}$, then we have $\mods _\zeta (J_L(x_1,\ldots,x_m))=0$.  Hence, by
\eqref{e49}, we have
\begin{equation*}
  \mods _\zeta (J_M)= \mods _\zeta (J_{L^0}(\omega _{r,-f_1},\omega _{r,-f_2},\ldots,\omega _{r,-f_m})),
\end{equation*}
where $\omega _{r,\pm 1}$ is the truncation of $\omega ^{\pm 1}$ given by
\begin{equation*}
  \omega _{r,\pm 1} =\sum_{n=0}^d(\pm 1)^n v^{\pm \hf n(n+3)}P'_n\in \modP .
\end{equation*}
Note that $\omega _{r,\pm 1} \in \Span_{A_r}\{P_0,P_1,\ldots,P_{d}\}\subset \modR _{A_r}$.
We also have $\Omega _{r,\pm 1}\in \modR _{A_r}$.

\begin{lemma}
  \label{r34}
  Let $L=L_1\cup \cdots\cup L_m\cup K$ be an $(m+1)$-component,
  algebraically-split, $0$-framed link.  Then, for
  $x_1,\ldots,x_m\in \modR _{A_r}$, we have
  \begin{equation}
    \label{e52}
    \mods _\zeta (J_{L\cup K}(x_1,\ldots,x_m,\omega _{r,\pm 1})) =
    \mods _\zeta (J_{L\cup K}(x_1,\ldots,x_m,\Omega _{r,\pm 1})).
  \end{equation}
\end{lemma}

\begin{proof}
  If $K$ is an unknot, then by Theorem \ref{thm:99} and the definition
  of $\omega _{r,\pm 1}$ we have
  \begin{equation}
    \label{e53}
    \mods _\zeta (J_{L\cup K}(x_1,\ldots,x_m,\omega _{r,\pm 1})) =
    \mods _\zeta (J_{L_{(K,\mp1)}}(x_1,\ldots,x_m)),
  \end{equation}
  where $L_{(K,\mp1)}$ is the result from $L$ by surgery along $K$
  with framing $\mp1$.  By \eqref{e46} and the standard fusing
  argument, we also have
  \begin{equation*}
    \mods _\zeta (J_{L\cup K}(x_1,\ldots,x_m,\Omega _{r,\pm 1})) = \mods _\zeta (J_{L_{(K,\mp1)}}(x_1,\ldots,x_m)).
  \end{equation*}
  Hence, we have \eqref{e52}.

  Now, we consider the general case.  A crossing change of two strands
  in $K$ can be done by a Hoste move.  Applying this move finitely
  many times, we get a framed link $L\cup K'\cup E_1\cup \cdots\cup E_t$, $t\ge 0$,
  such that all the linking numbers between two components of
  $L\cup K'\cup E_1\cup \cdots\cup E_t$ are zero, $K'$ is a $0$-framed unknot,
  $E_1\cup \cdots\cup E_t$ is an admissible framed unlink with framings
  $g_1,\ldots,g_t\in \{\pm 1\}$, and $(L\cup K')_{E_1\cup \cdots\cup E_t}\cong
  L\cup K$.  Using \eqref{e53} iteratively, we obtain
  \begin{equation*}
    \begin{split}
      &\mods _\zeta (J_{L\cup K}(x_1,\ldots,x_m,y))\\
      &\quad =\mods _\zeta (J_{L\cup K'\cup E_1^0\cup \cdots\cup E_m^0}
      (x_1,\ldots,x_m,y,\omega _{r,-g_1},\ldots,\omega _{r,-g_t})) ,
    \end{split}
  \end{equation*}
  where $E_i^0$ is $E_i$ with $0$-framing, and $y$ is either
  $\omega _{r,\pm 1}$ or $\Omega _{r,\pm 1}$.  Since $K'$ is an unknot, the
  assertion reduces to the special case proved above.
\end{proof}

\subsection{An alternative proof of Theorem \ref{thm:41}}
\label{sec:an-alternative-proof}

Here we give an alternative proof of Theorem \ref{thm:41} which does
not use Theorem \ref{r32}, but uses the existence of the
invariant $\tau _\zeta (M)$.

Let $M$ and $L$ be as in Section \ref{sec:defin-unif-sl2}.  Let
$(J_M)_L$ denote the right-hand side of \eqref{eq:42}.  Let $L'$ be
another link such that $(S^3)_{L'}=M$.  We have to show that
$(J_M)_L=(J_M)_{L'}$.  By Theorem \ref{thm:65}, we have
$\mods _\zeta ((J_M)_L)=\tau _\zeta (M)=\mods _\zeta ((J_M)_{L'})$ for each $\zeta \in \calZ$.
(Note that we do not need the invariance of $J_M$ in the proof of
Theorem \ref{thm:65}.)  Hence, Proposition \ref{r85} implies
$(J_M)_L=(J_M)_{L'}$.  This completes the proof.

Similar idea has been used in \cite{Beliakova-Blanchet-Le,Le4}.  See
Section \ref{sec:rati-homol-spher}.

\section{Some properties of $J_M$}
\label{sec12}

In this section, we give some properties of $J_M$ and some
applications.

Throughout this section, $M$ denotes an integral homology sphere.

\subsection{Connected sum and orientation reversal}
\label{sec:conn-sum-orient}

\begin{proposition}
  \label{thm:43}
  (1) The invariant $J_M$ is multiplicative under connected sum, i.e.,
  for any two integral homology spheres $M$ and $M'$ we have
  $J_{M\sharp M'}=J_M J_{M'}$, where $M\sharp M'$ denotes the connected sum of
  $M$ and $M'$.  We also have $J_{S^3}=1$.

  (2) If $-M$ denotes the mirror image of an integral homology sphere,
  then we have $J_{-M}=\overline{J_M}$.  Here $\overline{J_M}\in \Zqh$
  is the conjugate of $J_M$, i.e., the image of $J_M$ under the
  involutive ring automorphism $\Zqh\rightarrow \Zqh$, $f(q)\mapsto f(q^{-1})$.
\end{proposition}

\begin{proof}
  The assertions follows from the corresponding, well-known properties
  of $\tau _\zeta (M)$
  \begin{gather*}
    \tau _\zeta (M\sharp M')=\tau _\zeta (M)\tau _\zeta (M),\quad
    \tau _\zeta (S^3)=1,\quad
    \tau _\zeta (-M)=\tau _{\zeta ^{-1}}(M)
  \end{gather*}
  for any $\zeta \in \calZ$, and from Proposition \ref{r85}.  (Alternatively,
  one can use the similar properties for the Ohtsuki series $\tau ^O(M)$
  and the injectivity of $\iota _1\zzzcolon \Zqh\rightarrow \modZ [[q-1]]$.  One can also prove
  the assertions directly using properties of the colored Jones
  polynomials.)
\end{proof}

\subsection{On determination of $J_M$ by infinitely many WRT invariants}
\label{sec:determination-jm-wrt}

In this subsection, we prove the following.

\begin{proposition}
  \label{r56}
  If $\calZ'\subset \calZ$ has no limit point, then the ring homomorphism
  $\mods _{\calZ'}$ in \eqref{e61} is not injective.
\end{proposition}

Here, recall from Section \ref{sec:ring-zqh-cyclotomic} that a ``limit
point'' of a subset $\calZ'\subset \calZ$ is an element $\xi \in \calZ$ such
that there are infinitely many $\zeta \in \calZ'$ adjacent to $\xi $.  (Recall
that $\xi $ and $\zeta $ are adjacent if $\ord(\zeta \xi ^{-1})$ is a prime
power.)  Any finite subset $\calZ'\subset \calZ$ has no limit point.  There
are many infinite subset $\calZ'\subset \calZ$ without limit points, for
example $\{6^i\zzzvert i\ge 0\}$.

Proposition \ref{r56} shows that Conjecture 6.1 in our previous paper
\cite{H:cyclotomic} is false.  It follows that we can not prove, using
only properties of $\Zqh$, that the $\tau _\zeta (M)$ for any infinitely many
$\zeta \in \calZ$ would determine $J_M$ and hence $\tau (M)$.  It should be
remarked that it is still open whether the WRT invariants at any
infinitely many roots of unity determine $J_M$ or not.

Proposition \ref{r56} follows easily from Proposition \ref{r55} below.
Two elements $m,n\in \modN $ are said to be {\em adjacent} to each other
(written $m\Leftrightarrow n$) if $m/n=p^e$ with $p$ prime and $e\in \modZ $.  For a
subset $S\subset \modN $, a limit point of $S$ is defined to be an element
$m\in \modN $ which is adjacent to infinitely many elements of $S$.  Note
that if $\zeta ,\xi \in \calZ$ are adjacent, then $\ord(\zeta )$ and $\ord(\xi )$
are adjacent.  (One can define a topology on $\modN $ where $S\subset \modN $ is
open if for each $m\in S$, $S$ contains all but finitely many elements
adjacent to $m$.  The map $\ord\zzzcolon \calZ\rightarrow \modN $, $\zeta \mapsto\ord(\zeta )$,
turns out to be a continuous map.)

\begin{proposition}
  \label{r55}
  If $S\subset \modN $ is a subset with no limit point, then the ring
  homomorphism
  \begin{equation}
    \label{e8}
    \mods _S\zzzcolon \Zqh\rightarrow \prod_{m\in S}\modZ [q]/(\Phi _m(q)),\quad
    f(q)\mapsto (f(q)\mod(\Phi _m(q)))_{m\in S},
  \end{equation}
  is not injective.
\end{proposition}

\begin{proof}
  Let $S=\{m_1,m_2,\ldots \}$ with $1\le m_1<m_2<\cdots$.
  Choose $a\in \modN \setminus S$.  We
  construct $f(q)\in \Zqh$ such that $f(q)\in (\Phi _{m_i}(q))$ for all
  $i\ge 1$ and $f(q)\not\in (\Phi _a(q))$.
  For this purpose, it suffices to construct a sequence
  $f_i(q)\in \modZ [q]$, $i\ge 1$, such that
  \begin{enumerate}
  \item[(1)]
    $f_i(q)\rightarrow 1$ in $\Zqh$ as $i\rightarrow \infty $ (i.e., for each $j$ there is
    $i$ such that $f_k(q)-1$ is divisible by $(q)_j$ for all $k\ge i$),
  \item[(2)]
    $f_i(q)\in (\Phi _{m_i}(q))$ for $i\ge 1$,
  \item[(3)]
    $f_i(q)\not\in (\Phi _a(q))$ for $i\ge 1$.
  \end{enumerate}
  Suppose there is such a sequence $f_i(q)$.  Set
  $f(q)=\prod_{i\ge 1}f_i(q)$, which is well defined by (1).  By (2) we
  have $f(q)\in (\Phi _{m_i}(q))$ for all $i\ge 1$.  (1) and (3) implies that
  $f(q)\not\in (\Phi _a(q))$, hence $f(q)\neq0$.

  Now we show that there is a sequence $f_i(q)$ as above.
  For $i\ge 0$, set
  \begin{equation*}
    b_i=\max(\{0\}\cup \{b\in \modN \zzzvert b<m_i;\; b\not\Leftrightarrow m_j\;\text{for all $j\ge i$}\}).
  \end{equation*}
  The assumption in the statement implies:
  \begin{itemize}
  \item[(a)] $0\le b_1\le b_2\le \cdots$ and $\lim_{i}b_i=\infty $,
  \item[(b)] if $1\le n\le b_i$ and $i\le j$, then $m_j$ and $n$ are not adjacent.
  \end{itemize}
  It follows from (b) that, for $i\ge 1$, there is $u_i(q)\in \modZ [q]$ such
  that $(q)_{b_i}u_q(q)\equiv1\mod{(\Phi _{m_i}(q))}$.  Set
  \begin{equation*}
    f_i(q)=
    \begin{cases}
      \Phi _{m_i}(q)& \text{if $b_i<a$}, \\
      1-(q)_{b_i}u_i(q)& \text{if $b_i\ge a$},\\
    \end{cases}
  \end{equation*}
  (2) and (3) can be easily verified.  (1) follows from (a).
\end{proof}

\begin{proof}[Proof of Proposition \ref{r56}]
  Set $S=\{\ord(\zeta )\zzzvert \zeta \in \calZ'\}\subset \modN $.  By assumption, $S$ has no
  limit point.  Hence, by Proposition \ref{r55}, $\mods _S$ is not
  injective.  It follows that $\mods _{\calZ'}$ is not injective, since it
  factors through $\mods _S$.
\end{proof}

\begin{remark}
  \label{r43}
  By slightly modifying the definition of $f_i(q)$ in the proof of
  Proposition \ref{r55}, we see that if $S\subset \modN $ has no limit point and
  if $d\zzzcolon S\rightarrow \{0,1,\ldots \}$ is any function, then the ring homomorphism
  \begin{equation*}
    g_{S,d}\zzzcolon \Zqh\rightarrow \prod_{m\in S}\modZ [q]/(\Phi _m(q)^{d(m)}),\quad
    f(q)\mapsto (f(q)\mod(\Phi _m(q)^{d(m)}))_{m\in S},
  \end{equation*}
  is not injective.  (It suffices to replace $f_i(q)$ with the
  $d(m_i)$th power of the original definition of $f_i(q)$.)  Using
  this, one can show that if $\calZ'\subset \calZ$ has no limit point and if
  $d\zzzcolon \calZ'\rightarrow \{0,1,\ldots \}$ is any function, then the ring homomorphism
  \begin{equation*}
    g_{\calZ',d}\zzzcolon \Zqh\rightarrow \prod_{\zeta \in \calZ'}\modZ [q]/((q-\zeta )^{d(m)}),\quad
    f(q)\mapsto (f(q)\mod((q-\zeta )^{d(\zeta )}))_{\xi \in \calZ'},
  \end{equation*}
  is not injective.
\end{remark}

\begin{remark}
  \label{r74}
  In the situation of Proposition \ref{r55}, suppose moreover that the
  elements of $S$ are pairwise non-adjacent.  Then the homomorphism
  $\mods _S$ is surjective.  Since $S=\{6^n\zzzvert n\ge 0\}$ is such an example of
  $S$ and it is infinite, it follows that $\Zqh$ is {\em not}
  Noetherian.
\end{remark}

\subsection{Power series invariants}
\label{sec:power-seri-expans}

In this subsection, we study the power series expansions
$\iota _\zeta (J_M)\in \modZ [\zeta ][[q-\zeta ]]$, in particular the Ohtsuki series
$\tau ^O(M)=\iota _1(J_M)$ (see Theorem \ref{r52} below).

\subsubsection{Derivatives and power series expansions}
\label{sec:power-seri-expans-1}

Like polynomials and analytic functions, the elements of $\Zqh$ can be
differentiated arbitrarily many times.  In fact, the derivation
$\frac{d}{dq}\zzzcolon \modZ [q]\rightarrow \modZ [q]$ is continuous with respect to the
topology defined by the ideals $(q)_n\modZ [q]$, $n\ge 0$, since we have
\begin{equation}
  \frac{d}{dq}(q)_{2n}\in (q)_n\modZ [q].
\end{equation}
Therefore, $\frac d{dq}$ induces a continuous $\modZ $-linear derivation
\begin{equation}
  \frac{d}{dq}\zzzcolon \Zqh\rightarrow \Zqh.
\end{equation}
This implies that there are higher derivations
\begin{equation*}
  \frac{1}{k!}\frac{d^k}{dq^k}=\frac{1}{k!}
  (\frac{d}{dq})^k\zzzcolon \Zqh\rightarrow \Zqh.
\end{equation*}
for all $k\ge 0$.  (Recall that $(\frac{d}{dq})^k(\modZ [q])\subset k!\modZ [q]$.
This implies that $(\frac{d}{dq})^k(\Zqh)\subset k!\Zqh$.)
Suppose $f(q)\in \Zqh$ and $\zeta \in \calZ$.  Define the ``Taylor series at
$\zeta $'' of $f(q)$ by
\begin{equation*}
  T_\zeta (f(q))=\sum_{k\ge 0}\frac{1}{k!}\frac{d^kf}{dq^k}(\zeta )(q-\zeta )^k
  \in \modZ [\zeta ][[q-\zeta ]].
\end{equation*}
Then one can easily check $T_\zeta (f(q))=\iota _\zeta (f(q))$.

\subsubsection{The Ohtsuki series}
\label{sec:ohtsuki-series}
As mentioned in the introduction, according to Lawrence's conjecture
proved by Rozansky, the Ohtsuki series $\tau ^O(M)$ for an integral
homology sphere can be defined as the unique element
$\tau ^O(M)\in \modZ [[q-1]]$ such that for each root of unity $\zeta $ of odd
prime power order $r$, $\tau ^O(M)|_{q=\zeta }$ converges $p$-adically to
$\tau _\zeta (M)\in \modZ [\zeta ]$, i.e., we have
\begin{equation}
  \label{e59}
  \mods _\zeta (\tau ^O(M))=\tau _\zeta (M)
\end{equation}
in $\modZ _p[\zeta ]\cong\varprojlim_n\modZ [\zeta ]/(p^n)\cong\varprojlim_n\modZ [\zeta ]/((\zeta -1)^n)$,
where
\begin{equation*}
  \mods _\zeta \zzzcolon \modZ [[q-1]]\rightarrow \modZ _p[\zeta ]
\end{equation*}
denote the ring homomorphism induced by the evaluation map
\begin{equation*}
  \mods _\zeta \zzzcolon \modZ [q]\rightarrow \modZ [\zeta ],\quad  f(q)\mapsto f(\zeta ).
\end{equation*}
These properties characterize $\tau ^O(M)$.

\begin{theorem}
  \label{r52}
  For every integral homology sphere $M$, we have $\iota _1(J_M)=\tau ^O(M)$.
\end{theorem}

\begin{proof}
  For each root of unity $\zeta $ of prime power order $r=p^e$, we have
  the following commutative diagram
  \begin{equation*}
    \begin{CD}
      \Zqh @>\iota _1>> \modZ [[q-1]]\\
      @V\mods _\zeta VV @VV\mods _\zeta V\\
      \modZ [\zeta ] @>>i> \modZ _p[\zeta ].
    \end{CD}
  \end{equation*}
  Here $i$ is the inclusion.  By the commutativity of the diagram and
  by Theorem \ref{thm:65}, we have
  $\mods _\zeta \iota _1(J_M)=i\mods _\zeta (J_M)=i(\tau _\zeta (M))$.  Therefore, by \eqref{e59},
  we have $\mods _\zeta (\iota _1(J_M))=\mods _\zeta (\tau ^O(M))$.  The assertion follows from
  Lemma \ref{lem:8}.
\end{proof}

\begin{lemma}
  \label{lem:8}
  The map
  \begin{equation*}
    \modZ [[q-1]]\rightarrow \prod_{p:\text{odd prime}}\modZ _p[\zeta _p],\quad
    x\mapsto(\mods _{\zeta _p}(x))_p,
  \end{equation*}
  is injective.  Here $\zeta _p=\exp\frac{2\pi \sqrt{-1}}{p}$.
\end{lemma}

\begin{proof}
  Write $x=\sum_{n\ge 0}x_n(q-1)^n$.  We assume that $\mods _{\zeta _p}(x)=0$ for
  all $p$, and we see that $x_n=0$ by induction on $n$.  Let $n\ge 0$
  and suppose that $x_0=\cdots=x_{n-1}=0$.  Then, for every odd prime $p$,
  we have $\mods _{\zeta _p}(x)\equiv x_n(\zeta -1)^n\mod (\zeta -1)^{n+1}$.  Therefore,
  $x_n$ is divisible by $p$.  It follows that $x_n=0$.
\end{proof}

\subsubsection{Expansions at roots of unity}
\label{sec:expansions-at-roots}

For every root of unity $\zeta $, $\iota _\zeta (J_M)\in \modZ [\zeta ][[q-\zeta ]]$ may be
considered as a generalization of the Ohtsuki series
$\tau ^O(M)=\iota _1(J_M)$.  Define $\lambda _{\zeta ,k}(M)\in \modZ [\zeta ]$, $k\ge 0$, by
\begin{equation*}
  \iota _\zeta (J_M)=\sum_{k\ge 0}\lambda _{\zeta ,k}(M)(q-\zeta )^k.
\end{equation*}

\begin{question}
  \label{r41}
  Recall that the Ohtsuki series is considered to be related to the
  perturbative expansion of the Chern-Simons path integral.  Are there
  any ``physical'' interpretation of the other power series
  $\iota _\zeta (J_M)$?
\end{question}

Let us consider the case $\zeta =-1$.  We have
\begin{equation}
  \label{e75}
  \iota _{-1}(J_M)=\sum_{k\ge 0}\lambda _{-1,k}(M)(q+1)^k
  \in \modZ [[q+1]]
\end{equation}
where $\lambda _{-1,k}(M)\in \modZ $.  It follows from Lemma \ref{r47} that
$\lambda _{-1,0}(M)=1$.  Thus the first interesting coefficient is
$\lambda _{-1,1}(M)\in \modZ $, which is additive under connected sum, and
satisfies $\lambda _{-1,1}(-M)=-\lambda _{-1,1}(M)$ by Proposition \ref{thm:43}.
Some computations show that $\lambda _{-1,1}$ is nontrivial (i.e.,
non-vanishing for some $M$), and linearly independent over $\modQ $ to the
Casson invariant $\lambda (M)=\frac{1}{6}\lambda _1(M)$ of $M$.  (However, there
is a congruence relation between them, see Remark \ref{r63} below.)
See Remarks \ref{r61} and \ref{r63} for some other properties of
$\lambda _{-1,k}(M)$.  An interesting problem is to identify the topological
information carried by $\lambda _{-1,1}(M)$.

\subsubsection{Relations between two power series expansions at
  different roots of unity}
\label{sec:relat-betw-two}

Let $\zeta ,\xi \in \calZ$, where $\ord(\zeta /\xi )=p^e$ is a prime power.  We
consider how the two power series expansions
$\iota _\zeta (J_M)\in \modZ [\zeta ][[q-\zeta ]]$ and $\iota _\xi (J_M)\in \modZ [\xi ][[q-\xi ]]$ are
related.  Note that these two formal power series rings can be
regarded as subrings of
\begin{equation*}
  R_{\zeta ,\xi }:=\varprojlim_{i,j}\modZ [\zeta ,\xi ][q]/((q-\zeta )^i,(q-\xi )^j)\cong
  \modZ _p[\zeta ,\xi ][[q-\zeta ]]\cong\modZ _p[\zeta ,\xi ][[q-\xi ]].
\end{equation*}
In $R_{\zeta ,\xi }$ we have
\begin{equation*}
  \begin{split}
    \iota _\zeta (J_M)
    &=\sum_{k\ge 0}\lambda _{\zeta ,k}(q-\zeta )^k\\
    &=\sum_{k\ge 0}\lambda _{\zeta ,k}((q-\xi )+(\xi -\zeta ))^k\\
    &=\sum_{i\ge 0}\Bigl(\sum_{j\ge 0}\binom{i+j}{i}(\xi -\zeta )^j\lambda _{\zeta ,i+j}\Bigr)(q-\xi )^i\\
  \end{split}
\end{equation*}
Hence, we have
\begin{equation}
  \label{e65}
  \lambda _{\xi ,k}=\sum_{j\ge 0}\binom{k+j}{k}(\xi -\zeta )^j\lambda _{\zeta ,k+j}
\end{equation}
for $k\ge 0$.  This formula describes how the coefficients of the power
series expansion at $\zeta $ determines those at $\xi $.  Note that each
$\lambda _{\xi ,k}$ is determined modulo a given power of $p$ by finitely many
coefficients in the power series expansion at $\zeta $.

Setting $k=0$ in \eqref{e65}, we have
\begin{equation}
  \label{e66}
  \tau _\xi (M)=\sum_{j\ge 0}(\xi -\zeta )^j\lambda _{\zeta ,j}
\end{equation}
Note that if $\zeta =1$ and $\ord(\xi )$ is an odd prime power, then
\eqref{e66} implies that the Lawrence's $p$-adic convergence
conjecture mentioned in the introduction.  Thus \eqref{e66}
generalizes Lawrence's conjecture.

\eqref{e66} implies the following.

\begin{proposition}
  \label{r67}
  Let $\zeta ,\xi \in \calZ$ with $\ord(\zeta /\xi )$ a prime power.  Then we have
  \begin{equation*}
    \tau _\zeta (M)\equiv \tau _\xi (M)\pmod{(\xi -\zeta )}
  \end{equation*}
  in $\modZ [\zeta ,\xi ]$.
\end{proposition}

\begin{corollary}
  \label{r71}
  Let $\zeta \in \calZ$ with $\ord(\zeta )=dp^e$, where $d\in \{1,2,3,4,6\}$,
  $e\ge 0$, and $p$ is a prime.  Then we have $\tau _\zeta (M)\neq0$.
\end{corollary}

\begin{proof}
  It is well known that we have $\tau _\xi (M)\neq\pm 1$ for $\xi \in \calZ$ with
  $\ord(\xi )\in \{1,2,3,4,6\}$ (see Section
  \ref{sec:divisibility-results}).  This fact and Proposition
  \ref{r67} implies the assertion, since $(\xi -\zeta )\neq(1)$ in
  $\modZ [\zeta ,\xi ]$.
\end{proof}

The case $d=1$ and $p$ odd in Corollary \ref{r71} follows also from
Lawrence's conjecture proved by Rozansky, and presumably has been well
known.

\begin{conjecture}[Non-vanishing Conjecture]
  \label{r73}
  For any integral homology sphere $M$, we have $\tau _\zeta (M)\neq0$ for
  every root of unity $\zeta $.
\end{conjecture}

See Conjectures \ref{r33} and \ref{r31} for stronger statements.

\subsection{Divisibility properties and applications to power series
  invariants}
\label{sec:divisibility-jm-1}

\subsubsection{Divisibility in the ring $\Zqh$}
\label{sec:divis-ring-zqh}

Here, we make some basic observations about divisibility for elements
of $\Zqh$ by products of powers of cyclotomic polynomials, which we
will freely use in the rest of this section.

\begin{lemma}
  \label{r66}
  Let $f(q)\in \Zqh$, $\zeta \in \calZ$, $n=\ord\zeta $, and $k\ge 0$.  Then the
  following conditions are equivalent.
  \begin{enumerate}
  \item\label{item:9} $f(q)$ is divisible by $\Phi _n(q)^k$ in $\Zqh$.
  \item\label{item:10} $f(q)$ is divisible by $(q-\zeta )^k$ in
  $\Zqh\otimes _\modZ \modZ [\zeta ]$.
  \item\label{item:11} $\frac{d^if}{dq^i}(\zeta )=0$ for $i=0,1,\ldots,k-1$.
  \end{enumerate}
\end{lemma}

\begin{proof}
  It is easy to check that (\ref{item:10}) and (\ref{item:11}) are
  equivalent, and that (\ref{item:9}) implies (\ref{item:10}).  To prove
  (\ref{item:10}) implies (\ref{item:9}), use the fact that $f(q)$ is
  divisible by $(q-\zeta ')^k$ for each $\zeta '\in \calZ$ with $\ord(\zeta ')=n$
  (by Galois equivariance).  The details are left to the reader.
\end{proof}

Given $f(q)\in \Zqh$ with $f(q)\neq0$, let $\ord_n(f(q))\ge 0$ denote the
largest integer $l$ such that $f(q)$ is divisible by $\Phi _n(q)^l$.
(Existence of such integer follows from injectivity of the natural
homomorphism
\begin{equation*}
  \Zqh\rightarrow \varprojlim_{k}\modZ [q]/(\Phi _n(q)^k),
\end{equation*}
see \cite[Corollary 4.1]{H:cyclotomic}.)  If $f(q)=0$, then set
$\ord_n(f(q))=\infty $ for all $n\in \modN $.  It is not difficult to show that
if $d\zzzcolon \modN \rightarrow \{0,1,\ldots \}$ is a function such that $d(n)\le \ord_n(f(q))$
for all $n\in \modN $, and such that $d(n)=0$ for all but finitely many
$n\in \modN $, then $f(q)$ is divisible by $\prod_{n\in \modN }\Phi _n(q)^{d(n)}$ in
$\Zqh$.

A consequence of this fact is the following.
\begin{proposition}
  \label{r57}
  The ring $\Zqh$ is {\em not} a unique factorization domain.
\end{proposition}

\begin{proof}
  One can easily check that each cyclotomic polynomial $\Phi _n(q)$,
  $n\in \modN $, is a prime element in $\Zqh$.  (In fact, the principal
  ideal $(\Phi _n(q))$ in $\Zqh$ is a prime ideal.)

  Let $f(q)\in \Zqh$ be as defined in the proof of Proposition
  \ref{r55}.  Thus, there is a sequence $1\le m_1<m_2<\cdots$ such that
  $f(q)\in (\Phi _{m_i}(q))\subset \Zqh$ for all $i\ge 1$, where $\{m_1,m_2,\ldots \}$
  has no limit points (see the proof of Proposition \ref{r55}).  The
  observation above implies that $f(q)$ is divisible by
  $\Phi _{m_1}(q)\cdots\Phi _{m_i}(q)$ for all $i\ge 1$.  Since $\Phi _{m_i}(q)$
  are prime elements, it follows that $\Zqh$ is not a unique
  factorization domain.
\end{proof}

\subsubsection{Divisibility results}
\label{sec:divisibility-results}

\begin{proposition}
  \label{thm:44}
  For any integral homology sphere $M$, $J_M-1$ is divisible by
  $q^6-1$ in $\Zqh$.
\end{proposition}

\begin{proof}
  In Lemma \ref{r47}, we showed that $J_M-1$ is divisible by
  \begin{equation*}
    \frac{(q^2-1)(q^3-1)}{q-1}
    =\Phi _1(q)\Phi _2(q)\Phi _3(q)
    =\frac{q^6-1}{\Phi _6(q)}.
  \end{equation*}
  Kirby, Melvin and Zhang \cite{Kirby-Melvin-Zhang} proved that
  $\tau _6(M)=1$ for any integral homology sphere $M$.  Therefore, $J_M-1$
  is divisible by $\Phi _6(q)$, hence the assertion.
\end{proof}

Set
\begin{equation*}
  \tJ_M = q^{-\lambda _1(M)}J_M\in \Zqh,
\end{equation*}
where $\lambda _1(M)=6\lambda (M)\in 6\modZ $ is the first coefficient of the Ohtsuki
series $\tau ^O(M)\in \modZ [[q-1]]$.  Clearly,
$\iota _1(\tJ_M-1)\in (q-1)^2\modZ [[q-1]]$.  By Proposition \ref{thm:44},
$\tJ_M-1$ is divisible by $q^6-1=(q-1)(q+1)(q^2+q+1)(q^2-q+1)$.  Kirby
and Melvin \cite{Kirby-Melvin} proved that
$\tau _4(M)=\tau _{\sqrt{-1}}(M)=(-1)^{\chi (M)}\in \{\pm 1\}$, where
$\chi (M)\in \modZ /2\modZ $ is the Rochlin invariant of $M$, which is the mod $2$
reduction of $\lambda (M)$.  Therefore, we have $\mods _{\sqrt{-1}}(\tJ_M-1)=0$.
Hence, $\tJ_M-1$ is divisible by $q^2+1$.  These divisibility
properties of $\tJ_M-1$ imply the following.

\begin{proposition}
  \label{r10}
  For any integral homology sphere $M$, $\tJ_M-1$ is divisible by
  \begin{equation*}
    \Phi _1(q)^2\Phi _2(q)\Phi _3(q)\Phi _4(q)\Phi _6(q)
    =(q-1)(q^6-1)(q^2+1).
  \end{equation*}
\end{proposition}

\subsubsection{Applications to the Ohtsuki series}
\label{sec:appl-ohts-seri}

Here we apply Propositions \ref{thm:44} and \ref{r10} to the Ohtsuki
series $\tau ^O(M)$, and obtain some congruence relations on the
coefficients of $\tau ^O(M)$, which generalize well-known results by
Murakami \cite{Murakami1} and Ohtsuki \cite{Ohtsuki1}, and by Lin and
Wang \cite{Lin-Wang,Lin-Wang2} for the first and the second
coefficients.

We write
\begin{equation*}
  \tau ^O(M)=\iota _1(J_M)=1+\lambda _1(M)\modh +\lambda _2(M)\modh ^2+\cdots\in \modZ [[\modh ]],
\end{equation*}
where we set $\modh =q-1$.  We often write $\lambda _i=\lambda _i(M)$ in the
following.

The following gives an infinite set of congruence relations for the
coefficients $\lambda _i(M)$ of the Ohtsuki series of $M$.

\begin{proposition}
  \label{r68}
  For $k\ge 0$, we have
  \begin{equation}
    \label{e39}
    \sum_{i=0}^k a_i\lambda _{k-i+1}(M)\equiv0\pmod{\modZ }.
  \end{equation}
  Here $a_0=1/6$, $a_1=-1/4$, $a_2=17/72$, $a_3=-25/144$,
  $\dots\in \modZ [\frac16]$ are determined by
  \begin{equation*}
    \sum_{i\ge 0}a_i\modh ^i=\frac1{(q+1)(q^2+q+1)}=\frac{1}{6+9\modh +5\modh ^2+\modh ^3}.
  \end{equation*}
  In particular, each $\lambda _{k+1}(M)\mod6$, $k\ge 0$, is determined by
  $\lambda _i(M)$, $i=1,\ldots,k$.
\end{proposition}

\begin{proof}
  Set
  \begin{gather*}
    f(\modh )=\modh ^{-1}(\tau ^O(M)-1)=\sum_{k\ge 0}\lambda _{k+1}\modh ^k\in \modZ [[\modh ]],\\
    a(\modh )=(q+1)(q^2+q+1)=6+9\modh +5\modh ^2+\modh ^3.
  \end{gather*}
  Using Proposition \ref{thm:44}, we see that there is
  $g(\modh )=\sum_{k\ge 0}g_k\modh ^k\in \modZ [[\modh ]]$ such that $f(\modh )=a(\modh )g(\modh )$,
  hence $g(\modh )=f(\modh )a(\modh )^{-1}$.  (Here we used the divisibility of
  $J_M-1$ by $(q-1)(q+1)(q^2+q+1)$, not by $q^6-1$.)
  We have
  \begin{equation*}
    g_k= \sum_{i=0}^k a_i\lambda _{k-i+1}.
  \end{equation*}
  Since $g_k\in \modZ $ for all $k\ge 0$, we have \eqref{e39}.  Since
  $a_0=1/6$, the latter assertion follows.
\end{proof}

By \eqref{e39} and $a_0=\frac16$, we have
\begin{equation}
  \label{e71}
  \lambda _{k+1}\equiv - \sum_{i=1}^k 6a_i\lambda _{k-i+1}\pmod6.
\end{equation}
For small values of $k$, \eqref{e71} implies the following.
\begin{enumerate}
\item\label{item:5} $\lambda _1\equiv0\pmod6$.  This is known by Murakami
  \cite{Murakami1} and Ohtsuki \cite{Ohtsuki1}.
\item\label{item:6}
  $\lambda _2\equiv\frac{3}{2}\lambda _1\equiv\frac{1}{2}\lambda _1\pmod6$.  This is
  known by Lin and Wang \cite{Lin-Wang}.  See also Corollary \ref{r60}
  below.
\item\label{item:7} $\lambda _3\equiv\frac{3}{2}\lambda _2-\frac{17}{12}\lambda _1\pmod6$.
\item\label{item:8} $\lambda _4\equiv\frac{3}{2}\lambda _3-\frac{17}{12}\lambda _2+\frac{25}{24}\lambda _1\pmod6$.
\end{enumerate}

\begin{remark}
  \label{r69}
  One can use divisibility of $J_M-1$ by $q^6-1$ instead of
  $(q-1)(q+1)(q^2+q+1)=(q^6-1)/\Phi _6(q)$.  In this case, we obtain a
  result equivalent to Proposition \ref{r68} as a consequence of the
  fact that $\Phi _6(q)=q^2-q+1$ is invertible in $\modZ [[q-1]]$.
\end{remark}

The above-described idea can be applied also to Proposition \ref{r10}
to obtain better results.  We write
\begin{gather*}
  \iota _1(J_M-q^{\lambda _1(M)})=\tau ^O(M)-q^{\lambda _1(M)}=\sum_{k\ge 2}\lambda '_k(M)\modh ^k.
\end{gather*}
We have
\begin{equation}
  \label{e69}
  \lambda '_k(M) = \lambda _k(M)-\binom{\lambda _1(M)}{k}\in \modZ 
\end{equation}
for $k\ge 2$.

\begin{proposition}
  \label{r70}
  For $k\ge 0$, we have
  \begin{equation}
    \label{e68}
    \sum_{i=0}^k b_i\lambda '_{k-i+2}(M)\equiv0\pmod{\modZ }.
  \end{equation}
  Here $b_0=1/12$, $b_1=-5/24$, $b_2=41/144$, $b_3=-77/288$,
  $\dots\in \modZ [\frac16]$ are determined by
  \begin{equation*}
    \sum_{i\ge 0}b_i\modh ^i=\frac1{(q+1)(q^2+q+1)(q^2+1)}
    =\frac{1}{12+30\modh +34\modh ^2+21\modh ^3+7\modh ^4+\modh ^5}.
  \end{equation*}
  In particular, for $k\ge 0$, $\lambda '_{k+2}(M)\mod{12}$ is determined
  by $\lambda '_i(M)$, $i=2,\ldots,k+1$.
\end{proposition}

\begin{proof}
  The proof is similar to that for Proposition \ref{r68}, where we use
  Proposition \ref{r10} instead of Proposition \ref{thm:44}.  The
  details are left to the reader.
\end{proof}

\begin{corollary}
  \label{r59}
  For $k\ge 0$, $\lambda _{k+2}(M)\mod{12}$ is determined by $\lambda _i(M)$,
  $i=1,\ldots,k+1$.
\end{corollary}

\begin{proof}
  This immediately follows from Proposition \ref{r70} and \eqref{e69}.
\end{proof}

The case $k=0$ of \eqref{e68} is equivalent to the following.

\begin{corollary}[Lin and Wang \cite{Lin-Wang2}]
  \label{r60}
  We have
  \begin{equation*}
    \lambda _2(M)\equiv3\lambda (M)=\hf\lambda _1(M)\pmod{12}.
  \end{equation*}
\end{corollary}

\begin{proof}
  We have
  $\lambda _2(M)\equiv\binom{\lambda _1(M)}{2} =\binom{6\lambda (M)}{2}
  =18\lambda (M)^2-3\lambda (M) \equiv6\lambda (M)-3\lambda (M) \equiv3\lambda (M)\pmod{12}$.
\end{proof}

\begin{remark}
  \label{r61}
  Propositions \ref{thm:44} and \ref{r10} can be applied to the power
  series $\iota _{-1}(J_M)\in \modZ [[q+1]]$ and one can obtain results similar
  to Propositions \ref{r68} and \ref{r70}.  In particular, we have the
  following.
  \begin{enumerate}
  \item $\lambda _{-1,1}(M)\in 6\modZ $.
  \item $\lambda _{-1,2}(M)\equiv\hf\lambda _{-1,1}(M)\pmod{12}$.
  \item For $k\ge 2$, $\lambda _{-1,k}(M)\mod12$ is determined by
  $\lambda _{-1,i}(M)$, $i=1,2,\ldots,k-1$.
  \end{enumerate}
  Here $\lambda _{-1,k}(M)\in \modZ $ is defined as in \eqref{e75}.
\end{remark}

\begin{remark}
  \label{r63}
  As for the relation between $\lambda _1(M)$ and $\lambda _{-1,1}(M)$, we have
  \begin{equation}
    \label{e70}
    \lambda _1(M)\equiv\lambda _{-1,1}(M)\pmod{12}.
  \end{equation}
  Hence, like the Casson invariant, the mod $2$ reduction of the
  integer-valued invariant $\frac16\lambda _{-1,1}$ equals the Rochlin
  invariant of $M$.  \eqref{e70} follows from
  $\lambda _1(M)\equiv\lambda _{-1,1}(M)\equiv0\pmod3$ and
  $\lambda _1(M)\equiv\lambda _{-1,1}(M)\pmod{4}$.  (Here, the latter congruence
  follows from \eqref{e65} for $\zeta =1$, $\xi =-1$, $k=1$.)
\end{remark}

\subsubsection{The eighth WRT invariant $\tau _8(M)$}
\label{sec:wrt-invariant-at}

For each $r\in \modN $, it would be interesting to know the range of
$\tau _r(M)$.  Using the modified version
\begin{equation*}
  \tilde\tau _r(M)=\mods _{\zeta _r}(\tJ_M)=\zeta _r^{-6\lambda (M)}\tau _r(M)\in \modZ [\zeta _r],
\end{equation*}
one can expect to obtain sharper results than using $\tau _r(M)$.
Theorem \ref{r10} implies
\begin{equation}
  \label{e63}
  \tilde\tau _r(M)-1
  \in 
  (\zeta _r-1)(\zeta _r^6-1)(\zeta _r^2+1)\modZ [\zeta _r].
\end{equation}
If $r$ is of the form $ip^e$, where $i\in \{1,2,3,4,6\}$ and $p^e$ is a
prime power, \eqref{e63} gives some restriction on the value of
$\tilde\tau _r(M)$.

For example, we here consider $\tau _8(M)$.  We have $\zeta _8=\exp
\frac{2\pi \sqrt{-1}}{8}=\frac{1+\sqrt{-1}}{\sqrt{2}}$, and
\begin{equation*}
  \tilde\tau _8(M)=\sqrt{-1}^{\lambda (M)}\tau _8(M)\in \modZ [\zeta _8].
\end{equation*}

\begin{proposition}
  \label{r11}
  For any integral homology sphere $M$, we have
  \begin{equation*}
    \tilde\tau _8(M)-1
    \in 2(\zeta _8-1)\modZ [\zeta _8]
    =\Span_\modZ \{4,2\sqrt2,2+2\sqrt{-1},2+\sqrt2+\sqrt{-2}\}.
  \end{equation*}
\end{proposition}

\begin{proof}
  By \eqref{e63}, it follows that $\tilde\tau _8(M)-1$ is divisible in
  $\modZ [\zeta _8]$ by $(\zeta _8-1)(\zeta _8^4-1)=-2(\zeta _8-1)$.  Hence, we have the
  assertion.
\end{proof}

Note that $\modZ [\zeta _8]/2(\zeta _8-1)\modZ [\zeta _8]\cong(\modZ /4\modZ )\oplus(\modZ /2\modZ )^3$.
Computations suggest the following, a much stronger restriction on
$\tilde\tau _8$.

\begin{conjecture}
  \label{r7}
  For any integral homology sphere $M$, we have
  \begin{equation*}
    \tilde\tau _8(M)-1
    \in \Span_\modZ \{4,2\sqrt2\}.
  \end{equation*}
  (Consequently, $\tilde\tau _8(M)$ is real.  Hence, $\tau _8(M)$ is either
  real or purely imaginary.)
\end{conjecture}

\begin{remark}
  \label{r13}
  Conjecture \ref{r7} implies that the $\modQ $-vector subspace $G_8$ of
  $\modQ (\zeta _8)$ generated by $\tilde\tau _8(M)-1$ for all integral homology
  spheres $M$ is strictly smaller than $\modQ (\zeta _8)$.  More generally,
  let $G_r\subset \modQ (\zeta _r)$ be the $\modQ $-vector subspace generated by the
  $\tilde\tau _r(M)-1$.  Then $G_r\subsetneq\modQ (\zeta _r)$ if $r=1,2,3,4,6$.
  Conjecturally, $G_r=\modQ (\zeta _r)$ if $r\neq 1,2,3,4,6,8$.
\end{remark}

\section{The $p$-adic and the mod $p$ WRT functions}
\label{sec13}

In this section, we first observe that the invariant $J_M\in \Zqh$ can
be formally evaluated at any element in a commutative, unital ring,
and in particular at any complex number.  Then we introduce the
$p$-adic analytic versions and mod $p$ versions of the
WRT functions.

Most constructions in this section can be applied to any invariants of
links and $3$-manifolds which take values in $\Zqh$.

\subsection{Formal evaluations in rings}
\label{sec:form-eval-rings}

Let $R$ be a commutative, unital ring, and let $\alpha \in R$.  Let
$\mods _\alpha \zzzcolon \modZ [q]\rightarrow R$ be the ring homomorphism satisfying $\mods _\alpha (q)=\alpha $.
Then $\mods _\alpha $ induces a ring homomorphism
\begin{equation*}
  \hat{\mods }_\alpha \zzzcolon \Zqh\rightarrow \hat{R}^\alpha :=\varprojlim_{n} R/((\alpha )_n),
\end{equation*}
where
\begin{equation*}
  (\alpha )_n = (1-\alpha )(1-\alpha ^2)\cdots(1-\alpha ^n)\in R
\end{equation*}
for $n\ge 0$.

Using $\hat{\mods }_\alpha $, we can produce various ``reduced version''
$\hat{\mods }_\alpha (J_M)$ of the invariant $J_M\in \Zqh$, which may be regarded
as the ``value at $\alpha $'' of the WRT function of $M$.

\subsection{Evaluations at complex numbers}
\label{sec:expans-at-compl}

Let $\alpha \in \C$.  Then $\mods _\alpha \zzzcolon \modZ [q]\rightarrow \modZ [\alpha ]$ induces
\begin{equation*}
  \hat{\mods }_\alpha \zzzcolon \Zqh\rightarrow \Za=\varprojlim_{n} \modZ [\alpha ]/((\alpha )_n).
\end{equation*}
Thus, one can formally evaluate $J_M$ at any complex number.  We can
easily verify the following.
\begin{enumerate}
\item If $\alpha =0$, then $\Zx0=0$.
\item If $\alpha $ is transcendental, then $\hat{\mods }_\alpha \zzzcolon \Zqh\congto\Za$.
\item If $\alpha \in \calZ$, then $\Za\cong\modZ [\alpha ]$.
\end{enumerate}

The case $\alpha \in \modQ $ is of special interest.
\begin{proposition}
  \label{r12}
  If $a/b\in \modQ \setminus \{\pm 1\}$ with $a,b\in \modZ ,(a,b)=1$, then we have
  \begin{equation}
    \label{e58}
    \Zx{a/b}
    \cong \varprojlim_{m\in \modN ,(m,ab)=1}\modZ /m\modZ 
    \cong \prod_{p\in P_{ab}}\modZ _p,
  \end{equation}
  where, for $n\in \modN $, $P_n$ denote the set of primes $p$ coprime with
  $n$.
\end{proposition}

\begin{proof}
  We have
  \begin{equation*}
    \begin{split}
      \Zx{a/b}
      =\varprojlim_{n}\modZ [1/b]/((a/b)_n)
      \cong\varprojlim_{n}\modZ [1/b]/((b,a)_n),
    \end{split}
  \end{equation*}
  where we set
  \begin{equation*}
    (b,a)_n = b^{\hf n(n+1)}(a/b)_n =(b-a)(b^2-a^2)\cdots(b^n-a^n)\in \modZ .
  \end{equation*}
  The families of ideals in $\modZ [1/b]$, $X:=\{((b,a)_n)\}_{n\ge 0}$ and
  $Y:=\{(m)\}_{m\in \modN ,(m,ab)=1}$, are cofinal with each other.
  (Indeed, $a/b\neq1$ implies $X\subset Y$.  Conversely, for $m\in Y$, we
  have $b^n-a^n\equiv0\mod m$, where $n\ge 1$ is the least common
  multiple of the orders of $b$ and $a$ in $(\modZ /m\modZ )^\times $.  Then
  $(b,a)_n$ is divisible by $m$.)  Therefore,
  \begin{equation*}
    \begin{split}
      \Zx{a/b}
      \cong \varprojlim_{m\in \modN ,(m,ab)=1}\modZ [1/b]/(m)
      \cong \varprojlim_{m\in \modN ,(m,ab)=1}\modZ /m\modZ .
    \end{split}
  \end{equation*}
  Here the latter isomorphism follows from invertibility of $b$ in
  $\modZ /m\modZ $ for each $m$.  The latter isomorphism in \eqref{e58} is
  standard.
\end{proof}

By Proposition \ref{r12}, we may regard $\mods _{a/b}(J_M)$ as an element
in $\prod_{p\in P_{ab}}\modZ _p$.  Hence, for each $p\in P_{a,b}$, we can
extract a $\modZ _p$-valued function
\begin{equation}
  \label{e60}
  \tau ^p(M)\zzzcolon U_p(\modQ )\rightarrow \modZ _p,
\end{equation}
where $U_p(\modQ )$ is the set of $a/b\in \modQ $ with $a,b\in \modZ $, $(a,b)=1$,
such that $p$ is coprime with $ab$.

Note that $U_p(\modQ )=\modQ \cap U(\modZ _p)$, where $U(\modZ _p)$ is the group of units
in the ring $\modZ _p$.  We extend the domain of $\tau ^p(M)$ in the next
subsection.

\subsection{The $p$-adic WRT functions}
\label{sec:p-adic-wrt}

In this subsection, we show that the WRT
function of an integral homology sphere can be extended to a $p$-adic
analytic function on the unit circle in the field $\C_p$ of complex
$p$-adic numbers.

The field $\C_p$ is defined as follows (see \cite{Koblitz,Gouvea} for
details).  Let $\modQ _p$ be the field of $p$-adic numbers, and let
$|\cdot|_p\zzzcolon \modQ _p\rightarrow \R$ be the $p$-adic norm, defined by
\begin{equation*}
  |x|_p=\begin{cases}
  p^{-\ord_p(x)} \quad &\text{if $x\neq0$},\\
  0 &\text{if $x=0$},
  \end{cases}
\end{equation*}
where $\ord_p(x)\in \modZ $ is the $p$-adic ordinal.  Let $\bar\modQ _p$ denote
the algebraic closure of $\modQ _p$.  It is well known that the norm
$|\cdot|_p$ on $\modQ _p$ uniquely extends to $\bar\modQ _p$.  Then $\C_p$ is
defined to be the completion of $\bar\modQ _p$ with respect to the norm
$|\cdot|_p$.  The norm on $\C_p$ induced by $|\cdot|_p$ on $\bar\modQ _p$
is denoted again by $|\cdot|_p$.  It is well known that $\C_p$ is
algebraically closed.

Let $\modO _p$ and $\modP _p$ denote the valuation ring and the valuation
ideal, respectively, of $\C_p$, i.e.,
\begin{gather*}
  \modO _p=\{x\in \C_p\zzzcolon |x|_p\le 1\},\quad \modP _p=\{x\in \C_p\zzzcolon |x|_p<1\}.
\end{gather*}
Let
\begin{equation*}
  U(\modO _p)=\modO _p\setminus \modP _p=\{x\in \C_p\zzzcolon |x|_p=1\}
\end{equation*}
denote the unit circle of $\modO _p$, which consists of the multiplicative
units in $\modO _p$.  In particular, all the roots of unity in $\C_p$ are
contained in $U(\modO _p)$.  We may assume that the set
$\calZ(\subset \bar\modQ \subset \C)$ of complex roots of unity as a subset of
$U(\modO _p)$ (by choosing a field homomorphism
$\bar\modQ \hookrightarrow\bar\modQ _p(\subset \C_p$)).

\begin{lemma}
  \label{r46}
  For each $x\in U(\modO _p)$, the sequence $(x)_n$, $n\ge 0$, of elements in
  $\modO _p$ is converging to $0$ in $\modO _p$.
\end{lemma}

\begin{proof}
  Since $x\in U(\modO _p)$, there is $\zeta \in \calZ$ such that $x-\zeta \in \modP _p$.
  Set $r=\ord(\zeta )$.  Then
  \begin{equation*}
    \begin{split}
    1-x^r &= 1-(\zeta +(x-\zeta ))^r\\
    &=1-\sum_{i=0}^r\binom{r}{i}\zeta ^{r-i}(x-\zeta )^i\\
    &=-\sum_{i=1}^r\binom{r}{i}\zeta ^{r-i}(x-\zeta )^i\in \modP _p.
    \end{split}
  \end{equation*}
  Since $(x)_{rk}$ is divisible by $(1-x^r)^k$ for each $k\ge 0$, it
  follows that $(x)_n$ converges to $0$ as $n\rightarrow \infty $.
\end{proof}

Lemma \ref{r46} implies that for $x\in U(\modO _p)$, the ring homomorphism
$\modZ [q]\rightarrow \modO _p$, $f(q)\mapsto f(x)$, induces a continuous ring
homomorphism
\begin{equation*}
  \mods _x\zzzcolon \Zqh \rightarrow  \modO _p, \quad f(q)\mapsto f(x).
\end{equation*}
Gathering the $\mods _x$ for all $x\in U(\modO _p)$, we obtain a ring
homomorphism
\begin{equation*}
  \mods _{U(\modO _p)}\zzzcolon \Zqh \rightarrow  \Fun(U(\modO _p),\modO _p)
\end{equation*}
such that $\mods _{U(\modO _p)}(f)(x)=\mods _x(f(q))$.  For $f\in \Zqh$, the function
\begin{equation}
  \label{e79}
  \hat{f}=\mods _{U(\modO _p)}(f)\zzzcolon U(\modO _p)\rightarrow \modO _p
\end{equation}
is $p$-adic analytic of radius of convergence $\ge 1$ everywhere in
$U(\modO _p)$.  I.e., $\hat{f}$ has a converging power series expansion at
any point $x\in U(\modO _p)$ which coincides $\hat{f}$ for all $y\in U(\modO _p)$
with $|y-x|_p<1$.

It is easy to see that the restriction of $\hat{f}$ to
\begin{equation*}
  U(\modZ _p)=\{x\in \modZ _p\zzzcolon |x|_p=1\}=U(\modO _p)\cap \modQ _p
\end{equation*}
takes values in $\modZ _p$.

For any subset $V\subset U(\modO _p)$, $\mods _{U(\modO _p)}$ induces a ring homomorphism
\begin{equation*}
  \mods _V\zzzcolon \Zqh \rightarrow  \Fun(V,\modO _p).
\end{equation*}

\begin{proposition}
  \label{r54}
  Let $V\subset U(\modO _p)$ be an infinite subset.  Suppose that there is a
  point $x\in U(\modO _p)$ and a sequence $x_i\in V\setminus x$, $i=0,1,\ldots$, such that
  $\lim_{i\rightarrow \infty }x_i=x$ with respect to the $p$-adic topology.  Then
  $\mods _V$ is injective.
\end{proposition}

\begin{proof}
  Suppose $\mods _V(f)=0$, $f\in \Zqh$.  Define $\hat{f}$ as in \eqref{e79}.
  Set $D(x,1)=\{y\in \modO _p\zzzcolon |y-x|_p<1\}$.  The restriction $\hat{f}_x=\hat
  f|_{D(x,1)}\zzzcolon D(x,1)\rightarrow \modO _p$ can be expressed as a power series at
  $q=x$, convergent on $D(x,1)$.  The assumption implies that $\hat
  f_x=0$, i.e., $\hat{f}(x)=0$ for all $x\in D(x,1)$.

  Since $x\in U(\modO _p)$, we have $D(x,1)\cap \calZ\neq\emptyset$.  Moreover,
  if $\zeta \in D(x,1)\cap \calZ$ and $\xi \in \calZ$ with $\ord(\xi )$ a power of
  $p$, then $\zeta \xi \in D(x,1)\cap \calZ$.  Since for each $\zeta \in D(x,1)\cap \calZ$
  we have $f(\zeta )=\hat{f}(\zeta )=0$, it follows from Proposition \ref{r85} that
  $f=0$.
\end{proof}

Now we give applications of the above-mentioned facts to the
WRT invariants.  For each integral homology
sphere $M$, define the {\em $p$-adic WRT
function} of $M$ by
\begin{equation}
  \label{e77}
  \tau ^p(M)=\mods _{U(\modO _p)}(J_M)\zzzcolon U(\modO _p)\rightarrow \modO _p,
\end{equation}
which is $p$-adic analytic of radius of convergence $\ge 1$ everywhere
in $U(\modO _p)$.  $\tau ^p(M)$ restricts to
\begin{equation*}
  \tau ^p(M)\zzzcolon U(\modZ _p)\rightarrow \modZ _p.
\end{equation*}
Note that these invariants generalize the one given in \eqref{e60}.

If $V$ is as in Proposition \ref{r54}, then $J_M$ is determined by the
restriction of $\tau ^p(M)$ to $V$.  In particular, for $V=U(\modO _p)$,
$J_M$ (hence $\tau (M)$) is determined by $\tau ^p(M)\zzzcolon U(\modO _p)\rightarrow \modO _p$ and
also by its restriction to $U(\modZ _p)$.

The following generalizes Conjecture \ref{r73}.

\begin{conjecture}[$p$-adic Non-vanishing Conjecture]
  \label{r33}
  For any integral homology sphere $M$, for any prime $p$,
  $\tau ^p(M)\zzzcolon U(\modO _p)\rightarrow \modO _p$ is nowhere zero.
\end{conjecture}

If Conjecture \ref{r33} is true for at least one prime $p$, then
Conjecture \ref{r73} is true.

As a partial result of Conjecture \ref{r33}, we have the following
result, closely related to Corollary \ref{r71}.

\begin{proposition}
  \label{r50}
  If $\zeta \in \calZ\subset U(\modO _p)$ is a root of unity of order $1$, $2$, $3$, $4$ or
  $6$, then $\tau ^p(M)$ has no zero on $D(\zeta ,1)=\{x\in U(\modO _p)\zzzcolon |x-\zeta |_p<1\}$.
\end{proposition}

\begin{proof}
  The assertion follows from $\tau _\zeta (M)=\pm 1$ and the fact that the
  power series expansion of $\tau ^p(M)$ at $q=\zeta $ has coefficients in
  $\modO _p$.
\end{proof}

\begin{remark}
  \label{r53}
  By Corollary \ref{r71}, $\tau ^p(M)$ is non-vanishing at $q=\zeta \xi $,
  where $\xi \in \calZ$ is of any prime power order.  It follows from the
  analyticity that $\tau ^p(M)$ is non-vanishing on some neighborhood of
  $\zeta \xi $.
\end{remark}

\subsection{The mod $p$ WRT functions}
\label{sec:mod-p-wrt}

Let $p$ be a prime and let $\Fp$ be the field of $p$ elements.  Let
$\bFp$ denote the algebraic closure of $\Fp$, and set
$\bFp^\times =\Fp\setminus \{0\}$.

For $f(q)\in \Zqh$ and $x\in \bFp^\times $, the value $f(x)\in \bFp$ of $f(q)$ at
$x$ is well defined, since we have $(x)_m=0$ whenever $m\ge \ord(x)$.
(Recall that any $x\in \bFp^\times $ is a root of unity.)  Note that $f(x)$
is contained in the (finite) subfield of $\bFp$ generated by $x$ over
$\Fp$.  For $f(q)\in \Zqh$, define a function
\begin{equation*}
  \tilde{f}\zzzcolon \bFp^\times \rightarrow \bFp
\end{equation*}
by $\tilde{f}(x)=f(x)$.  Then $\tilde{f}$ is Galois equivariant, i.e.,
for any $\alpha \in \Gal(\bFp/\Fp)$, we have
\begin{equation*}
  f(\alpha (x))=\alpha (f(x)).
\end{equation*}
The correspondence $f(q)\mapsto \tilde{f}$ defines a ring homomorphism
\begin{equation*}
  \mods _{\bFp^\times }\zzzcolon \Zqh\rightarrow \Fung,
\end{equation*}
where $\Fung$ denotes the ring of all Galois-equivariant functions
from $\bFp^\times $ to $\bFp$.

Clearly, $\mods _{\bFp^\times }$ is not injective.  In fact, it factors
through the surjective homomorphism
\begin{equation*}
  \Zqh\rightarrow \Fpqh:=\varprojlim_{n}\Fpq/((q)_n) \bigl(\cong \Zqh\otimes _\modZ \Fp)
\end{equation*}
induced by $\modZ [q]\rightarrow \Fp[q]$.  Thus there is a ring homomorphism
\begin{equation*}
  s'_{\bFp^\times }\zzzcolon \Fpqh\rightarrow \Fung
\end{equation*}
induced by $\mods _{\bFp^\times }$.

\begin{proposition}
  \label{r35}
  $s'_{\bFp^\times }$ is surjective, but not injective.
  Moreover, we have an isomorphism
  \begin{equation}
    \label{e54}
    \Fung\cong \prod_{r\in \modN _p}\Fpq/(\Phi _r(q)).
  \end{equation}
\end{proposition}

\begin{proof}
  Set $\modN _p=\{r\in \modN \zzzvert \text{$r$ coprime with $p$}\}$.  We have the
  following isomorphisms
  \begin{equation*}
    \Fpqh
    \cong\varprojlim_{r\in \modN _p,k\ge 0}\Fpq/((q^r-1)^k)
    \cong\prod_{r\in \modN _p}\varprojlim_{k\ge 0}\Fpq/(\Phi _r(q)^k).
  \end{equation*}
  The first one follows, since in $\Fpq$ the families of ideals,
  $\{((q)_n)\}_{n\ge 0}$ and $\{((q^r-1)^k)\}_{r\in \modN _p,k\ge 0}$ are cofinal
  with each other.  (This is a consequence of the identity
  $q^{rp^e}-1=(q^r-1)^{p^e}$ in $\Fpq$ for $r\in \modN _p$ and $e\ge 1$.)  The
  second isomorphism can be verified using the Chinese Remainder
  Theorem in a way similar to the arguments in \cite[Section
  7.5]{H:cyclotomic}.

  For each $r\in \modN _p$, there is a surjective, non-injective
  homomorphism
  \begin{equation*}
    \varprojlim_{k\ge 0}\Fpq/(\Phi _r(q)^k)\rightarrow \Fpq/(\Phi _r(q)),\quad
  \end{equation*}
  induced by $\Fpq\rightarrow \Fpq/(\Phi _r(q))$.  Hence, there is a surjective,
  non-injective homomorphism
  \begin{equation*}
    \prod_{r\in \modN _p}\varprojlim_{k\ge 0}\Fpq/(\Phi _r(q)^k)\rightarrow \prod_{r\in \modN _p}\Fpq/(\Phi _r(q)).
  \end{equation*}

  Thus it remains to show the second assertion.  Note that $\Fung$ has
  a natural direct sum decomposition
  \begin{equation*}
    \Fung \cong \prod_{r\in \modN _p}\Fun_g(\calZ_r(\bFp),\bFp),
  \end{equation*}
  where $\calZ_r(\bFp)\subset \bFp$ is the set of primitive $r$th roots of
  unity, and $\Fun_g(\calZ_r(\bFp),\bFp)$ is the ring of
  Galois-equivariant functions from $\calZ_r(\bFp)$ into $\bFp$.  Since
  there is a natural $\Fp$-algebra isomorphism
  \begin{equation*}
    \Fp[q]/(\Phi _r(q))\cong \Fun_g(\calZ_r(\bFp),\bFp), \quad
    f(q)\mod{(\Phi _r(q))}\mapsto (x\mapsto f(x)),
  \end{equation*}
  we have  the isomorphism \eqref{e54}.
\end{proof}

Now, set
\begin{equation}
  \label{e22}
  \tmodp(M)=\mods _{\bFp^\times }(J_M)\in \Fung.
\end{equation}
We call it the {\em mod $p$ WRT function} of
$M$.  For $x\in \bFp^\times $,
\begin{equation*}
  \tau _x(M):=\tmodp(M)(x)
\end{equation*}
may be regarded as the WRT invariant at $x$.

$\tmodp(M)$ is determined by $J_M\mod p\in \Fpqh$, but Proposition
\ref{r35} implies that this fact gives no information on the range of
$\tmodp(M)$.

However, the $\tmodp(M)$ for infinitely many primes $p$ can
determine $J_M$ (and hence the WRT invariants)
by the following.

\begin{proposition}
  \label{r1}
  If $P$ is an infinite set of primes, then the ring homomorphism
  \begin{equation*}
    (\mods _{\bFp^\times })_{p\in P}
    \zzzcolon \Zqh\rightarrow \prod_{p\in P}\Fun(\bFp^\times ,\bFp)
  \end{equation*}
  is injective.
\end{proposition}

\begin{proof}
  Suppose $(\mods _{\bFp^\times })_{p\in P}(f(q))=0$, $f(q)\in \Zqh$.  It
  suffices to prove that for any $\zeta \in \calZ$, we have $\mods _\zeta (f(q))=0$.
  For each $p\in P$ coprime with $n:=\ord(\zeta )$, choose a ring
  homomorphism $\pi _{\zeta ,p}\zzzcolon \modZ [\zeta ]\rightarrow \bFp$, which maps $\zeta $ to a
  primitive $n$th root of unity in $\bFp$.  By assumption,
  $\pi _{\zeta ,p}\mods _\zeta (f(q))=0$.  This implies that $\mods _\zeta (f(q))\in \modZ [\zeta ]$ is
  divisible by $p$.  Since there are infinitely many such $p$, it
  follows that $\mods _\zeta (f(q))=0$.  Hence, the assertion follows from
  Proposition \ref{r85}.
\end{proof}

It is well known that $\bFp$ is isomorphic to the residue class
field $\modO _p/\modP _p$.  Thus, understanding the property of the function
$\tmodp(M)$ may be considered as the first step in understanding the
properties of $\tau ^p(M)$.  Let
\begin{equation*}
  \pi \zzzcolon \modO _p\rightarrow \modO _p/\modP _p\cong\bFp,
\end{equation*}
denote the projection.  Then we have
\begin{equation}
  \label{e84}
  \tau _{\pi (y)}(M)=\pi (\tau _y(M))
\end{equation}
for $y\in U(\modO _p)$.

The following conjecture is supported by some computer calculations.

\begin{conjecture}[mod $p$ Non-vanishing Conjecture]
  \label{r31}
  For any integral homology sphere $M$ and for any prime $p$,
  $\tmodp(M)\zzzcolon \bFp^\times \rightarrow \bFp$ is nowhere zero.
\end{conjecture}

Using \eqref{e84}, one can easily see that Conjecture \ref{r31}
implies Conjecture \ref{r33}.

\begin{remark}
  \label{r65}
  One can apply the construction in this section to any invariants of
  $3$-manifolds and links with values in $\Zqh$ to obtain functions
  like \eqref{e77} and \eqref{e22}.  For the case of the normalized
  Jones polynomials $\theta _2(J_T) = J_K(\V'_1)\in \modZ [q,q^{-1}]\subset \Zqh$ for
  $i\ge 1$ and the unified Kashaev invariant $\theta _0(J_T)\in \Zqh$, the
  analogues of Conjectures \ref{r33} and \ref{r31} do not hold.
  Probably, these conjecture do not hold in general for the normalized
  colored Jones polynomials $\theta _n(J_T) = J_K(\V'_{n-1})$ for $n\ge 3$.
\end{remark}

\section{Examples}
\label{sec14}

In this section, we compute the invariants of some links and integral
homology spheres.

\subsection{Borromean rings}
\label{sec:borromean-rings}

As in Section \ref{sec:invariant}, let $B\in \BT^0_3$ denote the
Borromean tangle.  The closure $A=\cl(B)$ is the Borromean rings.

\begin{proposition}
  \label{r4}
  For $i,i'\ge 0$, we have
  \begin{equation*}
    (1\otimes \trq^{P'_{i}}\otimes \trq^{P'_{i'}})(J_B) =\delta _{i,i'}(-1)^i\sigma _i .
  \end{equation*}
\end{proposition}

\begin{proof}
  Let $U_h^0\cong\modQ [H][[h]]$ denote the $h$-adic closure of the
  $\modQ [[h]]$-subalgebra of $U_h$ generated by $H$.  Recall that the
  {\em Harish-Chandra map} of $U_h$ is the continuous $\modQ [[h]]$-module
  homomorphism $\varphi\zzzcolon U_h \rightarrow  U_h^0$ determined by
  $\varphi(F^iH^jE^k)=\delta _{i,0}\delta _{k,0}H^j$ for $i,j,k\ge 0$.

  By \eqref{eq:7}, for any $x,x'\in \modR _{\modQ [[h]]}$, we have
  \begin{equation*}
    (\varphi\otimes \trq^x\otimes \trq^{x'})(J_B)
       =(\varphi\otimes \trq^x\otimes \trq^{x'})
      (\sum_{m_1,m_2,m_3,n_1,n_2,n_3\ge 0}
	B_{m_1,m_2,m_3,n_1,n_2,n_3}),
  \end{equation*}
  where $B_{m_1,m_2,m_3,n_1,n_2,n_3}$ denotes the summand in
  \eqref{eq:7}.  It follows from properties of $\varphi$ and $\trq^x$
  that $(\varphi\otimes \trq^x\otimes \trq^{x'})(B_{m_1,m_2,m_3,n_1,n_2,n_3})=0$
  unless $n_1=n_3=0$ and $m_1=m_3=m_2+n_2$.  Hence, setting $k=m_1$
  and $l=n_2$, we have
  \begin{equation*}
    \begin{split}
      &(\varphi\otimes \trq^x\otimes \trq^{x'})(J_B)\\
      &=\sum_{k\ge 0}\sum_{l=0}^k
      (\varphi\otimes \trq^x\otimes \trq^{x'})(B_{k,l,k,0,k-l,0})\\
      &=\sum_{k\ge 0}\sum_{l=0}^k
      (-1)^{l}	q^{ -\hf (k-l)(k-l+1)-l -4kl+2k^2 }\\
	&\quad
	\varphi(e^k   \tF^{(k)} K^{-2k+2l})
	\trq^x(e^{k-l}    \tF^{(k)} e^lK^{-2k})
	\trq^{x'}(\tF^{(l)} e^k    \tF^{(k-l)}K^{-2k} ).
    \end{split}
  \end{equation*}
  Using the identity $\trq^x(yy')=\trq^x(y'S^2(y))$ for $y,y\in U_h$, we
  obtain
  \begin{equation*}
    \begin{split}
      &\trq^x(e^{k-l}    \tF^{(k)} e^lK^{-2k})
      \trq^{x'}(\tF^{(l)} e^k    \tF^{(k-l)}K^{-2k} )\\
      &=q^{-2k+2l+3kl+l^2}\bbq{k}{l}
      \trq^x(\tF^{(k)}K^{-2k}e^k)
      \trq^{x'}(\tF^{(k)}K^{-2k}e^k).
    \end{split}
  \end{equation*}
  Now, set $x=P_{i}$ and $x'=P_{i'}$.  By Lemma \ref{lem:6}, we
  have for $i\ge 0$
  \begin{equation*}
      \trq^{P_i}(\tF^{(k)}K^{-2k}e^k)
      =\delta _{i,k}v^{k} q^{-k^2} \{k\}_q!.
  \end{equation*}
  Moreover, using \eqref{eq:60}, we obtain
  \begin{equation*}
    \varphi(e^k   \tF^{(k)} K^{-2k+2l})
    =\varphi(e^k   \tF^{(k)}) K^{-2k+2l}
    =\{H\}_{q,k}  K^{-2k+2l}.
  \end{equation*}
  Using the above formulas, we obtain
  \begin{equation*}
    \begin{split}
      &(\varphi\otimes \trq^{P_{i}}\otimes \trq^{P_{i'}})(J_B)
      =\delta _{i,i'}
	\Bigl(\sum_{l=0}^{i}
	(-1)^{l}
	q^{\hf l(l+3)-\hf i(i+3)}
	\bbq{i}{l}
	K^{2l-2i}
	\Bigr)
	\{H\}_{q,i}
	 (\{i\}_q!)^2.
    \end{split}
  \end{equation*}
  Using the identities
  \begin{gather*}
      \sum_{l=0}^{i}(-1)^{l}q^{\hf l(l+3)-\hf i(i+3)}
      \bbq{i}{l}K^{2l-2i}
      =\BBq{-H-2}{i},\\
      \varphi(\sigma _i)
      = \BB{H}{i} \BB{H+i+1}{i}
      = (-1)^i\BB{H}{i} \BB{-H-2}{i},
  \end{gather*}
  we obtain
  \begin{equation*}
    \begin{split}
      &(\varphi\otimes \trq^{P_{i}}\otimes \trq^{P_{i'}})(J_B)
      =\delta _{i,i'}
	\BBq{-H-2}{i}\BBq{H}{i}
	 (\{i\}_q!)^2\\
	&\quad =\delta _{i,i'} \BB{-H-2}{i}\BB{H}{i}
	 (\{i\}!)^2
	=\delta _{i,i'}(-1)^i\varphi(\sigma _i) (\{i\}!)^2.
    \end{split}
  \end{equation*}
  Hence, we have
  $(\varphi\otimes \trq^{P'_{i}}\otimes \trq^{P'_{i'}})(J_B)=\delta _{i,i'}(-1)^i\varphi(\sigma _i)$.
  Then the assertion follows from the well-known fact that $\varphi$
  is injective on the center $Z(U_h)$.
\end{proof}

Propositions \ref{r4} implies the following.

\begin{corollary}
  \label{r2}
  For $i,j,k\ge 0$, we have
  \begin{gather}
    \label{e3}
    J_{A}(P'_i, P'_j, P'_k)
    =\begin{cases}
    (-1)^i{\BB{2i+1}{i+1}}/{\{1\}} & \text{if $i=j=k$},\\
    0 & \text{otherwise},
    \end{cases}\\
    \label{e10}
    J_{A}(\V_i, \V_j, \V_k)
    =\sum_{p=0}^{\min(i,j,k)}(-1)^p
    \bb{i+1+p}{2p+1} \bb{j+1+p}{2p+1}\bb{k+1+p}{2p+1}
    (\{p\}!)^2\{2p+1\}_{2p}.
  \end{gather}
\end{corollary}

\begin{proof}
  For \eqref{e3}, use Proposition \ref{r22}.  For \eqref{e10}, use
  \eqref{e37}.
\end{proof}

Corollary \ref{r2} was first stated using the Kauffman bracket in
\cite{H:iias2000} without proof.  One can easily modify the arguments
in \cite{Masbaum} to obtain a skein-theoretic proof of Corollary
\ref{r2}.

\subsection{Powers of the ribbon element and the twist element}
\label{sec:powers-ribb-elem}
Here we give formulas for the powers of the ribbon element $\modr $ and
the twist element $\omega $.

For $p,n\ge 0$, set
\begin{equation*}
  S(n,p) = \{(i_1,\ldots,i_p)\zzzvert i_1,\ldots,i_p\ge 0,\quad i_1+\cdots+i_p=n\}.
\end{equation*}
For $\modi =(i_1,\ldots,i_p)\in S(n,p)$, set
$\bbq{n}{\modi }=\frac{[n]_q!}{[n_1]_q!\cdots[n_p]_q!}$,
$f(\modi )=\sum_{j=1}^{p-1}(s_j^2+s_j)$ with $s_j=\sum_{k=1}^ji_k$, and
$g(\modi )=\sum_{j=1}^p(p-j)i_j$.

\begin{proposition}
  \label{r8}
  For $p\ge 0$, we have
  \begin{gather*}
    \modr ^{-p}
    = v^{\frac p2 H(H+2)}\sum_{n\ge 0}v^{\hf n(n+3)}
    \Bigl(\sum_{\modi \in S(n,p)}\bbq{n}{\modi }q^{f(\modi )}K^{2g(\modi )}\Bigr)
    K^nF^{(n)}e^n,\\
    \modr ^p
    = v^{-\frac p2 H(H+2)}
    \sum_{n\ge 0}(-1)^nv^{-\hf n(n+3)}
    \Bigl(\sum_{\modi \in S(n,p)}\bb{n}{\modi }_{q^{-1}}q^{-f(\modi )}K^{-2g(\modi )}\Bigr)
    K^{-n}F^{(n)}e^n,
  \end{gather*}
  where
  $\bb{n}{\modi }_{q^{-1}}\Bigl(=q^{-\sum_{1\le j<k\le p}i_ji_k}\bbq{n}\modi \Bigr)$
  is the conjugate of $\bbq{n}{\modi }$.
\end{proposition}

\begin{proof}
  We only prove the first formula.  The other can be similarly proved.

  Define $a_{p,n}\in U_h^0$ by
  \begin{equation}
    \label{e14}
    \modr ^{-p}
    = v^{\frac p2 H(H+2)}\sum_{n\ge 0}v^{\hf n(n+3)}a_{p,n}
    K^nF^{(n)}e^n.
  \end{equation}
  We have $a_{0,n}=\delta _{n,0}$ and $a_{1,n}=1$.

  Since $\modr $ is central, we have
  \begin{equation*}
    \begin{split}
      \modr ^{-p-1}
      &=v^{\frac p2 H(H+2)}\sum_{i\ge 0}v^{\hf i(i+3)}a_{p,i}
      K^iF^{(i)}
      \Bigl(v^{\hf H(H+2)}\sum_{j\ge 0}v^{\hf j(j+3)}
      K^jF^{(j)}e^j\Bigr) e^i\\
      &=v^{\frac {p+1}2 H(H+2)}
      \sum_{k\ge 0}v^{\hf k(k+3)}
      \Bigl(
      \sum_{i=0}^k a_{p,i} \bbq{k}{i} q^{i^2+i} K^{2i}
      \Bigr)
      K^kF^{(k)}e^{k}.
    \end{split}
  \end{equation*}
  Hence, we have
  \begin{equation*}
    a_{p+1,k}=\sum_{i=0}^k a_{p,i} \bbq{k}{i} q^{i^2+i} K^{2i}.
  \end{equation*}
  By induction using this formula, we can verify for $p\ge 0$
  \begin{equation}
    \label{e16}
    a_{p,n}=\sum_{\modi \in S(n,p)}\bbq{n}{\modi }q^{f(\modi )}K^{2g(\modi )}.
  \end{equation}
\end{proof}

We can derive formulas for the powers of $\omega $ from Proposition
\ref{r8} as follows.

\begin{proposition}
  \label{r5}
  For $p\in \modZ $, we have
  \begin{equation}
    \label{e13}
    \omega ^p = \sum_{n\ge 0} \omega _{p,n}P'_n,
  \end{equation}
  where
  \begin{gather*}
    \omega _{p,n} =
    \begin{cases}
      v^{\hf n(n+3)}
      \sum_{\;\modi \in S(n,p)}\bbq{n}{\modi }q^{f(\modi )}&{\text{for $p\ge 0$}},\\
      (-1)^nv^{-\hf n(n+3)}
      \sum_{\;\modi \in S(n,-p)}\bb{n}{\modi }_{q^{-1}}q^{-f(\modi )}&{\text{for $p\le 0$.}}
    \end{cases}
  \end{gather*}
\end{proposition}

\begin{proof}
  Since the proof is similar to that of Proposition \ref{thm:38}, we
  give only a sketch proof.  We consider only the case $p\ge 0$; the
  other case follows by symmetry.

  We can express $\omega ^p$ as in \eqref{e13}, where the $\omega _{p,n}$ are in
  $\Zvv$.  For each $k\ge 0$, we have
  \begin{equation*}
    \begin{split}
      \langle \omega ^p,\V'_{2k}\rangle 
      =\sum_{n\ge 0}\omega _{p,n}\langle P'_n,\V'_{2k}\rangle 
      =\sum_{n=0}^k \omega _{p,n}\{k+n\}_{2n}/\{n\}!.
    \end{split}
  \end{equation*}
  By computing the action of $\modr ^{-p}$ on the weight $0$ vector
  $\modv ^{2k}_k\in \V_{2k}$ using Proposition \ref{r8}, we obtain
  \begin{equation*}
    \trq^{\V'_{2k}}(\modr ^{-p})
    =\sum_{n=0}^k v^{\hf n(n+3)}\epsilon (a_{p,n})\BB{k+n}{2n}/\{n\}!,
  \end{equation*}
  where $a_{p,n}$ is given by \eqref{e14} and \eqref{e16}.  Since
  $\langle \omega ^p,\V'_{2k}\rangle =\trq^{\V'_{2k}}(\modr ^{-p})$, we have
  \begin{equation*}
    \sum_{n=0}^k \omega _{p,n}\{k+n\}_{2n}/\{n\}!
    =\sum_{n=0}^k v^{\hf n(n+3)}\epsilon (a_{p,n})\BB{k+n}{2n}/\{n\}!
  \end{equation*}
  for all $k\ge 0$.  Since $\BB{k+n}{2n}/\{n\}!\neq0$ whenever
  $0\le n\le k$, it follows that $\omega _{p,n}=v^{\hf n(n+3)}\epsilon (a_{p,n})$
  for all $n\ge 0$, hence the assertion.
\end{proof}

A skein-theoretic version of Proposition \ref{r5} has been given by
Masbaum \cite{Masbaum}.

\subsection{Surgeries along the Borromean rings}
\label{sec:surg-borr-rings}

For $i,j,k\in \modZ $, we define $L_i$, $K_{i,j}$ and $M_{i,j,k}$ as
follows.  Let $A_1,A_2,A_3$ denote the components of the Borromean
rings $A$.
\begin{itemize}
\item $L_i=(A_1\cup A_2)_{A_3;-1/i}\subset S^3$, the two component link
 obtained from $A_1\cup A_2$ by surgery along $A_3$ with framing $-1/i$.
\item $K_{i,j}=(A_1)_{A_2\cup A_3;-1/i,-1/j}\subset S^3$, the knot obtained
from $A_1$ by surgery along $A_2$ and $A_3$ with framings $-1/i$ and
$-1/j$, respectively.
\item $M_{i,j,k}=(S^3)_{A_1\cup A_2\cup A_3;-1/i,-1/j,-1/k}$, the integral
homology sphere obtained from $S^3$ by surgery along $A_1$, $A_2$,
$A_3$ with framings $-1/i$, $-1/j$, $-1/k$, respectively.
\end{itemize}

\begin{proposition}
  \label{r3}
  For $i,j,k\in \modZ $ and $l,m\ge 0$, we have
  \begin{gather}
    \label{e5}
    J_{L_i}(P'_l,P'_m)
    =\delta _{l,m}\omega _{i,l}(-1)^l{\BB{2l+1}{l+1}}/{\{1\}},\\
    \label{e15}
    J_{K_{i,j}}(P'_l)
    =\omega _{i,l}\omega _{j,l}(-1)^l{\BB{2l+1}{l+1}}/{\{1\}},\\
    \label{e20}
    J_{M_{i,j,k}}=
    \sum_{l\ge 0}\omega _{i,l}\omega _{j,l}\omega _{k,l}(-1)^l{\BB{2l+1}{l+1}}/{\{1\}},\\
    \label{e18}
    J_{L_i}(\V_l,\V_m)
    =\sum_{s=0}^{\min(l,m)} \bb{l+s+1}{2s+1}\bb{m+s+1}{2s+1}\omega _{i,s}
    (-1)^s\{s\}!\{2s+1\}!/{\{1\}},
    \\
    \label{e19}
    \begin{split}
      J_{K_{i,j}}(\V_l)
      &=\sum_{s=0}^{l} \bb{l+s+1}{2s+1}\omega _{i,s}\omega _{j,s}
      (-1)^s\{2s+1\}!/{\{1\}}\\
      &=\sum_{s=0}^{l} (-1)^s\omega _{i,s}\omega _{j,s}
      \frac{\{l+s+1\}\{l+s\}\cdots\{l-s+1\}}{\{1\}}.
    \end{split}
  \end{gather}
\end{proposition}

\begin{proof}
  The first three identities follows from
  $J_{L_i}(P'_l,P'_m)=J_A(P'_l,P'_m,\omega ^i)$,
  $J_{K_{i,j}}(P'_l)=J_A(P'_l,\omega ^i,\omega ^j)$, and
  $J_{M_{i,j,k}}=J_A(\omega ^i,\omega ^j,\omega ^k)$.
  Use \eqref{eq:100} for the others.
\end{proof}

By \eqref{e15}, we have
\begin{equation*}
  J_{K_{i,j}}(P''_l) =(-1)^l\omega _{i,l}\omega _{j,l}.
\end{equation*}

The special case of Proposition \ref{r3} for $K_{i,j}$, $M_{i,j,k}$
with $i,j,k\in \{\pm 1\}$ appeared in \cite{H:iias2000,H:rims2001,Le2}
(without proofs).  The case for $K_{i,\pm 1}$, $i\in \modZ $, (and essentially
$L_{\pm 1}$) has been proved by Masbaum \cite{Masbaum} using skein
theory.

\section{Knots in integral homology spheres}
\label{sec15}

Let $K$ be a $0$-framed knot in an integral homology sphere $M$.  In
this section, we will define an invariant $J_{(M,K)}\in Z(\hUqv)$ of the
pair $(M,K)$.  Here $Z(\hUqv)$ denotes the $\Zqh$-subalgebra of
$Z(U_h)$ consisting of the infinite sums $\sum_{n\ge 0}a_n\sigma _n$ with
$a_n\in \Zqh$ for $n\ge 0$.  (The notation $Z(\hUqv)$ comes from the fact
that $Z(\hUqv)$ is equal to the center of the even part $\hUqv$ of the
completion $\hUq$ of $\Uq$ defined in \cite{H:uqsl2}.  We will not
need this fact in what follows.)  Clearly, we have
$Z(\tUqv)\subset Z(\hUqv)$.

In what follows, the invariant $J_{(M,K)}$ is constructed using a
generalization $J_{(M,K)}$ of the invariant $J_K(P''_m)$ of knots in
$S^3$ studied in Section~\ref{sec6}.  $J_{(M,K)}(P''_m)$ and its
generalization to algebraically-split links have been studied by
Garoufalidis and Le \cite{Garoufalidis-Le3}.  (In \cite[Section
4]{H:rims2001}, we announced a ``universal $sl_2$ invariant for links
in integral homology spheres''.  The invariant $J_{(M,K)}$ may be
considered as a special case of this.  We do not consider the link
case here.)

\subsection{Bottom knots with links}
\label{sec:bottom-knots-with}

Let $T=T_0\cup L_1\cup \cdots\cup L_l$ be a tangle in a cube consisting of a bottom
knot $T_0$ and an algebraically-split, $0$-framed link $L_1\cup \cdots\cup L_l$,
where $\lk(T_0,L_i)=0$ for $i=1,\ldots,l$.  For $x_1,\ldots,x_l\in \modP $, we
define $J_T(x_1,\ldots,x_l)$ as follows.  Let $T'=T_0\cup T_1\cup \cdots\cup T_l$
denote an $(l+1)$-component bottom tangle such that $T$ is obtained
from $T'$ by closing the components $T_1,\ldots,T_l$.  Set
\begin{equation*}
  J_T(x_1,\ldots,x_l) = (1\otimes \trq^{x_1}\otimes \cdots\otimes \trq^{x_l})(J_{T'}),
\end{equation*}
which is well-defined, i.e., does not depend on the choice of $T'$.

For $m\ge 0$, define $F_m(Z(\Uqv))$ as the $\modZ [q,q^{-1}]$-subalgebra of
$Z(\Uqv)$ consisting of the elements $\sum_{k\ge 0}a_k\sigma _k$, where
$a_k\in \modZ [q,q^{-1}]$ for $k\ge 0$ and if $0\le k\le m$ then $a_k$ is
divisible by $\frac{\{2m+1\}_{q,m+1}}{\{2k+1\}_{q,k+1}}$.  The
$F_m(Z(\Uqv))$ are a descending filtration of ideals in $Z(\Uqv)$.  We
have a natural isomorphism
\begin{equation*}
  \varprojlim_{m} Z(\Uqv)/F_{m}(Z(\Uqv)) \cong Z(\hUqv).
\end{equation*}
(Indeed, both are naturally isomorphic to $\varprojlim_{k,l}Z(\Uqv)/((q)_k,\sigma _l)$.)

\begin{proposition}
  \label{r76}
  Let $T$ be as above.  Let $m_1,\ldots,m_l\ge 0$.  Then we have
  \begin{equation*}
    J_{T}(\tP'_{m_1},\ldots,\tP'_{m_l})\in  F_m(Z(\Uqv)),
  \end{equation*}
  where $m=\max(m_1,\ldots,m_l)$.
\end{proposition}

\begin{proof}
  Let $g_k\in \modZ [q,q^{-1}]$, $k\ge 0$, be such that
  $\sum_{k\ge 0}g_k\sigma _k=J_T(\tP'_{m_1},\ldots,\tP'_{m_l})$.  By an argument
  similar to the proof of Theorem \ref{r18}, we have
  $g_k=\trq^{P''_k}(J_T(\tP'_{m_1},\ldots,\tP'_{m_l}))$.  Let $L$ denote
  the link obtained from $T$ by closing $K$.  Then we have
  \begin{equation*}
    g_k=\trq^{P''_k}(J_T(\tP'_{m_1},\ldots,\tP'_{m_l}))=J_L(P''_k,\tP'_{m_1},\ldots,\tP'_{m_l}).
  \end{equation*}
  It follows from Theorem \ref{r29} that if $0\le k\le m$ then $g_k$ is
  divisible by $\frac{\{2m+1\}_{q,m+1}}{\{2k+1\}_{q,k+1}}$.  Hence, we
  have the assertion.
\end{proof}

\begin{corollary}
  \label{r77}
  Let $T$ be as above.  Then $J_T\zzzcolon \modP \times \cdots\times \modP \rightarrow Z(\Uqv)$,
  $(x_1,\ldots,x_l)\mapsto J_T(x_1,\ldots,x_l)$, induces a well-defined
  $\modZ [q,q^{-1}]$-multilinear map
  \begin{equation*}
    J_T\zzzcolon \hP\times \cdots\times \hP\rightarrow Z(\hUqv).
  \end{equation*}
\end{corollary}

Let $L=\cl(T_0)\cup L_1\cup \cdots\cup L_l$ be the link obtained from $T$ closing
$T_0$.  Then we set
\begin{equation*}
  J_L(*,x_1,\ldots,x_l)=J_T(x_1,\ldots,x_l)\in Z(\hUqv)
\end{equation*}
for $x_1,\ldots,x_l\in \hP$.  Thus ``$*$'' means ``keep the element in this
place''.

\subsection{The invariant $J_{(M,K)}$}
\label{sec:invariant-jm-k}

Let $(M,K)$ be a pair of an integral homology sphere $M$ and a
$0$-framed knot $K\subset M$.  Choose an algebraically-split framed link
$L=K'\cup L_1\cup \cdots\cup L_l$ in $S^3$ such that
\begin{enumerate}
\item the framing of $L_i$ is $f_i\in \{\pm 1\}$,
\item $(S^3,K')_{L_1\cup \cdots\cup L_l}\cong (M,K)$.
\end{enumerate}
$L$ is called a {\em surgery presentation} for $(M,K)$.

\begin{theorem}
  \label{r79}
  Set
  \begin{equation}
    \label{e87}
    J_{(M,K)} = J_L(*,\omega ^{-f_1},\ldots,\omega ^{-f_l})\in Z(\hUqv).
  \end{equation}
  Then $J_{(M,K)}$ is a well-defined invariant of the pair $(M,K)$.
\end{theorem}

\begin{proof} (Sketch)
  The shortest way to prove the well-definedness of $J_{(M,K)}$ (i.e.,
  independence from the choice of $L$) may be to consider the elements
  $\mods _\zeta (J_{(M,K)}(\V'_k))$ for $k\ge 0$ and $\zeta \in \calZ$, where
  $\V'_k=\V_k/[k+1]$.  One can prove that this quantity is equal to
  the normalized WRT invariant of $(M,K)$ at
  $\zeta $ with $K$ colored by the $\V'_k$.  It follows from Proposition
  \ref{r85} that $J_{(M,K)}(\V'_k)$ is a well-defined invariant of
  $(M,K)$.  Therefore, $J_{(M,K)}(P''_k)$ for $k\ge 0$ are also
  well-defined, hence so is $J_{(M,K)}$.
\end{proof}

\begin{remark}
  \label{r78}
  Another way to prove Theorem \ref{r79} is to use a ``pair version''
  of Theorem \ref{r32}: {\it Two surgery presentations (in the above
  sense) of pairs of integral homology spheres and knots yield
  orientation-preserving homeomorphic results of surgery if and only
  if they are related by a sequence of stabilizations and Hoste moves.
  Here a Hoste move is defined to be a Fenn-Rourke move (see Figure
  \ref{fig:FR}) between two surgery presentations.}  (One can
  generalize this result to to links and tangles (in a homology ball).)

  In view of this result, one has only to show that the right-hand
  side of \eqref{e87} is invariant under a Hoste move.  To prove this
  invariance, one has to consider the ``evaluations'' at $\V_k$ for
  $k\ge 0$ as in the proof of Theorem \ref{r79}, but one does not have
  to consider the evaluations at roots of unity.
\end{remark}

\begin{remark}
  \label{r37}
  Yet another, conceptually more interesting, way to prove Theorem
  \ref{r79} is to extend the domain of definition of the invariant to
  ``bottom tangles in homology handlebodies''.  See Section
  \ref{sec:homology-cylinders} for some explanation.
\end{remark}

\subsection{Surgery along a knot with framing $1/m$}
\label{sec:surgery-along-knot}

Kirby and Melvin proved periodicity results for the WRT invariant of
integral homology spheres obtained from $S^3$ by surgery along a knot
with framings $1/m$ (Corollaries 8.15 and 8.26 of
\cite{Kirby-Melvin}).  We have the following generalization.

\begin{theorem}
  \label{r9}
  Let $K$ be a knot in an integral homology sphere $M$.  Let
  $M'=M_{K;1/m}$ be the result of surgery from $M$ along $K$ with
  framing $1/m$ with $m\in \modZ $.  Then we have
  \begin{equation*}
    J_{M'}\equiv J_M \pmod{(q^{2m}-1)}.
  \end{equation*}
\end{theorem}

\begin{proof}
  One way to prove the assertion is first to prove
  $\tau _\zeta (M')=\tau _\zeta (M)$ for all $\zeta \in \calZ$ with $\ord(\zeta )|2m$ and then
  to use Section \ref{sec:divis-ring-zqh}.  The former is a natural
  generalization of Kirby and Melvin's result, and is not difficult to
  prove.  Here we give a direct proof.

  We express $M$ as the result of surgery $S^3_{L}$ along an
  admissible framed link $L=L_1\cup \cdots\cup L_l$ in $S^3$, with framings
  $f_1,\ldots,f_l\in \{\pm 1\}$.  We may assume that $K$ corresponds to a knot
  $K'\subset S^3\setminus L$ such that $K'\cup L$ is algebraically split.  By
  puncturing $S^3$ at a point of $K'$, we obtain a bottom knot $T$ and
  a link $L$ in a cube, so that closing $T$ yields $K'\cup L$.  We have
  \begin{equation*}
    \begin{split}
      J_{M'}-J_M
      &= J_{K\cup L^0}(\omega ^{-m},\omega ^{-f_1},\ldots,\omega ^{-f_l})
      -J_{K\cup L^0}(1,\omega ^{-f_1},\ldots,\omega ^{-f_l})\\
      &= J_{K\cup L^0}(\omega ^{-m}-1,\omega ^{-f_1},\ldots,\omega ^{-f_l})\\
      &= \trq^{\omega ^{-m}-1}(J_{T\cup L^0}(\omega ^{-f_1},\ldots,\omega ^{-f_l})),
    \end{split}
  \end{equation*}
  where $L^0$ is $L$ with $0$ framings.  In view of the last
  subsection, we have
  \begin{equation*}
    J_{T\cup L^0}(\omega ^{-f_1},\ldots,\omega ^{-f_l})=\sum_{k\ge 0}a_k\sigma _k,
  \end{equation*}
  where $a_k\in \Zqh$ for $k\ge 0$.  It suffices to prove that
  $\trq^{\omega ^{-m}-1}(\sigma _k)\in (q^{2n}-1)\Zqh$ for each $k\ge 0$.  We have
  \begin{equation*}
    \trq^{\omega ^{-m}-1}(\sigma _k)=\langle \omega ^{-m}-1,S_k\rangle 
    =\langle \omega ^{-m}-1,\sum_{i=0}^{k}b_{k,i}\V_{2i}\rangle ,
  \end{equation*}
  where $b_{k,i}\in \Zqq$.  We have
  \begin{equation*}
    \langle \omega ^{-m}-1,\V_{2i}\rangle =(q^{-i(i+1)m}-1)[2i+1]\in (q^{2m}-1)\Zqq.
  \end{equation*}
  Hence, the assertion follows.
\end{proof}

\subsection{Cyclotomic finite type invariants of integral homology spheres}
\label{sec:cycl-finite-type}

Kricker and Spence \cite{Kricker-Spence} (see also \cite{Ohtsuki4})
proved that for integral homology spheres the Ohtsuki series $\tau ^O(M)$
modulo $(q-1)^{k+1}$ is a finite type invariant of degree $3k$ in the
sense of Ohtsuki \cite{Ohtsuki6}.  It would be natural to ask the
topological meaning of $J_M\mod \Phi _d(q)$ for $d\in \Funz$, where $\Funz$
denote the set of functions from $\modN =\{1,2,\ldots\}$ into
$\modZ _+=\{0,1,\ldots\}$ vanishing for all but finitely many elements of
$\modN $.  $J_M\mod \Phi _d(q)$ is not an Ohtsuki finite type invariant in
general.  Thus it would be interesting to have a variant of the notion
of Ohtsuki finite type invariants in which $J_M\mod \Phi _d(q)$ is of
``finite type''.  In what follows, we give a conjectural definition of
such a notion of finite type.

Let $\modM $ denote the set of orientation-preserving homeomorphism
classes of integral homology spheres.  Let $\modZ \modM $ denote the free
abelian group generated by $\modM $, which has a standard ring structure
with multiplication induced by connected sum.  The {\em Ohtsuki
filtration}
\begin{equation*}
  \modZ \modM =F_0\supset F_1\supset F_2\supset \cdots
\end{equation*}
is defined as follows.  For $k\ge 0$, let $F_k$ denote the
$\modZ $-submodule of $\modZ \modM $ spanned by the alternating sums
\begin{equation*}
  [M;L_1,\ldots,L_k]=\sum_{L'\subset \{L_1,\ldots,L_k\}}(-1)^{|L'|}M_{L'}
\end{equation*}
for all pairs of $M\in \modM $ and algebraically-split, $\pm 1$-framed links
$L=\{L_1,\ldots,L_k\}$ in $M$.  Here the sum runs over all the sublinks $L'$
of $L$, $|L'|$ denotes the number of components of $L'$, and $M_{L'}$
denotes the result from $M$ of surgery along $L'$.  A $\modZ $-linear map
$f\zzzcolon \modZ \modM \rightarrow A$ with $A$ an abelian group is called an {\em Ohtsuki
invariant of type $k$} if $f|_{F_{k+1}}=0$.

Now we generalize the notion of finite type.  For $d\in \Funz$, let $F_d$
denote the $\modZ $-submodule of $\modZ \modM $ generated by the alternating sums
$[M;L_1,\ldots,L_k]$, where $M\in \modM $ and $L=\{L_1,\ldots,L_k\}$ is an
algebraically-split framed link in $M$ with framings in
$\{1/m\zzzvert m\in \modZ ,m\neq0\}$ such that for each $n\in \modN $ there is exactly
$d(n)$ components of $L$ of framing $\pm n$.  One sees easily that the
$F_d$, $d\in \Funz$, form an inverse system of $\modZ $-submodules of $\modZ \modM $.
Note that for $k\ge 0$, $F_k$ is equal to $F_{d_k}$, where
$d_k(n)=k\delta _{n,1}$.

A $\modZ $-linear map $f\zzzcolon \modZ \modM \rightarrow A$ with $A$ an abelian group is said to be
of {\em cyclotomic finite type} if $f|F_d=0$ for some $d$.

Theorem \ref{r9} implies that for all $m\in \modN $, $J_M\mod (q^m-1)$, is
of cyclotomic finite type.  It follows that for all $\zeta \in \calZ$, the
WRT invariant $\tau _\zeta (M)$ is of cyclotomic finite type.

\begin{conjecture}
  \label{r75}
  For each $d\in \Funz$, there is $d'\in \Funz$ such that the ring
  homomorphism
  \begin{equation*}
    J\zzzcolon \modZ \modM \rightarrow \Zqh,\quad M\mapsto J_M,
  \end{equation*}
  maps $F_{d'}$ into $\Phi _d(q)\Zqh$.  Thus $J_M\mod(\Phi _d(q))$ is of
  cyclotomic finite type.  Consequently, $J$ induces a continuous ring
  homomorphism
  \begin{equation*}
    J\zzzcolon \widehat{\modZ \modM } \rightarrow  \Zqh,
  \end{equation*}
  where
  \begin{equation*}
    \widehat{\modZ \modM }= \varprojlim_{d\in \Funz}\modZ \modM /F_d.
  \end{equation*}
\end{conjecture}

If $M$ is an integral homology sphere and if $d\in \Funz$, then there is
an integral homology sphere $M'$ not homeomorphic to $M$, such that
$M-M'\in F_d$.  Indeed, if $L$ is an algebraically-split, Brunnian
framed link in $M$ with framings in $\{1/m\zzzvert m\in \modZ ,m\neq0\}$ such that
for each $n\in \modN $ there are at least $d(n)$ components of $L$ of
framings $\pm n$, then one can show that $M_L-M\in F_d$ by modifying a
well-known arguments in the study of Ohtsuki finite type invariants.
For a certain type of Brunnian links (e.g., iterated Bing doubles of
Borromean rings), one can show that $M_L$ and $M$ are not
homeomorphic, using the Le-Murakami-Ohtsuki invariant \cite{LMO}.

Recall that for integral homology spheres $M$ the Le-Murakami-Ohtsuki
invariant $Z^{\LMO}(M)$ is a ``universal finite type invariant''
\cite{Le6}, i.e., every Ohtsuki finite type invariant factors through
$Z^{\LMO}(M)$.  It would be interesting to have a ``cyclotomic finite
type'' version $Z^{\operatorname{cyc}}(M)$ of the Le-Murakami-Ohtsuki
invariant which recovers {\em via ring homomorphisms} the unified
$sl_2$ invariant $J_M$ (as well as the invariant for the other simple
Lie algebras, see Section \ref{sec:general} below) and the
Le-Murakami-Ohtsuki invariant $Z^{\LMO}(M)$.  As we mentioned in the
introduction, $Z^{\LMO}(M)$ determines $J_M$.  However, this
determination is not via a ring homomorphism, i.e., there is no
natural homomorphism from the ring $\mathcal A(\emptyset)$, in which the
Le-Murakami-Ohtsuki invariant takes values, into $\Zqh$.

\begin{remark}
  \label{r38}
  Similarly to the Cochran and Melvin's generalization
  \cite{Cochran-Melvin:finitetype} of the Ohtsuki finite type
  invariants, one can generalize the above definition of ``cyclotomic
  finite type'' to oriented, connected $3$-manifolds by using
  surgeries along null-homologous, algebraically-split framed links
  with framings in $\{1/m\zzzvert m\in \modZ ,m\neq0\}$.  It would be natural to
  expect that the generalization of $J_M$ to rational homology spheres
  defined in \cite{Beliakova-Blanchet-Le,Le4} (see Section
  \ref{sec:rati-homol-spher} below) modulo $(\Phi _d(q))$, $d\in \Funz$,
  would be of ``cyclotomic finite type'' in the above-explained sense.
\end{remark}

\section{Concluding remarks}
\label{sec16}

In this section, we give some remarks and discussions.

\subsection{The WRT invariants as limiting values of holomorphic
  functions}
\label{sec:wrt-invariants-as}
As explained in the introduction, the invariant $J_M$ for an integral homology
sphere $M$ unifies the WRT invariants $\tau _\zeta (M)$ for all the roots of
unity in an algebraic way.  It should be remarked that there is
another {\em analytic} approach to ``unify'' the WRT invariants by
realizing the radial limiting values of a holomorphic function defined
on the unit disk $\{z\in \C;\;|z|<1\}$, studied in
\cite{Lawrence:integrality,Lawrence:holomorphic,Lawrence-Rozansky,Lawrence-Zagier,Hikami1,Hikami3,Hikami4,Hikami5}.

Unfortunately, the invariant $J_M\in \Zqh$ does not immediately
determine such a holomorphic function, i.e., there is no natural
homomorphism from the ring $\Zqh$ to the ring of holomorphic function
on $|q|<1$.  However, some expressions for $J_M$ may define
well-defined holomorphic functions for on $|q|<1$ as we explain below
for the Poincar\'e homology sphere $\Sigma (2,3,5)$.

We have
\begin{equation}
  \label{e11}
  J_{\Sigma (2,3,5)}
  =\sum_{n\ge 0}q^n\frac{(1-q^{n+1})(1-q^{n+2})\cdots(1-q^{2n+1})}{1-q}\in \Zqh.
\end{equation}
The modified WRT invariant $W(\zeta )=\zeta (\zeta -1)\tau _\zeta (\Sigma (2,3,5))$,
$\zeta \in \calZ$, defined by Lawrence and Zagier \cite{Lawrence-Zagier}, is
unified into
\begin{equation*}
  W(q)=q(q-1)J_{\Sigma (2,3,5)}\in \Zqh.
\end{equation*}
For our purpose, it is useful to modify it further as
\begin{equation}
  \label{e7}
  \begin{split}
    B(q)&=1-W(q) = 1+q(1-q)J_{\Sigma (2,3,5)}\\
    &=\sum_{n\ge 0}q^n(q^n;q)_n,
  \end{split}
\end{equation}
where $(q^n;q)_n=(1-q^n)(1-q^{n+1})\cdots (1-q^{2n-1})$.  A beautiful
observation by Lawrence and Zagier \cite{Lawrence-Zagier} is that the
series
\begin{equation*}
  A(q)=\sum_{n\ge 1}\chi _+(n)q^{\frac{n^2-1}{120}}=1+q+q^3+q^7-q^8-\dots \in \modZ [[q]],
\end{equation*}
which converges on $|q|<1$, converges radially at each root of
unity $\zeta $ to $2(1-W(\zeta ))=2B(\zeta )$.  They also showed that the radial
asymptotic expansion of $A(q)$ at $q=1$ gives $2$ times the Ohtsuki
series $\tau ^O(M)$.

Note that the formula \eqref{e7} for $B(q)$ can also define a power
series in $q$.  (Actually, \eqref{e7} defines an element of
$\varprojlim_{n}\modZ [q]/(q^n(q)_n)\cong\Zqh\times \modZ [[q]]$.)  Let
$B(q)_{\modZ [[q]]}\in \modZ [[q]]$ denote this element.  Hikami \cite{Hikami4} proved
\begin{equation*}
  A(q)=B(q)_{\modZ [[q]]}\in \modZ [[q]].
\end{equation*}
Thus, there are two ways to obtain the value of the (modified) WRT
invariants from $B(q)$.  One is just to evaluate at roots of unity,
and the other is to expand it in $q$, take the radial limit at roots
of unity and then multiply by $\half$.  The meaning of the factor
$\half$ is not clear to the author.  It would be natural to expect
that the radial asymptotic expansion at each root of unity $\zeta $ is $2$
times the power series $\iota _\zeta (B(q))\in \modZ [\zeta ][[q-\zeta ]]$.

It would be natural to expect that certain formulas \eqref{e11} for
$J_M$ of integral homology spheres $M$ may give a power series in $q$
whose radial asymptotic behavior at a root of unity $\zeta $ is closely
related to the ``algebraic behavior near $\zeta $'' obtained by the map
$\iota _\zeta \zzzcolon \Zqh\rightarrow \modZ [\zeta ][[q-\zeta ]]$.  For example, for an integral homology
sphere $M_{i,j,k}$ with $i,j,k>0$ defined in Section
\ref{sec:surg-borr-rings}, the formula \ref{e20} of $J_{M_{i,j,k}}$
defines an element of $\modZ [[q]]$ which converges on $|q|<1$.  Note
that the Poincar\'e homology sphere $\Sigma (2,3,5)$ is the special case
$M_{1,1,1}$.

\subsection{Generalizations to simple Lie algebra}
\label{sec:general}

In a joint work with Le \cite{H-Le:in-preparation}, for each simple
Lie algebra $\g$, we will generalize the invariant $J_M\in \Zqh$ of
integral homology sphere $M$ to an invariant $J^\g_M\in \Zqh$.  This
invariant can be characterized by the property that for each root of
unity $\zeta $ such that the quantum $\g$ invariant $\tau ^\g_M$ is
well-defined, we have $\mods _\zeta (J^\g_M)=\tau ^\g_\zeta $.  This specialization
property holds also for the projective quantum $\g$ invariant.  See
\cite[Conjecture 7.29]{Ohtsuki5} for a precise statement.

\subsection{Rational homology spheres}
\label{sec:rati-homol-spher}

The invariant $J_M$ has been generalized by Beliakova, Blanchet and Le
\cite{Beliakova-Blanchet-Le} to closed $3$-manifold whose first
homology group is isomorphic to $\modZ /2\modZ $, and by Le \cite{Le4} to
rational homology spheres.  These invariants take values in
modifications of the ring $\Zqh$.  We briefly describe these
invariants below.  See the original papers for the details.

In \cite{Beliakova-Blanchet-Le}, it is proved that there is an
invariant of closed $3$-manifolds $M$ with $H_1M\cong\modZ /2\modZ $ with
values in the completion
\begin{equation*}
  \widehat{\modZ [v]}_2
  :=\varprojlim_n\modZ [q^{1/2}]/\Bigl(\prod_{i=1}^n(1+(-q^{1/2})^i)\Bigr),
\end{equation*}
which specializes to the a normalized version of the $SO(3)$ quantum
invariant for all $q^{1/2}$ a root of unity of order $\not\equiv2(4)$.
(In the notation of \cite{H:cyclotomic}, we have
$\widehat{\modZ [v]}_2=\modZ [q^{1/2}]^{\{n\in \modN \zzzvert n\not\equiv2(4)\}}$.)  A
generalization to closed spin manifolds with $H_1\cong\modZ /2\modZ $ is also
given.

In \cite{Le4}, it is proved that for a rational homology sphere $M$
with $\max\{\ord(g)\zzzvert g\in H_1M\}=d$, there is an invariant $I_M$ of $M$
with values in
\begin{equation*}
  \hat\Lambda _d
  :=\varprojlim_{g\in \operatorname{Fun}^0(\modN _d,\modZ _+)}\modZ [1/d][q^{1/d}]/
  \Bigl(\prod_{n\in \modN _d}\bigl(\Phi _n(q^{1/d})^{g(n)}\bigr)\Bigr).
\end{equation*}
where $\modN _d=\{n\in \modN \zzzcolon \text{$n$ coprime with $d$}\}$, and
$\operatorname{Fun}^0(\modN _d,\modZ _+)$ denotes the set of functions from
$\modN _d$ to $\modZ _+=\{0,1,\ldots \}$ vanishing for all but finitely many
elements of $\modN _d$.  The invariant $I_M$ specializes (essentially) to
the $SO(3)$ quantum invariant for $q$ any root of unity of order odd
and coprime with $d$.

For $M$ with $H_1(M;\modZ )\cong\modZ /2\modZ $, the Beliakova-Blanchet-Le
invariant is better than Le's invariant in the sense that the former
is defined in a smaller subring, i.e., has stronger integrality.  It
would be interesting to generalize these two invariants to an
invariant of all rational homology spheres which is defined in smaller
rings than the $\hat\Delta _d$.

\subsection{Bottom tangles in homology handlebodies}
\label{sec:homology-cylinders}

An interesting generalization of integral homology spheres (more
precisely, integral homology $3$-balls) are {\em homology cylinders}
over a surface \cite{Goussarov,H:claspers} (see also
\cite{Garoufalidis-Levine,Habegger}).  (In \cite{H:claspers} they are
called ``homologically-trivial homology cobordisms''.)  For $g\ge 0$,
let ${\mathcal H}_{g,1}$ denote the monoid of homology cylinders over
a surface $\Sigma _{g,1}$ of genus $1$ with one boundary component.  For
$g=0$, ${\mathcal H}_{0,1}$ is identified with the monoid of integral
homology spheres.  The Torelli group $\modI _{g,1}$ for $\Sigma _{g,1}$ is
regarded as a subgroup of ${\mathcal H}_{g,1}$ via the mapping
cylinder construction.

In \cite[Section 14.5]{H:universal}, we defined the ``category of
bottom tangles in handlebodies'' $\modB $ and the ``category of bottom
tangles in homology handlebodies'' $\bar\modB $.  Here $\modB $ is a
subcategory of $\bar\modB $, and $\bar\modB $ may be regarded as a
subcategory of the category $\mathcal C$ of cobordisms of surfaces with
connected boundaries as introduced by Crane and Yetter
\cite{Crane-Yetter} and by Kerler \cite{Kerler:99}.  These categories
are braided categories and generated by a braided Hopf algebra and
some additional morphisms.  Homology cylinders are contained in
$\bar\modB $ as morphisms.

In a future paper, using a completion of the algebra $\Uqv$, we will
construct a braided functor $\bar J$ from $\bar\modB $ to a certain
braided category defined over $\Zqh$.  This functor maps each integral
homology $3$-ball $M'$ to the multiplication map $\Zqh\rightarrow \Zqh$,
$x\mapsto xJ_M$, where $M$ is the integral homology sphere obtained
from $M'$ by capping off the boundary with a ball.  Thus, $\bar J$ may
be regarded as a generalization of $J_M$ into a large class of
$3$-manifolds with boundary, which contains all homology cylinders.
This functor restricts to a representation of the Torelli group of
each~$\Sigma _{g,1}$.  It is expected that the functor $\bar J$ would lead
to some refinements of the integral structures in topological quantum
field theories \cite{Gilmer,GMW,GM}.

\begin{acknowledgments}
  The author would like to thank Thang Le for numerous, helpful
  conversations and correspondences.  The author would also like to
  thank Gregor Masbaum, Hitoshi Murakami, and Tomotada Ohtsuki for
  helpful conversations.
\end{acknowledgments}

\end{document}